\documentclass[12pt,a4paper,twoside]{amsart}
 \usepackage{amssymb,latexsym,amsmath,graphicx,epsf,color}
\usepackage[table]{xcolor}

 \numberwithin{equation}{section}

\newcommand{\argmin}{\mathop{\mathrm{arg\,min}}}

\def\weakto{\ \longrightarrow^{\hspace{-.8cm}\hbox{\tiny BVw}}}

\DefineNamedColor{named}{VeryLightOrange}{rgb}{1,0.8824,0.3333}

\newcommand{\mnewtext}{}

\newtheorem{theorem}{Theorem}[section]
\newtheorem{definition}[theorem]{Definition}
\newtheorem{lemma}[theorem]{Lemma}

\def \uMAP {\widetilde{u}}
\def \e {\varepsilon}

\def \a {\alpha}
\def \ba {\begin {eqnarray*}}
\def \ea {\end  {eqnarray*}}

\def \beq {\begin {eqnarray}}
\def \eeq {\end {eqnarray}}
\newcommand{\R}{\mathbb{R}}
\newcommand{\ra}{\rightarrow}

\newcommand{\f}{\mathbf{f}}
\newcommand{\g}{\mathbf{g}}

\newcommand{\hh}{\mathbf{h}}
\newcommand{\vv}{\mathbf{v}}
\newcommand{\zz}{\mathbf{z}}
\newcommand{\bb}{\mathbf{b}}
\newcommand{\yy}{\mathbf{y}}
\newcommand{\cc}{\mathbf{c}}
\newcommand{\atten}{u}
\newcommand{\dd}{d}
\newcommand{\xx}{\mathbf{x}}
\newcommand{\pp}{\mathbf{p}}

\newcommand{\greal}{\widetilde{\mathbf{g}}}
\newcommand{\Areal}{\mathcal{A}}

\newcommand{\fan}{\widehat{\mathbf{f}}_{\alpha,n}}

\title[TV parameter choice]{Multi-resolution parameter choice method\\for total variation regularized tomography}

\author[Niinim\"aki {\em et al.}]{Kati Niinim\"aki \and Matti Lassas \and Keijo H\"am\"al\"ainen \and Aki Kallonen \and Ville Kolehmainen \and Esa Niemi \and Samuli Siltanen }

\begin{document}

\maketitle
\centerline{\today}

\begin{abstract}
A computational method is introduced for choosing the regularization parameter for total variation (TV) regularization. The approach is based on computing reconstructions at a few different resolutions and various values of regularization parameter. The chosen parameter is the smallest one resulting in approximately discretization-invariant TV norms of the reconstructions. The method is tested with X-ray tomography data measured from a walnut and compared to the S-curve method. The proposed method seems to automatically adapt to the desired resolution and noise level, and it yields useful results in the tests. The results are comparable to those of the S-curve method; however, the S-curve method needs {\em a priori} information about the sparsity of the unknown, while the proposed method does not need any {\em a priori} information (apart from the choice of a desired resolution). Mathematical analysis is presented for (partial) understanding of the properties of the proposed parameter choice method. It is rigorously proven that the TV norms of the reconstructions converge with any choice of regularization parameter. 
\end{abstract}

\clearpage

\tableofcontents

\section{Introduction}

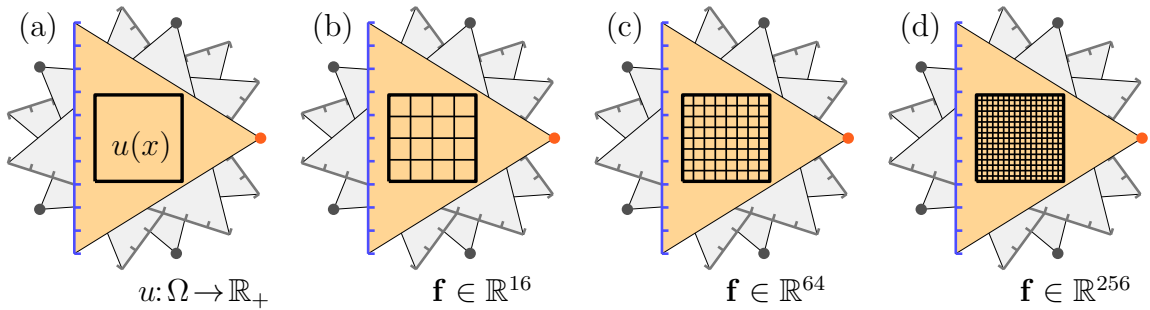
\begin{figure}
\centering
\begin{picture}(390,120)
\put(-20,10){\includegraphics[width=110pt]{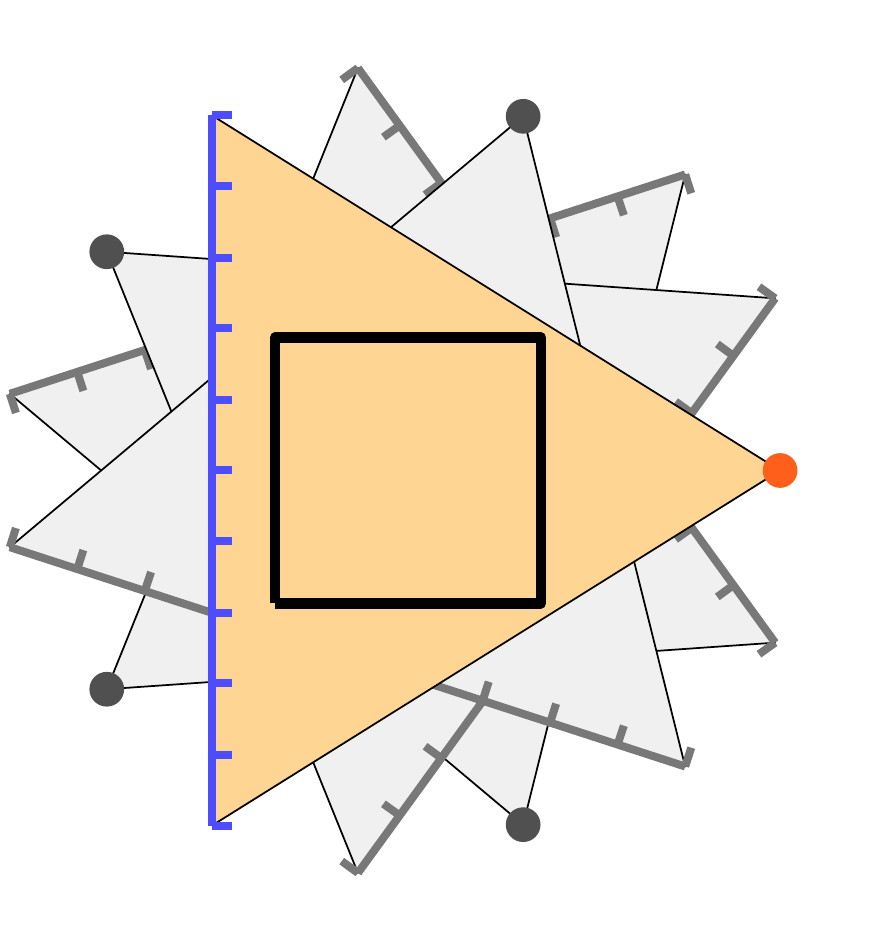}}
\put(-15,105){(a)}
\put(20,60){$\atten(x)$}
\put(30,5){$\atten\!\!:\!\Omega\!\rightarrow\!\R_+$}
\put(90,10){\includegraphics[width=110pt]{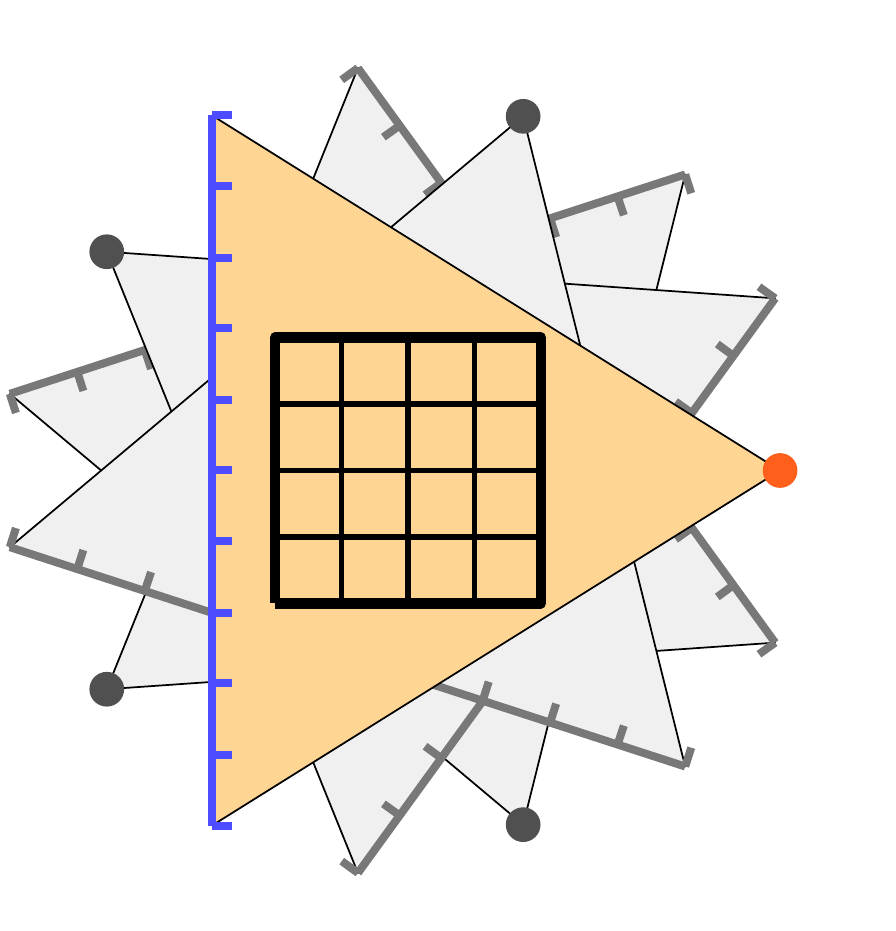}}
\put(95,105){(b)}
\put(140,5){$\f\in\R^{16}$}
\put(200,10){\includegraphics[width=110pt]{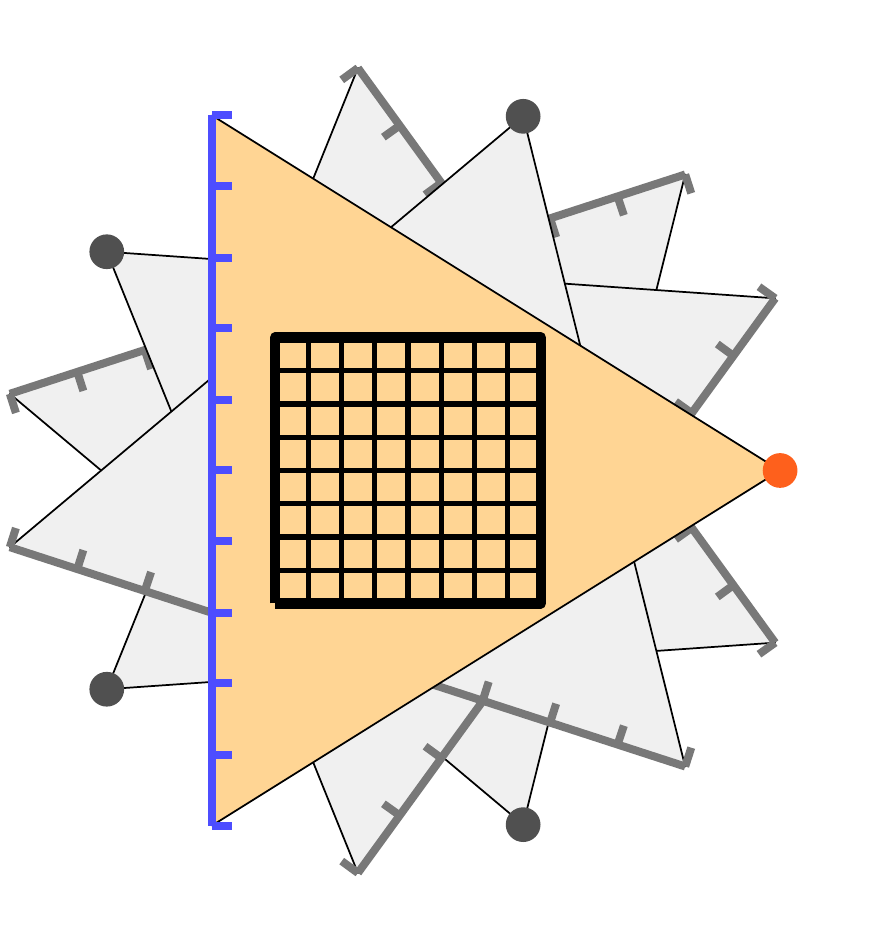}}
\put(205,105){(c)}
\put(250,5){$\f\in\R^{64}$}
\put(310,10){\includegraphics[width=110pt]{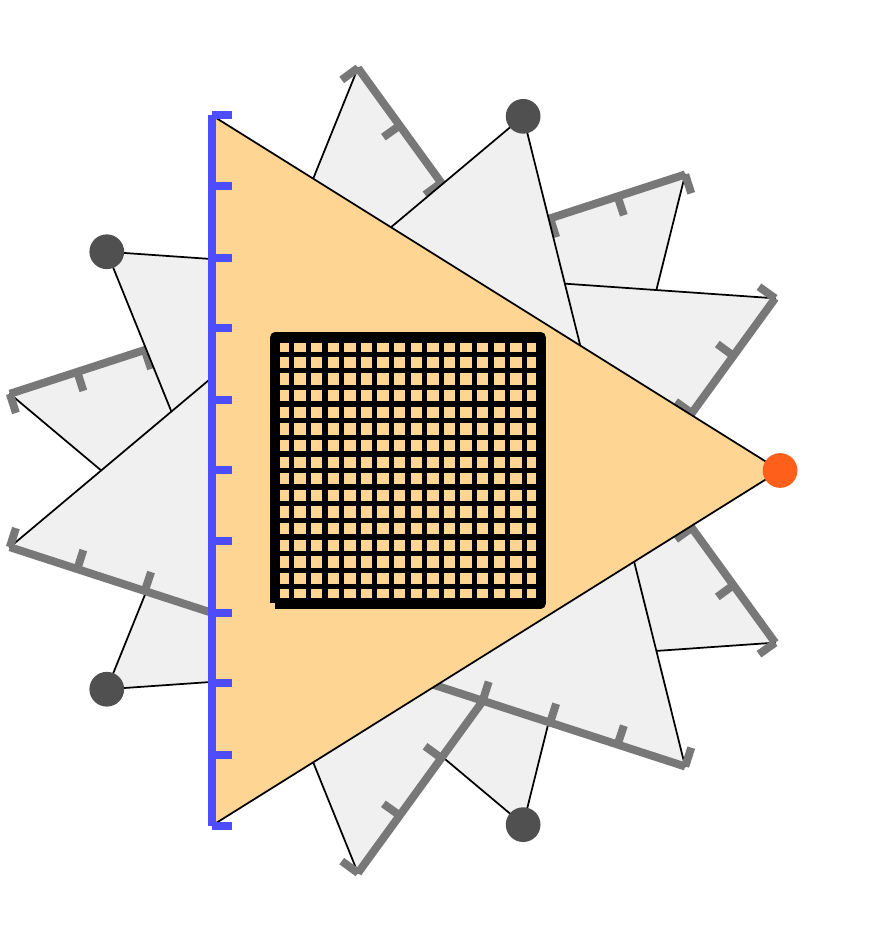}}
\put(315,105){(d)}
\put(360,5){$\f\in\R^{256}$}
\end{picture}
\caption{\label{fig:Xcontinuum}Schematic illustration about the relationship between continuous and discrete tomographic models. (a) Continuous model for two-dimensional X-ray tomography. The function $\atten(x)$ of equation (\ref{eqn:realmodelNoisy}) gives the point-wise X-ray attenuation coefficient inside the domain $\Omega$ shown here as a square. The dots show the five locations of the X-ray source. The detector has 10 pixels, so there is a total of 50 data points. (b-d) Discrete models with various resolutions. The attenuation coefficient is modelled as constant inside each pixel. The point of this figure is to demonstrate how the resolution of our computational model can be freely chosen while number of data points is fixed.}
\end{figure}

\noindent
In X-ray tomography one collects projection images of an unknown two- or three-dimensional body from different directions. For each direction, an X-ray source is placed on one side of the body and an X-ray detector on the opposite side.  After a calibration step involving a logarithm of the measured projection images, 
the X-ray tomography experiment can be
quite accurately modelled by 
\begin{equation}\label{eqn:realmodelNoisy}
 \greal=\Areal \atten+\mathbf{e},
\end{equation}
where vector $\greal\in\R^M$ is the data and $\atten:\Omega\rightarrow\R_+$ is the X-ray attenuation function defined in a bounded domain $\Omega\subset\R^n$. In this work, we take $\Omega$ to be the unit square $[0,1]^2$. The operator $\Areal$ represents integrals over $u$ over the areas or volumes connecting the X-ray source and a pixel in the detector (pencil beam model). See Figure \ref{fig:Xcontinuum}(a). Furthermore, $\mathbf{e}\in\R^M$ is a Gaussian random vector modelling measurement noise. 

The inverse problem is {\em  given the measurement geometry and noisy data $\greal$, recover a discrete approximation $\f\in\R^N$ to $\atten$.} Depending on the number and geometry of projection directions, this inverse problem is either mildly or strongly ill-posed \cite{Natterer1986,Mueller2012}. Ill-posedness means high sensitivity to inaccuracies in modelling and to errors in measurement. Overcoming this sensitivity requires {\em regularization}.

In this work, we study regularizing two-dimensional tomographic problems using Total Variation (TV) regularization introduced in \cite{Rudin1992}. We need a finite approximation to the continuous measurement model (\ref{eqn:realmodelNoisy}). Consider the discrete model
\begin{equation}\label{eqn:modelNoisy}
 \g=A\f,
\end{equation}
where $\f\in\R^N$ represents the discretized two-dimensional body under imaging. We divide the square $\Omega$ into $N=n\times n$ square pixels with side length $1/n$. The function $\atten$ is approximated by a function which is constant inside each pixel; those pixel values are listed as elements of the vector $\f$ (see Appendix \ref{Appendix_simple}). Further, $\g\in\R^M$ denotes the measurement data consisting of pixel values in digital X-ray projection images, and $A\in\R^{M\times N}$ is a matrix approximation to $\Areal$ (see \cite[Section 2.3.4]{Mueller2012}).

The discrete TV-regularized solution is defined by
\begin{equation}\label{def:discreteTV}
\fan  
:=\argmin_{\f\in \R^N_+}\left\{ \|A\f-\greal\|_2^2 + \frac{\alpha}{n}  \sum_{\kappa\sim \kappa^\prime}| f_\kappa- f_{\kappa^\prime}|\right\},
\end{equation}
where $\alpha>0$ is a so-called regularization parameter. We use the {\em anisotropic} definition of total variation: in (\ref{def:discreteTV}) the relation $\kappa\sim \kappa^\prime$ is true whenever elements $f_\kappa$ and $f_{\kappa^\prime}$ of the vector $\f$ correspond to values at either horizontally or vertically neighboring pixels.

How to choose the regularization parameter $\alpha$ in (\ref{def:discreteTV})? To our best knowledge, these are all the currently available methods in the literature:
\begin{enumerate}
\item The classical L-curve method \cite{Hansen1997}.
\item Discrepancy principle, introduced for TV regularization by Wen and Chan in \cite{Wen2012}  in the context of image restoration. This method needs {\em a priori} information about the measurement noise amplitude.
\item Two approaches (quasi-optimality principle and Hanke-Raus rules) described by Kindermann, Mutimbu and Resmerita in \cite{Kindermann2013}.
\item The so-called {\em S-curve method}, introduced in \cite{Kolehmainen2012} and implemented for wavelet-based tomography in \cite{Hamalainen2013}, can be extended to total variation regularization. Such extension is done in this paper for the first time, making use of {\em a priori} information about the sparsity of the gradient of the unknown.
\end{enumerate}

\noindent 
There is so far no parameter choice method that would be the best (or even perform acceptably) for {\em all} applications. Rather, it seems that for any given application there is a quest for finding a method that gives consistent, robust and useful results. For this reason we think that there needs to be a large collection of methods for parameter choice, based on different principles, for scientists and engineers to have a better chance of finding a good option for each practical situation.

We propose a new parameter choice method, based on solving the inverse problem at multiple resolutions. Namely, while the number $M$ of measurement points is fixed (it is simply the number of detector elements times the number of projection images), the number $N=n\times n$ of pixels in the reconstruction can be freely chosen. See Figure \ref{fig:Xcontinuum} for an illustration. 

Given a tomographic dataset, we can use definition (\ref{def:discreteTV}) to compute reconstructions in the same square domain but at varying resolutions $n$. We prove in Section \ref{sec:convergence} that if we keep $\alpha>0$  in equation (\ref{def:discreteTV}) fixed and let $n\rightarrow\infty$, then $\fan$ converges to a limit image.
The proof is an extension of the one-dimensional result \cite{Lassas2004} to dimension two. 

We propose the following practical method for choosing the parameter $\alpha>0$:\\
{\bf Parameter choice method.} Given a noisy dataset $\greal$, compute $\fan$ defined by (\ref{def:discreteTV}) for a set of $\alpha$ values and at two or more resolutions $n$, including the resolution intended for showing the final result. Calculate the discrete TV norms of the reconstructions and define the optimal $\alpha$ to be the smallest value that leads to TV norms that do not depend significantly on $n$.

At this point we do not claim that the above algorithm is a parameter choice rule in the strict sense as defined in \cite{Engl1996}. However, the numerical evidence presented in Section \ref{sec:results} suggests that the method works well and, as expected, leads to larger $\alpha$ when the noise level is increased.

In light of our new theorem, the TV norms of the reconstructions converge for any choice of $\alpha$ as $n\ra \infty$. How exactly can the above approach work at all? It seems that the ill-posedness of the tomography problem shows up as significant $n$-dependence of the TV norms of reconstructions whenever $\alpha$ is too small to regularize away the instability caused by the term $ \|A\f-\greal\|_2^2$. Since the final intended resolution is included in the test resolutions, our choice rule adapts to the relevant resolution and the noise level.

The reader may wonder why $\alpha$ is divided by $n$ in (\ref{def:discreteTV}). The TV regularization formulation for smooth functions takes the form
\begin{equation}\label{continuumTV}
 \argmin_{\atten\geq 0} \left\{\|\Areal \atten-\greal\|_2^2 + \int_\Omega |\nabla \atten| dx\right\}
\end{equation}
to the discrete minimization problem (\ref{def:discreteTV}). The division by $n$ comes from approximations with a finite difference quotient and a midpoint rule integration quadrature:
\begin{equation} \label{Riemann}
\int_\Omega \left|\frac{\partial \atten}{\partial x_1}\right| dx 
\approx 
\int_\Omega \left|\frac{\atten(x_1+\frac{1}{n},x_2)-\atten(x_1,x_2)}{\frac{1}{n}}\right| dx
\approx 
\left(\frac{1}{n}\right)^2 \sum_{\kappa\sim \kappa^\prime} \left|\frac{f_\kappa-f_{\kappa^\prime}}{\frac{1}{n}}\right|.
\end{equation}
In (\ref{Riemann}) the relation $\kappa\sim \kappa^\prime$ is true whenever elements $f_\kappa$ and $f_{\kappa^\prime}$ of $\f$ correspond to values at horizontally neighboring pixels. A similar computation using vertical differences is needed as well for the analysis regarding (\ref{continuumTV}).

Applying total variation regularization to tomographic problems dates back to 1998  \cite{Delaney1998} in case of simulated data and to 2003 \cite{Kolehmainen2003} in case of measured data. Since then, there have been many further studies: \cite{Kolehmainen2006,Liao2008,Sidky2008,Herman2008,Tang2009,
Duan2009,Bian2010b,Jensen2011,Tian2011, Jorgensen2011,Sidky2013,Chartrand2013,Chen2013}. 

There are many computational approaches for minimizing the total variation regularized least squares functional (\ref{def:discreteTV}). These include quadratic programming \cite{Lassas2004,Hintermuller2004,Goldfarb2005,Kolehmainen2012,Hamalainen2013}, lagged diffusivity method \cite{Dobson1997}, domain decomposition methods \cite{Fornasier2009,Fornasier2010}, Bregman distance methods \cite{Osher2005,Yin2008,Goldstein2009,Cai2009,Zhang2011}, primal-dual methods \cite{Chan1999,Chan2006,Esser2010,Nesterov2011}, finite element methods \cite{Feng2003,Bartels2012}, discontinuous Galerkin methods \cite{Hamalainen2014} and other methods \cite{Li1996,Vogel1998,Wang2008,Hale2010,Chambolle2011}. See also the books \cite{Vogel2002,Osher2003,Chan2005,Scherzer2009,Hansen2010}. In this paper we use a primal-dual interior-point quadratic programming method based on dividing the unknown to non-negative and non-positive parts in the spirit of \cite[Section 6.2]{Mueller2012}. One advantage of this approach is natural and effective enforcement of non-negativity of the attenuation coefficient.

Scientific novelties of this paper include 
\begin{itemize}
\item Proving the convergence of TV reconstructions as the resolution is increased.
\item Testing the S-curve method in the context of TV regularized tomography.
\item Introducing a novel multi-resolution parameter choice method and testing it with measured X-ray data.
\end{itemize}
While we discuss here the new multi-resolution parameter choice method for two-dimensional tomographic problems only, the computational approach generalizes directly to three-dimensional settings and other inverse problems as well. However, the convergence proof is currently available only in dimensions one \cite{Lassas2004} and two (Section \ref{sec:convergence} below).

\newpage

\section{The space $BV(D)$: definitions and approximations}\label{sec:BV}

\noindent
Let $D$ be the square $[0,1]^2\subset\R^2$. 
We use the following anisotropic definition for the space $BV(D)$ of functions of bounded variation defined on $D$.
\begin{definition}
The anisotropic space $BV(D)$ consists of functions $u\in L^1(D)$ for which $V(u)<\infty$. Here $V(u)$ is given for smooth functions by
\begin{equation}\label{def:V}
V(u) = V(u;D) = \int_D\Big(\left|\frac{\partial u(x)}{\partial x_1}\right|+\left|\frac{\partial u(x)}{\partial x_2}\right|\Big)dx.
\end{equation}
For general $u\in BV(D)$, $\left|\frac{\partial u(x)}{\partial x_j}\right|$ in formula (\ref{def:V}) is defined as total variation of the distributional derivative interpreted as a measure \cite{Attouch2014}. The BV norm is defined by
\begin{equation}\label{def:BVnorm}
\|u\|_{BV(D)} := \|u\|_{L^1(D)} + V(u;D).
\end{equation}
\end{definition}

\noindent
We will use two different concepts of weak convergence in the space  $BV(D)$.

\begin{definition}\label{def:weakBVconv}
A sequence $u_j\in BV(D)$ converges weakly in $BV(D)$ to limit $u$ as $j\to \infty$,
that is, $u_j\weakto u$, if 
$u_j$ and their distributional derivatives $\nabla u_j$ satisfy
\begin{equation}\label{weakBVconv}
\left\{\begin{array}{l}
\displaystyle \lim_{j\to \infty} \|u_j-u\|_{L^1(D)}=0,\quad\hbox{and }\\\\
\nabla u_j\to \nabla u \hbox{ weakly in ${\bf M}(D,\R^2)$ as $j\to \infty$}.
\end{array}\right.
\end{equation}
\end{definition}

\noindent
Above  ${\bf M}(D,\R^2)$ denotes $\R^2$-valued Borel measures on $D$. Note
that ${\bf M}(D,\R^2)$ is the dual space of $C(D;\R^2)$. 
Recall now from \cite[Def. 10.1.3]{Attouch2014} the definition of another useful type of weak convergence:
\begin{definition}\label{def:intermediate}
Let $u_j$ be a sequence in $BV(D)$ and $u\in BV(D)$. Then $u_j$ converges to $u$ as $j\ra\infty$ in the sense of {\em intermediate convergence} if and only if 
\begin{eqnarray}
\lim_{j\to \infty}\|u_j -u\|_{L^1(D)}
&=& \label{def:intermediate1}
0,\qquad \mbox{and}\\
\lim_{j\to \infty}\int_D |\nabla u_j(x)|dx 
&=& \label{def:intermediate2}
\int_D |\nabla u(x)| dx.
\end{eqnarray}
\end{definition}

\noindent
According to \cite[Proposition 10.1.2]{Attouch2014}, intermediate convergence of Definition \ref{def:intermediate} implies weak convergence of Definition \ref{def:weakBVconv}.

\begin{figure}
\begin{picture}(300,110)
\put(10,0){\includegraphics[height=3.5cm]{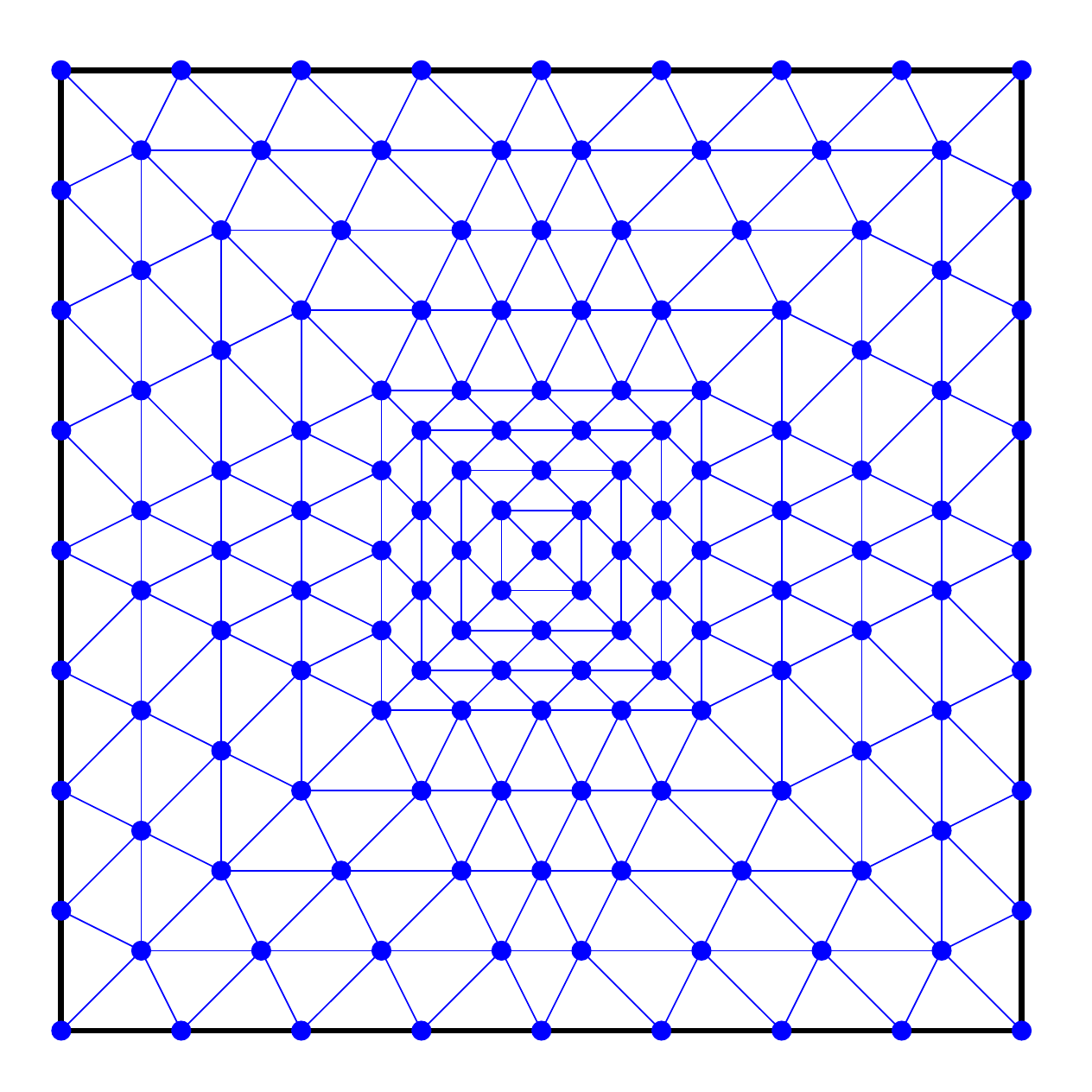}}
\put(170,4){\includegraphics[height=3.2cm]{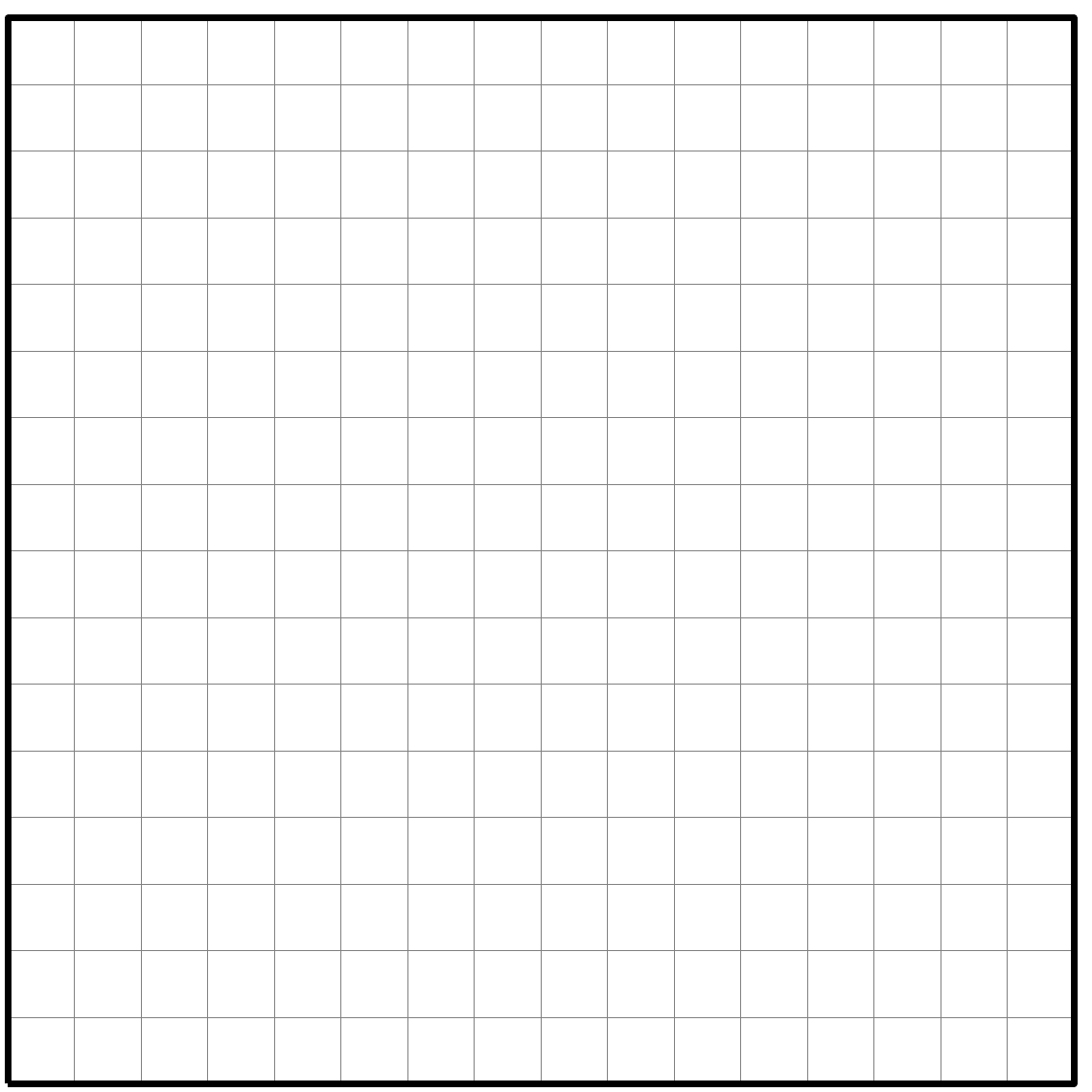}}
\put(0,95){(a)}
\put(155,95){(b)}
\end{picture}
\caption{\label{fig:triangulation}
Examples of different discretizations of the square domain $D=[0,1]^2$ used below overlaid in the same square. (a) Finite triangulation. (b) Square grid corresponding to $j=4$ and $n=16$. }
\end{figure}

We need a result on approximating elements of $BV(D)$ by functions that are piecewise constant on square pixels of fine enough resolution. For each $j>0$, consider the closed dyadic squares 
\begin{equation}\label{def:dyadicsquares}
  S_{j,\vec k}=S_{j,k_1,k_2}=[k_1 2^{-j},(k_1+1) 2^{-j}]\times [k_2 2^{-j},(k_2+1) 2^{-j}]\subset D,
\end{equation}
where $0\leq k_1\leq 2^j-1$ and $0\leq k_2\leq 2^j-1$. See Figure \ref{fig:triangulation}(b) for an illustration of a dyadic grid formed by such squares. Let 
$
  T_j:L^2(D)\to L^2(D)
$ 
be the orthogonal projections defined by averaging over squares:
\begin{equation}\label{def_Tn}
  (T_j u)\big|_{\mbox{int}(S_{j,\vec k})}=\frac{1}{|S_{j,\vec k}|}\int_{S_{j,\vec k}}u(x)\,dx.
\end{equation}
For any $j>0$, the set $Y_j:=\mbox{Range}(T_j)$ consists of $L^2(D)$ functions that are piecewise constant in the open squares belonging to the rectangular $2^j{\times}2^j$ grid (\ref{def:dyadicsquares}).

\begin{lemma}\label{approxlemma}
For all  $u\in BV(D)$ and $\e>0$ there is $j_1>0$  
such that there is a $v_2$ that is piecewise constant in the dyadic $2^{j_1}{\times}2^{j_1}$ grid (\ref{def:dyadicsquares}) and satisfies
\ba
\rho_{1,1}(u,v_2)<\e,
\ea
where the $\rho_{1,1}$-distance  is defined by 
\begin{equation}\label{def:rho11}
\rho_{1,1}(u,v)=\|u-v\|_{L^1(D)}+ | V(u)-V(v) |.
\end{equation}
\end{lemma}

\noindent
{\bf Proof.} Let $u\in BV$ and $\e>0$.  
We know from \cite{bvelik2003} that there is a function $v_0$ which is piecewise constant in a finite triangularization of $D$ and satisfies
\begin{equation}\label{u_ineq_1}
\rho_{1,1}(u,v_0)<\e/3.
\end{equation}
Set $M=\|v_0\|_{L^\infty(D)}$. Let us choose some notation for the concepts related to the fixed triangularization corresponding to $v_0$.
\begin{enumerate} 
\item Triangles denoted by $T_\nu\subset D$ with $\nu=1,2,\dots,T$.
\item Edges of the triangles are line segments $\gamma_\ell\subset D$ indexed by $\ell=1,2,\dots,L$. 
\item Vertices $p_r\in D$ with $r=1,2,\dots,R$. 
\item The minimum distance between two separate vertices is denoted by $\delta_0>0$.
\item The minimum angle between any two edges sharing an endpoint is called $\alpha_0>0$.
\item The length of the longest edge is denoted by $G$.
\end{enumerate}
See Figure \ref{fig:triangulation}(a) for an illustration of a finite triangularization.

\begin{figure}
\begin{picture}(350,110)
\put(0,0){\includegraphics[height=3.1cm]{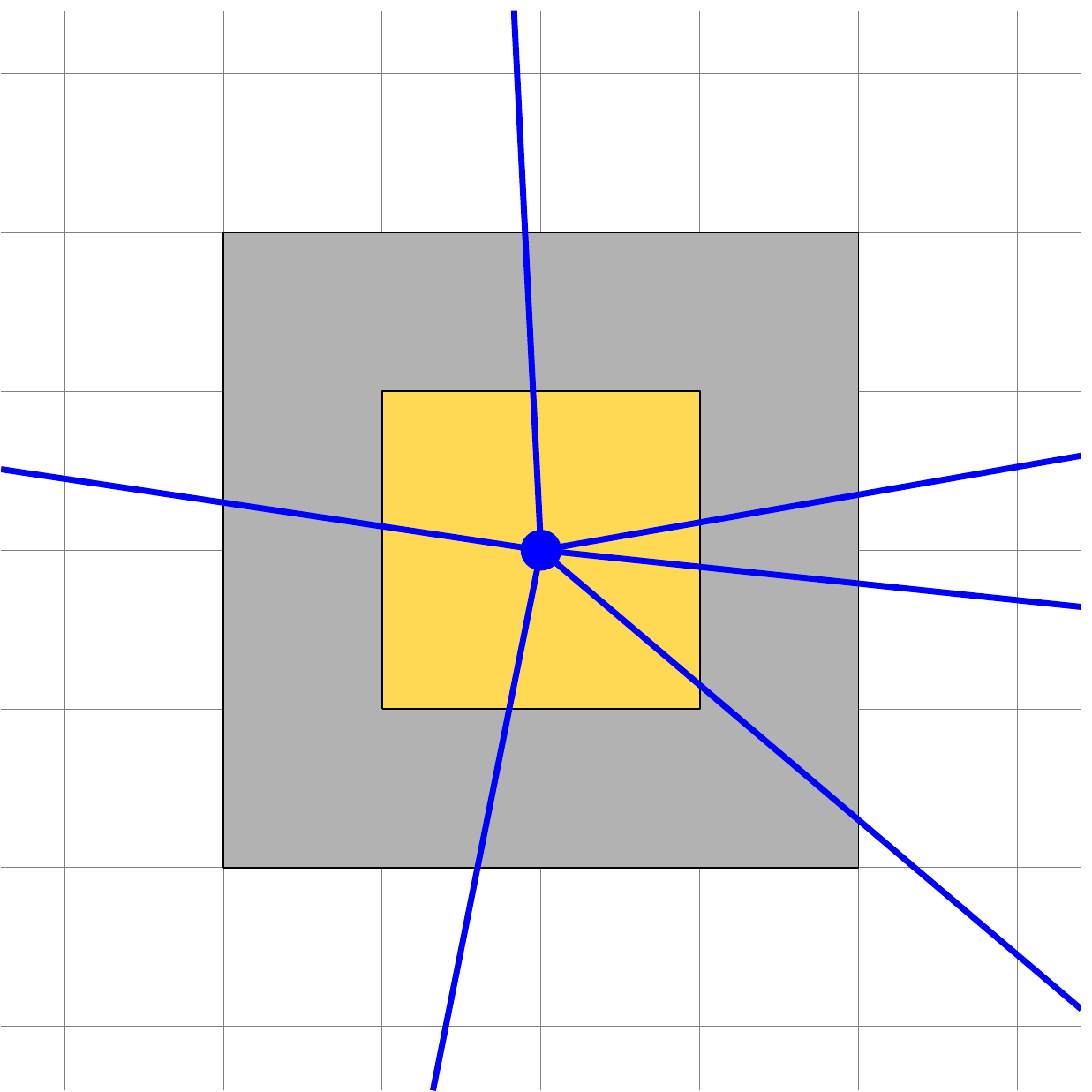}}
\put(95,0){\includegraphics[height=3.1cm]{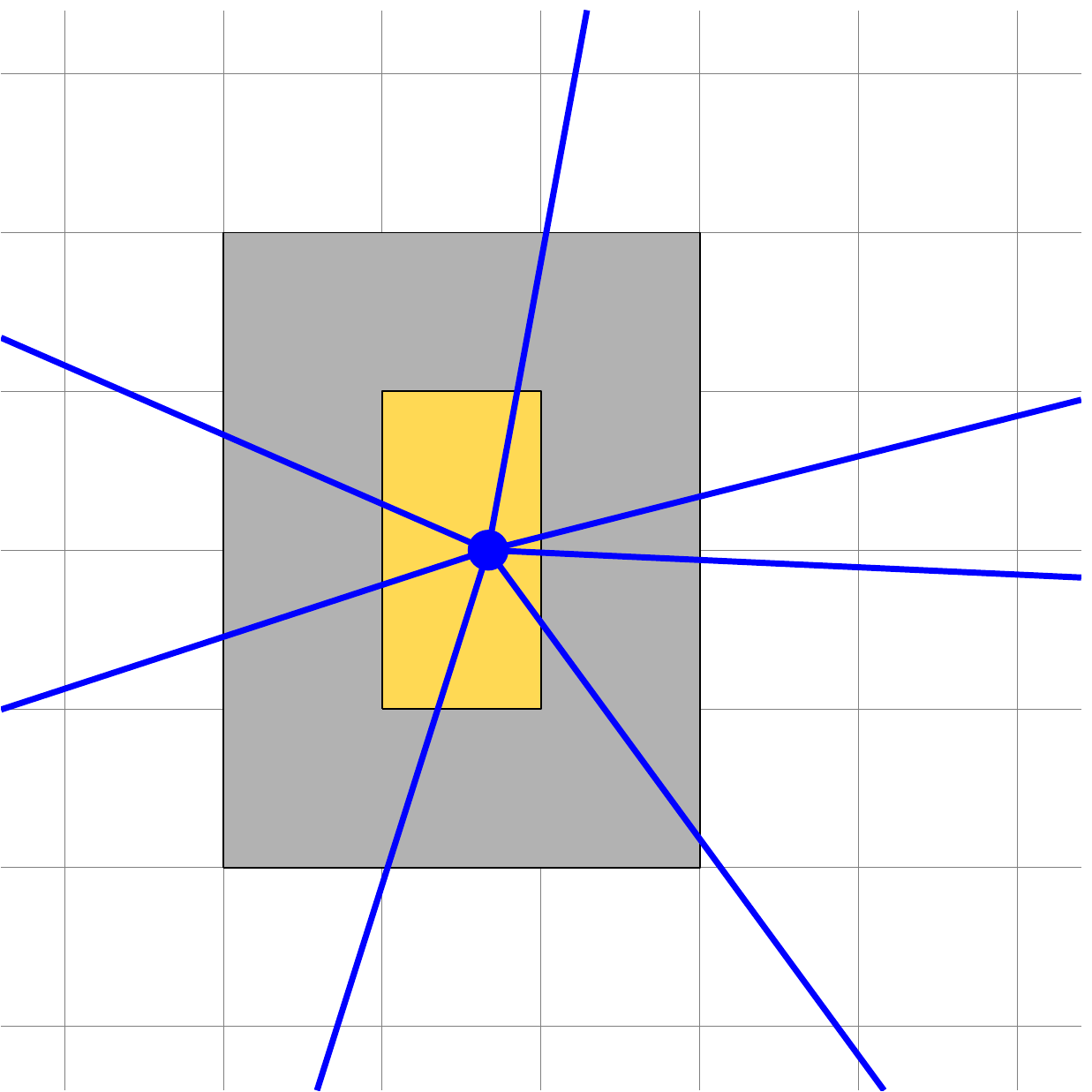}}
\put(190,0){\includegraphics[height=3.1cm]{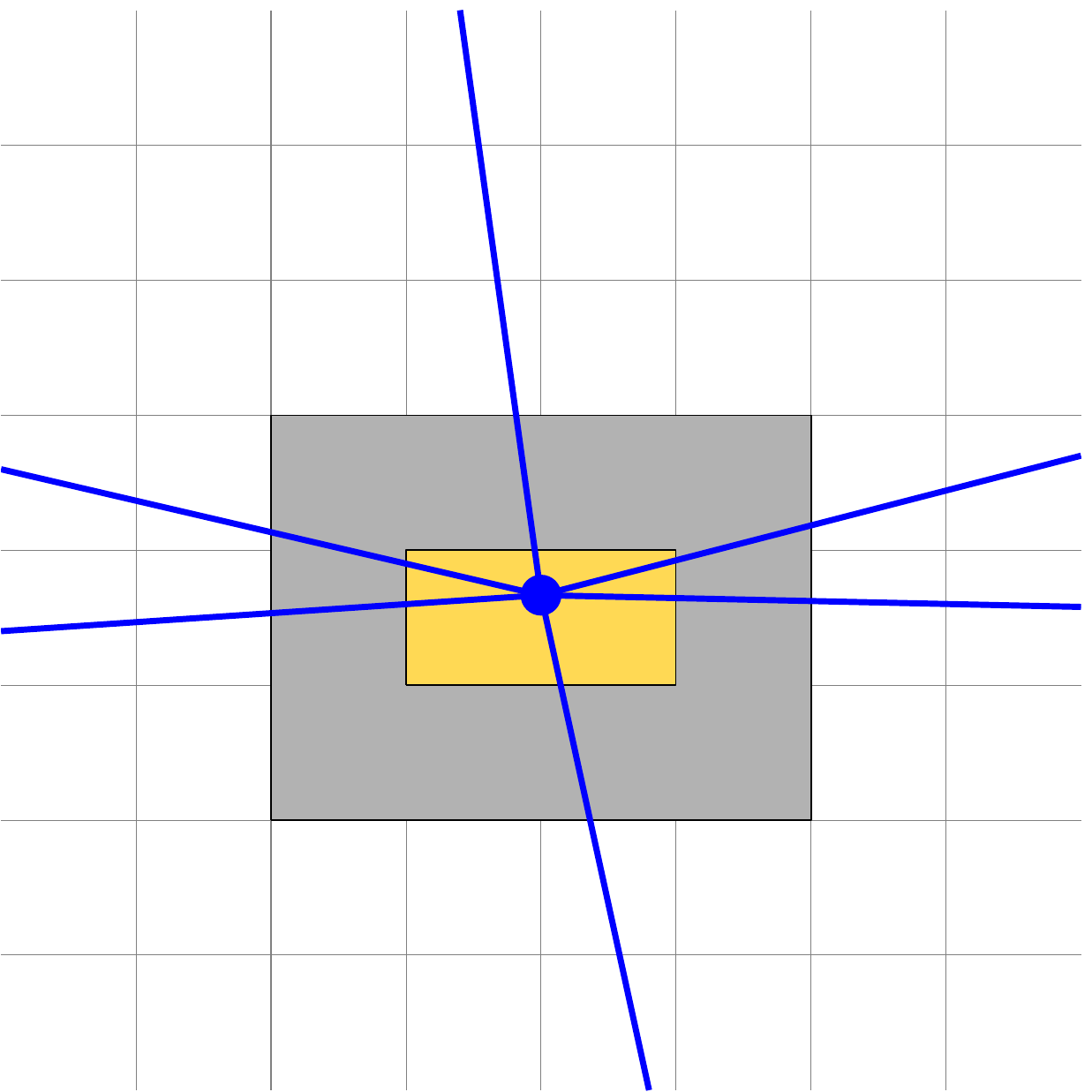}}
\put(285,0){\includegraphics[height=3.1cm]{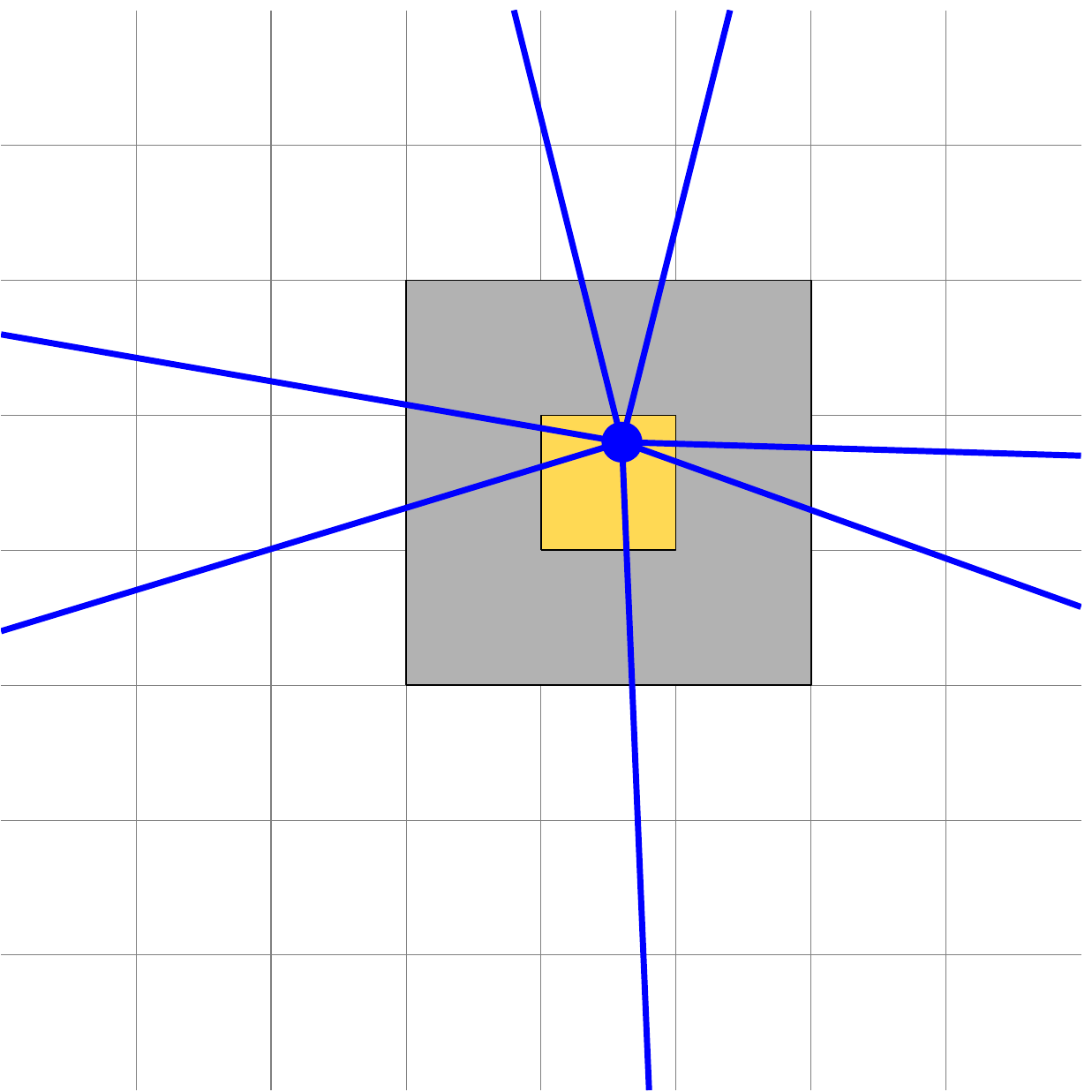}}
\put(0,95){(a)}
\put(95,95){(b)}
\put(190,95){(c)}
\put(285,95){(d)}
\put(18.5,59){\tiny $B_r$}
\put(32,48){\tiny $B^\prime_r$}
\end{picture}
\caption{\label{fig:vertex_in_grid}
Possible positions of a triangle vertex (blue dot) in the dyadic square grid. The triangle edges are shown in blue. The sets $B^\prime_r$ are shown as yellow and the clusters $B_r$ as light gray. (a) vertex at the corner of squares. (b) vertex at a horizontal edge between squares. (c) vertex at a vertical edge between squares. (d) vertex in the interior of a square.  }
\end{figure}
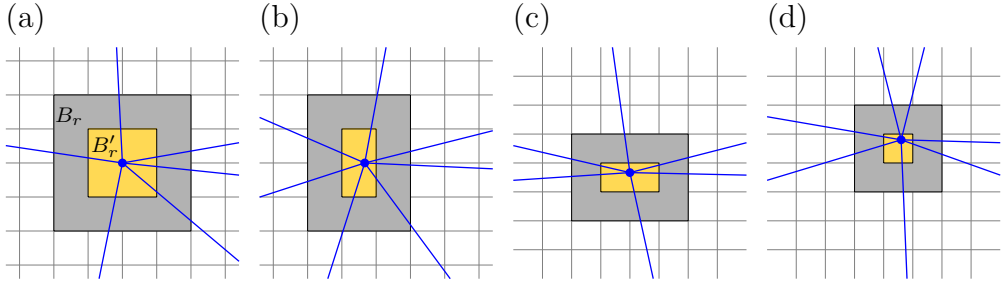

\begin{figure}
\begin{picture}(300,300)
\put(0,280){(a)}
\put(33,123){\includegraphics[height=5.7cm]{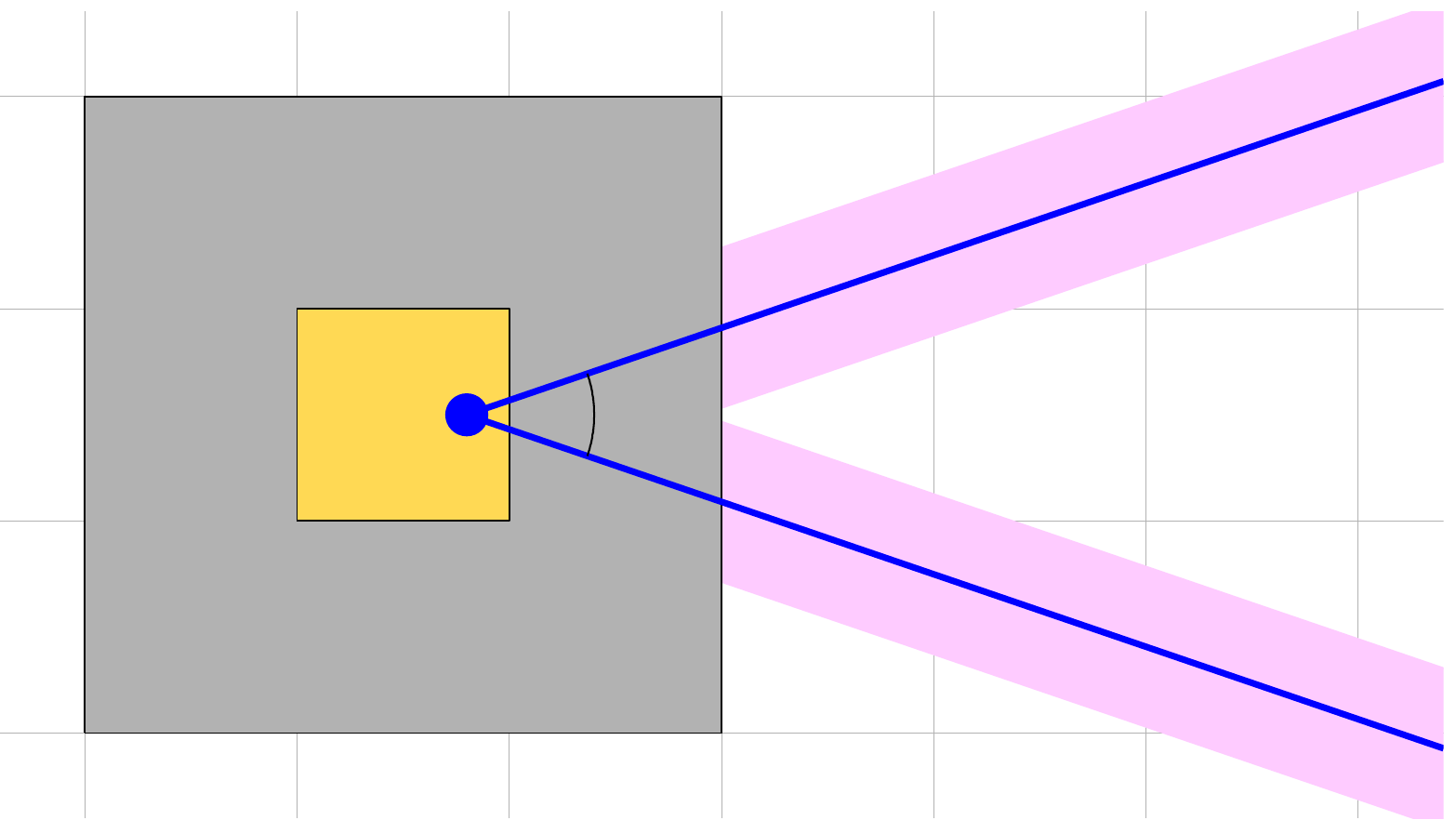}}
\put(113,208){\textcolor{blue}{$p_r$}}
\put(155, 200){\textcolor{black}{$\alpha_0$}}
\put(55,243){ $B_r$}
\put(95,187){$B^\prime_r$}
\put(323,266){\textcolor{blue}{$\widetilde{\gamma}_1$}}
\put(285,247){\rotatebox{20}{$\Gamma_{1,\varrho}$}}
\put(323,133){\textcolor{blue}{$\widetilde{\gamma}_2$}}
\put(285,154){\rotatebox{-20}{$\Gamma_{2,\varrho}$}}
\put(0,100){(b)}
\put(70,-5){\includegraphics[width=8.5cm]{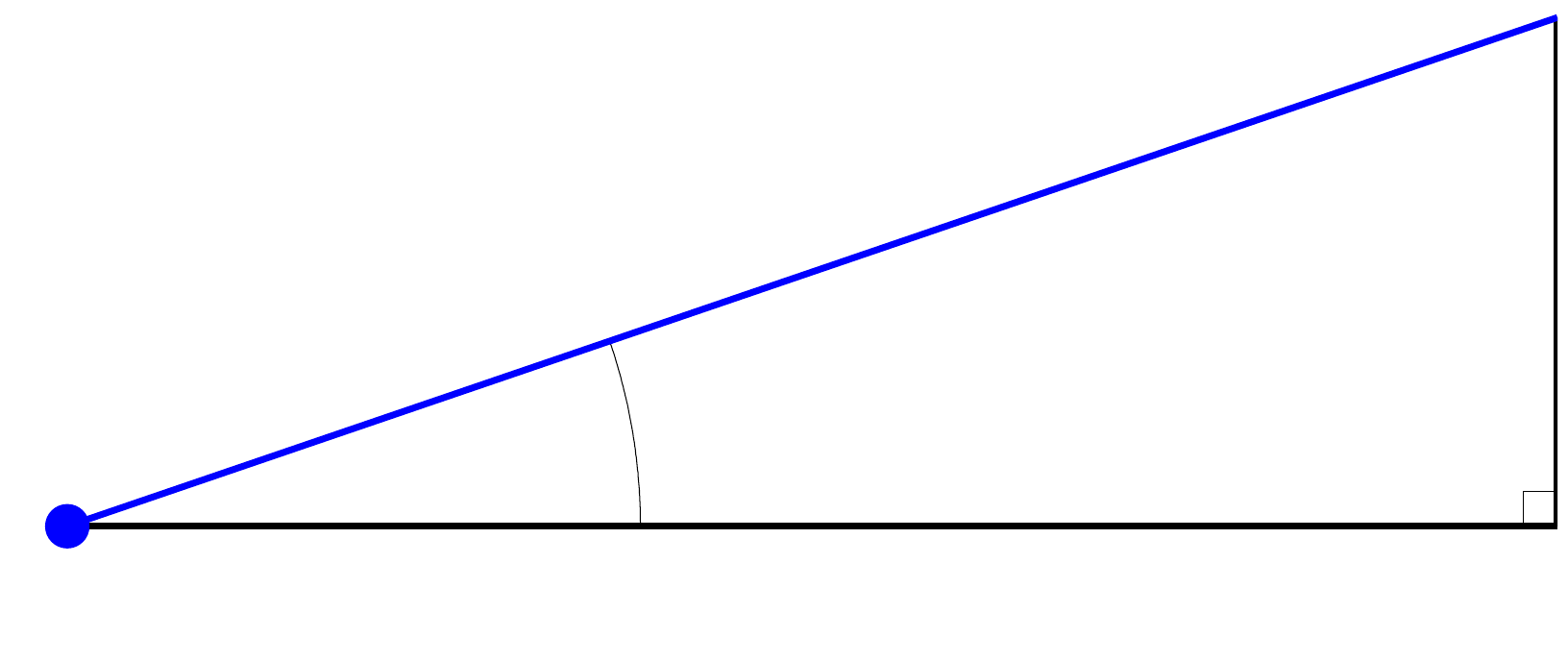}}
\put(220,0){$2^{-j_0}$}
\put(175,27){$\displaystyle\frac{\alpha_0}{2}$}
\put(315,53){$h_0$}
\end{picture}
\caption{\label{fig:vertex_worst}
(a) Illustration of two disjoint sets $\Gamma_{1,\varrho}$ and $\Gamma_{2,\varrho}$ (pink) as defined in (\ref{def:tubular}). The two edges $\widetilde{\gamma}_1$ and $\widetilde{\gamma}_2$, shown in blue, have the minimal angle $\alpha_0$ between them. The construction of $B_r$ and $B^\prime_r$ is done in the square grid of size $2^{-j_0}{\times}2^{-j_0}$. (b) The worst-case limit situation when $p_r$ is at the boundary of $B^\prime_r$ is where the definition (\ref{def:h_0}) comes from. }
\end{figure}
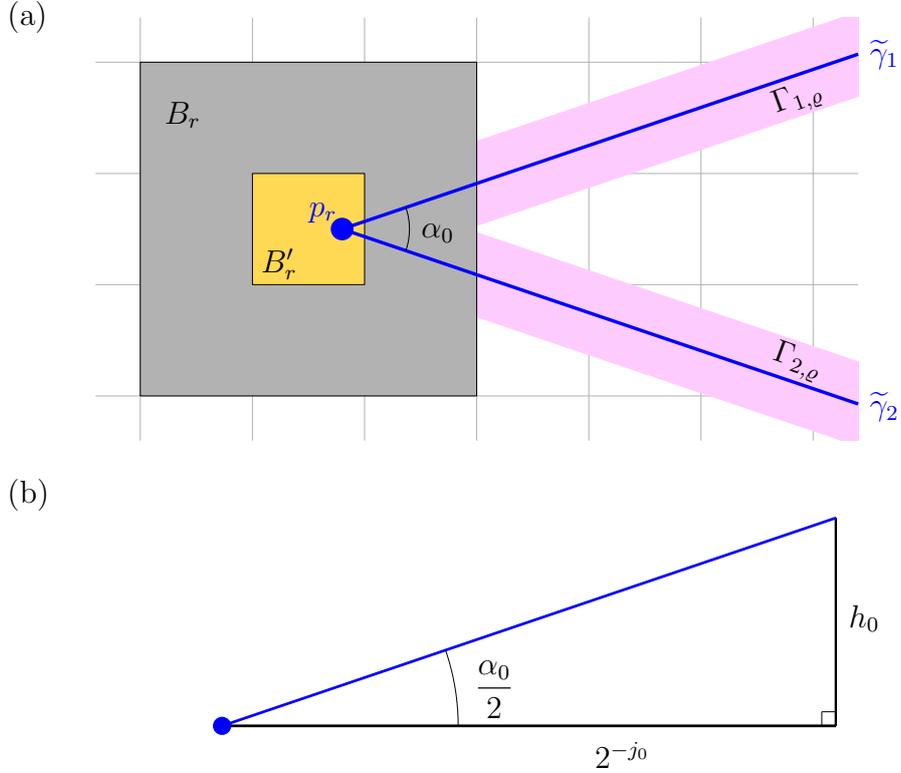

Next we define {\em clusters.}
Given $j>1$, let $B^\prime_r=B^\prime_{r,j}$ be union of all the closed squares $S_{j,\vec k}$ containing the vertex $p_r$. See the orange sets in Figure \ref{fig:vertex_in_grid}. Furthermore, define {clusters} $B_r=B_{r,j}$ by
\begin{equation}\label{def:Br}
  B_r = \bigcup \{S_{j,\vec k} \,|\, S_{j,\vec k}\cap B^\prime_r \not=\emptyset\}.
\end{equation}
Then $p_r\in B^\prime_r\subset B_r$. The clusters are shown in gray in Figures \ref{fig:vertex_in_grid}, \ref{fig:vertex_worst} and \ref{fig:skewdomain}. 

In the following, consider integers $j_0$ large enough for the following to hold:
\begin{equation}\label{def:j0}
  4\sqrt{2}\cdot 2^{-j_0}<\delta_0.
\end{equation}
Now (\ref{def:j0}) implies that the sets $B_{r}$ are disjoint and contain at most one vertex each.

Let $v_1$ be the function obtained from $v_0$ by replacing its value in each cluster $B_r$ by $0$. Then we have
\begin{eqnarray*}
  \|v_0-v_1\|_{L^1(D)} \leq \sum_{r=1}^R \int_{B_{r}}|v_0-v_1|dm\leq 16RM(2^{-j_0})^2
\end{eqnarray*}
since each cluster consists of at most 16 squares of area $(2^{-j_0})^2$. Furthermore, the $V$ expression (\ref{def:V}) can be very simply computed for functions that are piecewise constant with respect to either a finite triangulation or to a square grid. Namely, (\ref{def:V}) is the sum of (Manhattan distance) lengths of jump curve segments multiplied by the jump in function value over those curve segments. Therefore we can estimate
\begin{eqnarray*}
  |V(v_0)- V(v_1)| \leq 8LMR2^{-j_0} + 32MR2^{-j_0},
\end{eqnarray*}
where the first term in the right hand side comes from the triangularization. We used the facts that the Manhattan length of any edge restricted inside a cluster is at most $4\cdot 2^{-j_0}$ and that $v_0$ can jump at most $2M$ at an edge. The second term in the right hand side comes from the jump of $v_1$, bounded by $2M$, at the cluster boundary whose maximal length is $16\cdot 2^{-j_0}$. 

Now we can fix so large $j_0$ that (\ref{def:j0}) holds and that we have
\begin{equation}\label{u_ineq_2}
\rho_{1,1}(v_0,v_1)<\e/3.
\end{equation}

Next we construct tubular neighborhoods.
Denote the union of all clusters by 
$$
  \mathcal{B} = \bigcup_{r=1}^R B_r.
$$
Note that $ \mathcal{B}\subset \R^2$ is a closed set. Further, consider the parts of triangle edges that are outside the clusters:
\begin{equation}\label{def:gammatilde}
 \widetilde{\gamma}_\ell := \gamma_\ell\setminus \mathcal{B}.
\end{equation}

\begin{figure}
\begin{picture}(350,280)
\put(0,0){\includegraphics[height=10cm]{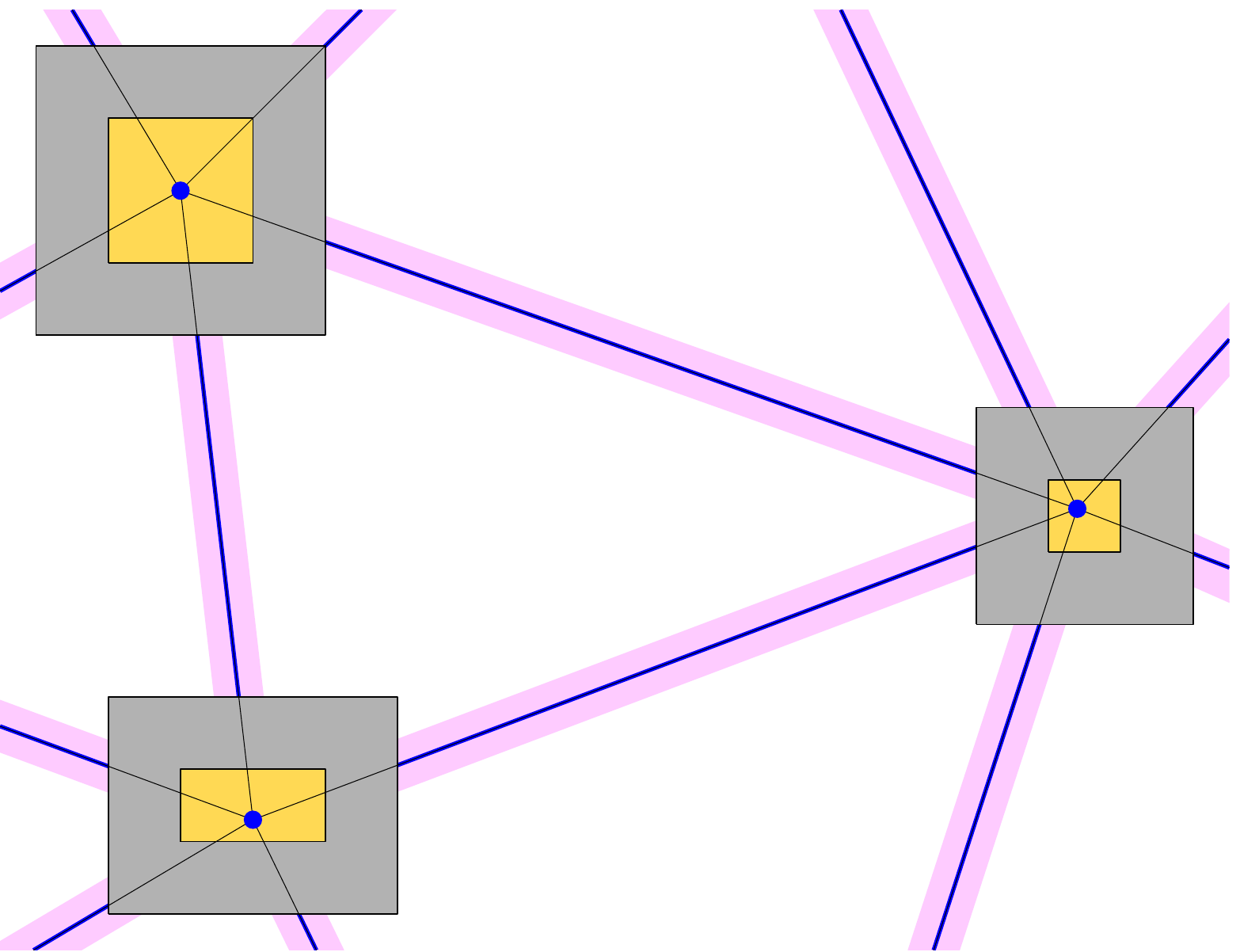}}
\end{picture}
\caption{\label{fig:rectangle_grid}
This is an illustration of the tubular neighborhoods  $\Gamma_{\ell,\varrho}$ (pink), defined in (\ref{def:tubular}), which are disjoint whenever $0<\varrho<h_0$ with $h_0$  given by (\ref{def:h_0}). The vertices $p_r$ are shown as blue dots and the ``triangle edges outside clusters,'' or line segments $ \widetilde{\gamma}_\ell$ defined in (\ref{def:gammatilde}), shown as blue lines.}
\end{figure}

For all $\varrho>0$, define 
\begin{equation}\label{def:tubular}
  \Gamma_{\ell,\varrho} := \{x\in\R^2 \,|\, \mbox{dist}(x, \widetilde{\gamma}_\ell)<\varrho\}\setminus  \mathcal{B},
\end{equation}
where ``dist'' refers to Euclidean distance. The sets $\Gamma_{\ell,\varrho}$ are open tubular neighborhoods of $ \widetilde{\gamma}_\ell$ in $D\setminus\mathcal{B}$. See Figures \ref{fig:vertex_worst}(a) and \ref{fig:rectangle_grid}.

We need to determine an upper bound for $\varrho>0$ ensuring that the sets $\Gamma_{\ell,\varrho}$ are disjoint.
It is clear from the construction of the sets $B_r$ and $B^\prime_r$ that (Euclidean) $\mbox{dist}(p_r,D\setminus\mathcal{B})>2^{-j_0}$. As illustrated in Figure \ref{fig:vertex_worst}(a), the minimal angle $\alpha_0$ between edges determines a lower bound for the distance between the intersection points of two edges with the boundary $\partial B_r$. A trigonometric consideration shown in Figure \ref{fig:vertex_worst}(b) implies that the tubular neighborhoods $\Gamma_{\ell,\varrho}$ are disjoint whenever $\varrho$ satisfies $0<\varrho<h_0$ with 
\begin{equation}\label{def:h_0}
  h_0 := 2^{-j_0}\tan\frac{\alpha_0}{2}.
\end{equation}
 See Figure \ref{fig:rectangle_grid} for an illustration of the tubular neighborhoods in a triangulation.


\begin{figure}
\begin{picture}(350,160)
\put(0,0){\includegraphics[width=13.5cm]{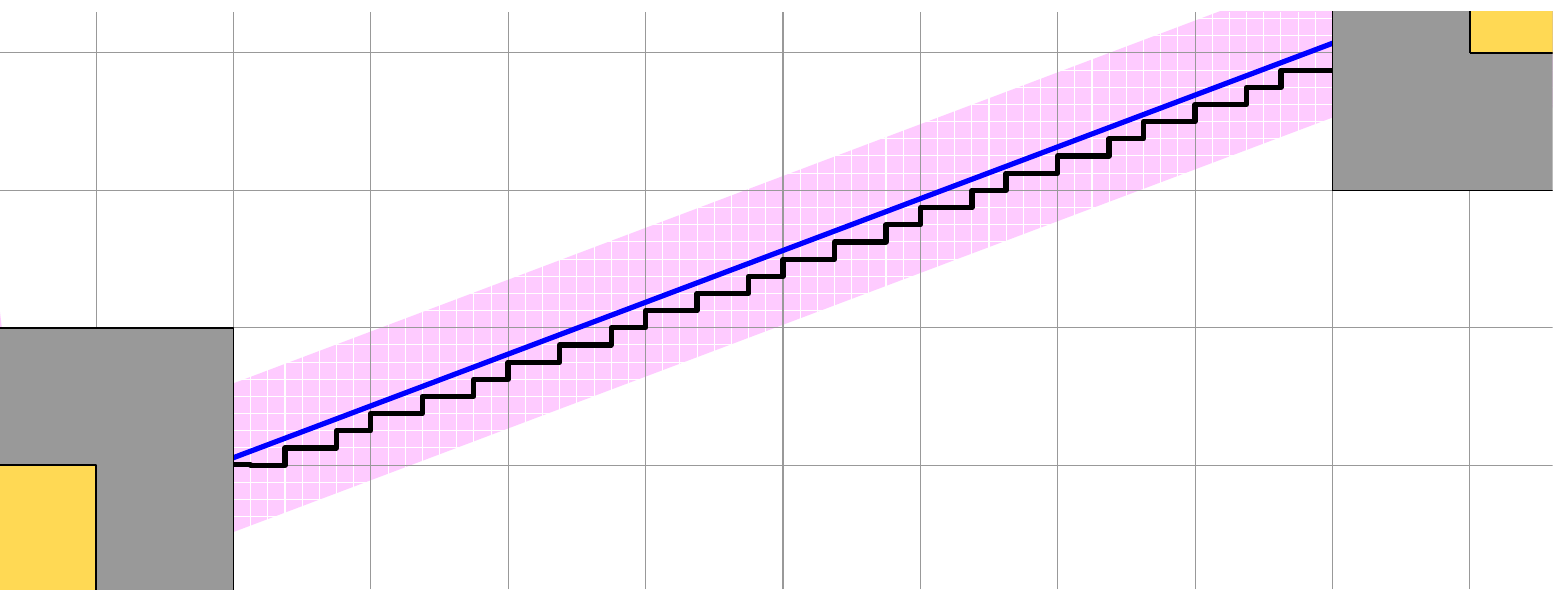}}
\put(40,50){$B_1$}
\put(70,45){\textcolor{blue}{$\widetilde{\gamma}_1$}}
\put(330,110){$B_2$}
\end{picture}
\caption{\label{fig:rectangle_grid3}
Refining the grid outside clusters. The bigger squares have size $2^{-j_0}{\times} 2^{-j_0}$ and the smaller squares have size $2^{-j_1}{\times} 2^{-j_1}$. A dyadic polygonal chain inside the tubular neighborhood $\Gamma_{1,\varrho}$ connects the clusters $B_1$ and $B_2$.}
\end{figure}
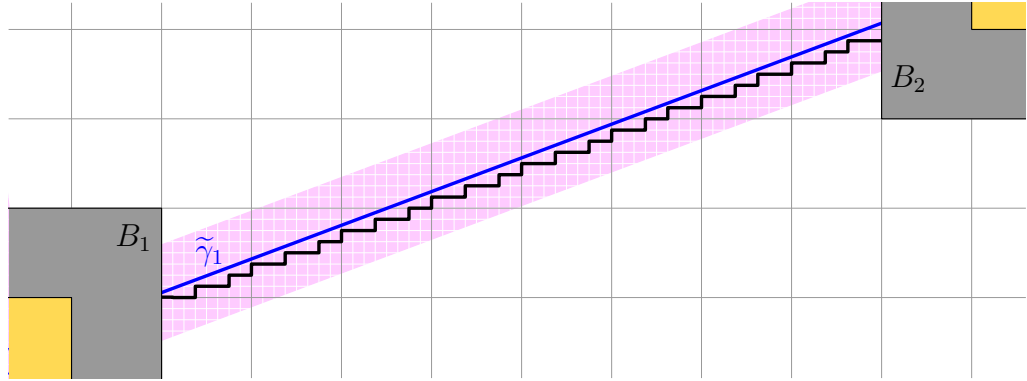

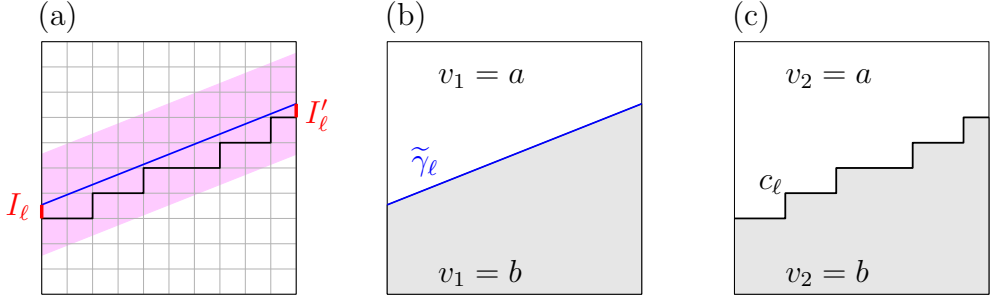
\begin{figure}
\begin{picture}(350,110)
\put(0,0){\includegraphics[height=3.4cm]{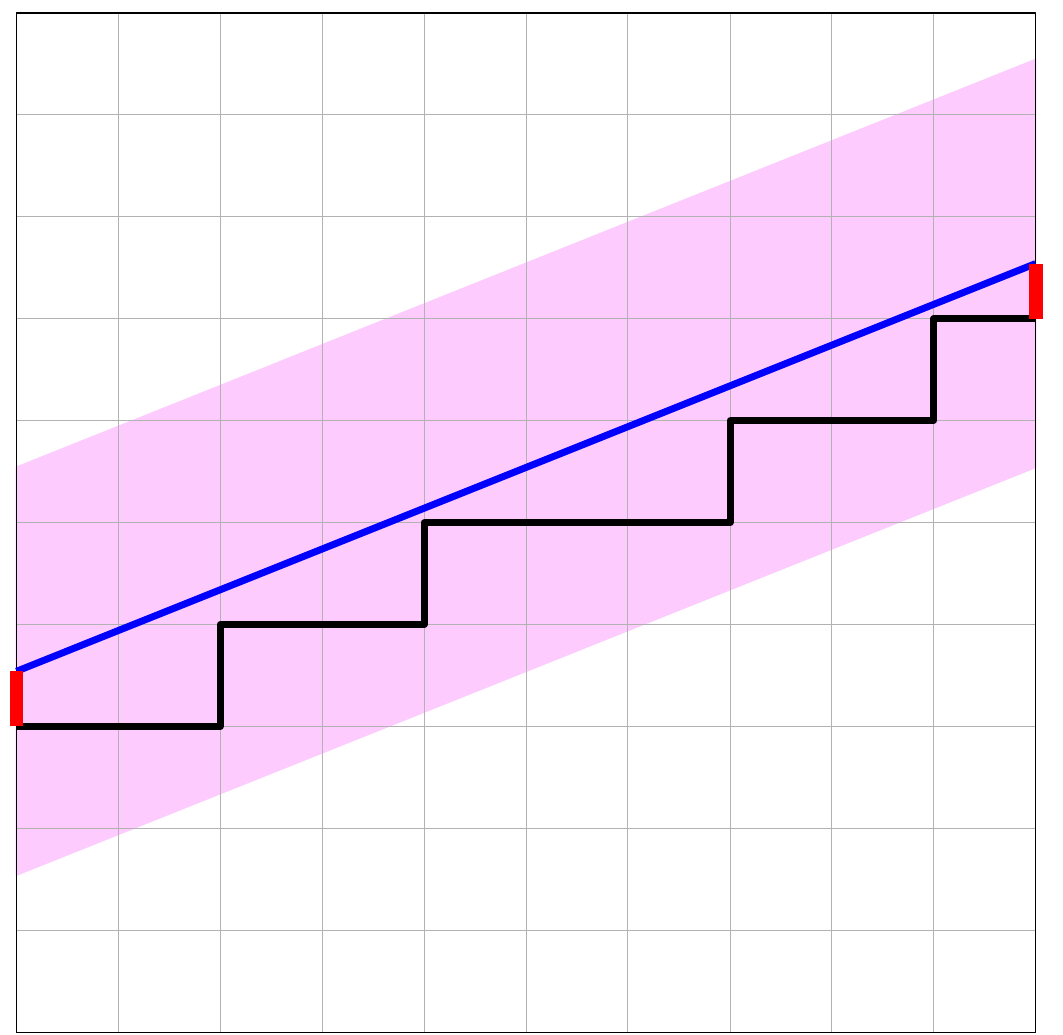}}
\put(130,0){\includegraphics[height=3.4cm]{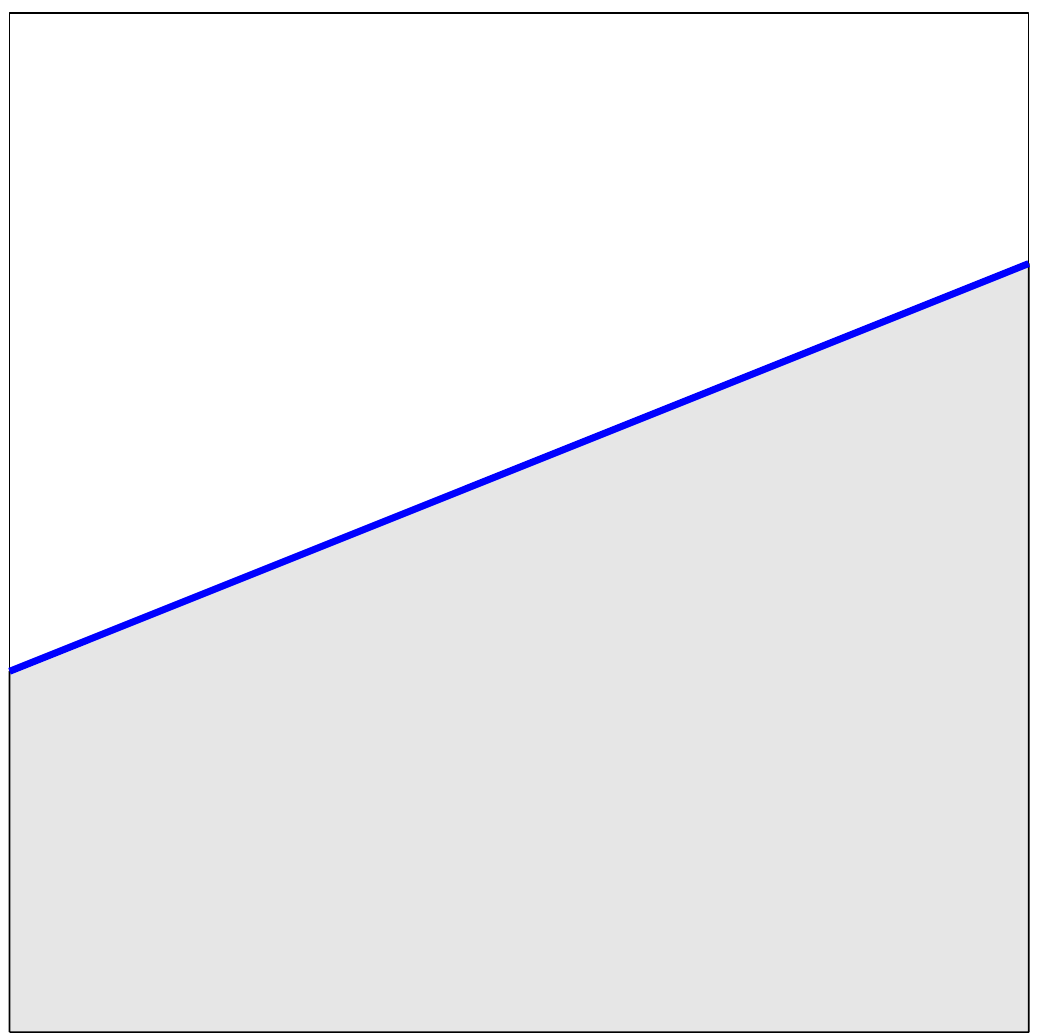}}
\put(260,0){\includegraphics[height=3.4cm]{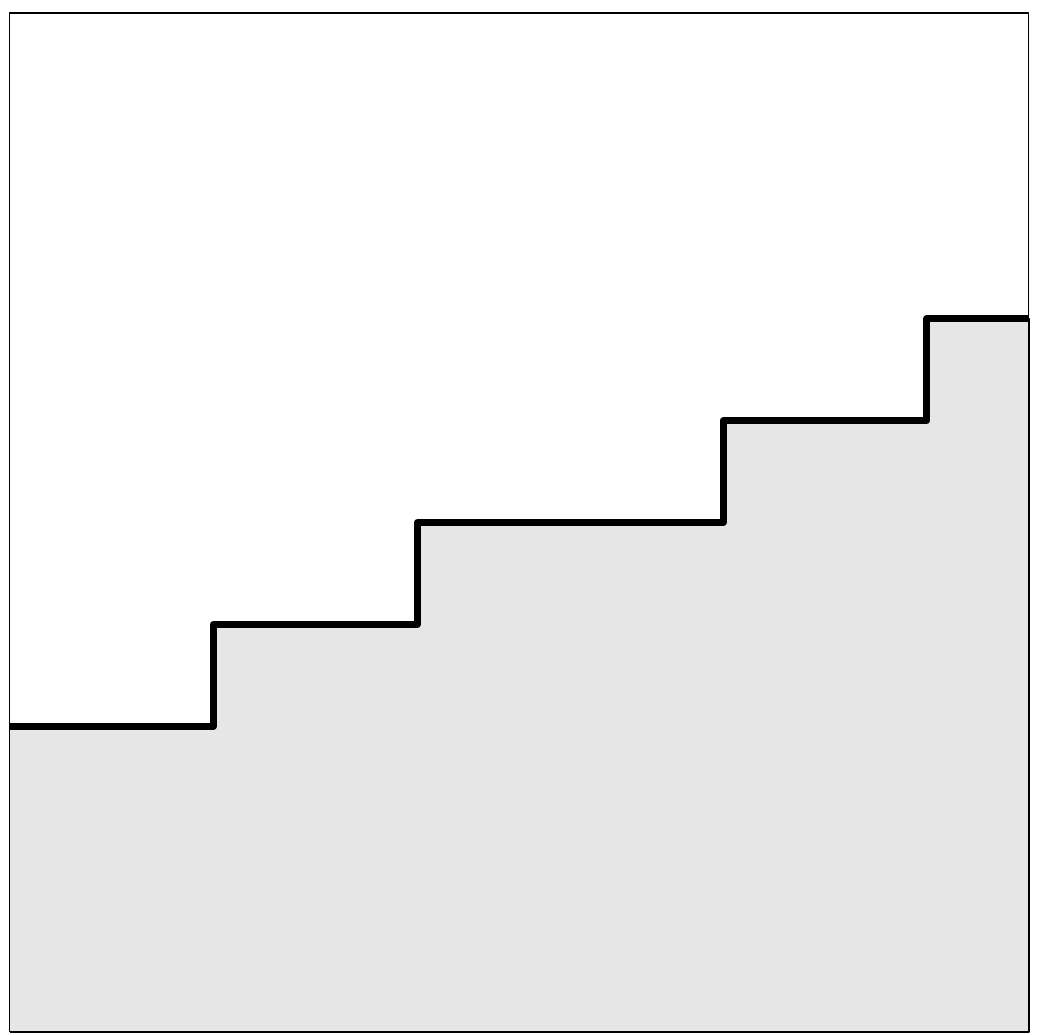}}
\put(0,102){(a)}
\put(-12,30){\textcolor{red}{$I_\ell$}}
\put(100,65){\textcolor{red}{$I_\ell^\prime$}}
\put(130,102){(b)}
\put(150,80){$v_1=a$}
\put(150,5){$v_1=b$}
\put(140,47){\textcolor{blue}{$\widetilde{\gamma}_\ell$}}
\put(260,102){(c)}
\put(280,80){$v_2=a$}
\put(280,5){$v_2=b$}
\put(270,40){$c_\ell$}
\end{picture}
\caption{\label{fig:skewdomain}
Square-grid based approximation of functions that are piecewise constant on a triangulation. (a) Detail from a situation similar to that shown in Figure \ref{fig:rectangle_grid3}. The straight blue line is the edge $\widetilde{\gamma}_\ell$, and the black jagged line is the dyadic polygonal chain $c_\ell$. The red intervals $I_\ell$ and $I_\ell^\prime$ are used for estimating the difference between the $BV$ norms of $v_1$ and $v_2$. (b) A function $v_1$ which is piecewise constant on a triangulation. (c) Approximation of $v_1$ by another function $v_2$ which is piecewise constant on a square grid. The difference of the BV norms of $v_1$ and $v_2$ comes entirely from jumps over the two red vertical intervals shown in (a).}
\end{figure}

Next we consider refining the square grid outside the clusters as shown in Figure \ref{fig:rectangle_grid3}. The smaller squares of size $2^{-j_1}{\times} 2^{-j_1}$ are defined by the formula (\ref{def:dyadicsquares}) with a large enough $j_1>j_0$.The idea is to take $j_1$ so large that we can fit a dyadic polygonal chains inside the tubular neighborhoods as shown in Figure \ref{fig:skewdomain}(a). 

Consider two vertices $p_1$ and $p_2$ connected by an edge $\gamma_1$. (There is no loss of generality since we can always renumber the vertices and edges.) Take $j_1$ so large that 
\begin{equation}\label{j1smallness1}
  2^{-j_1}<\frac{\varrho}{4}.
\end{equation}
Then we can connect the clusters $B_1$ and $B_2$ using a piecewise linear curve 
$$
  c_1:[0,1]\ra \overline{\Gamma}_{1,\varrho}
$$ 
with the following properties:
\begin{itemize}
\item The edge $\gamma_1$ and the dyadic polygonal chain $c_1$ are disjoint: $\gamma_1\cap c_1=\emptyset$.
\item The curve starts from the boundary of $B_1$ and ends at the boundary of $B_2$; in other words $c_1(0)\in \partial B_1$ and $c_1(1)\in \partial B_2$. 
\item The curve $c_1$ is close to the edge $\gamma_1$: 
\begin{equation}\label{curveclose}
  \min_{x\in\gamma_1}\|c_1(t)-x\|_2<3\cdot 2^{-j_1} \mbox{ for all } t\in[0,1].
\end{equation}
\item The range $c_1([0,1])$ of the curve consists of boundary segments of squares belonging to the $2^{-j_1}{\times} 2^{-j_1}$ grid.
\item The curve is injective: if $t_1\not= t_2$ then $c_1(t_1)\not= c_1(t_2)$.
\item The dyadic polygonal approximation $c_1$ is optimal in the following sense: if $c_1^\prime:[0,1]\ra \overline{\Gamma}_{1,\varrho}$ is another curve satisfying all of the above conditions for $c_1$, then $\mbox{length}(c_1^\prime([0,1]))\geq \mbox{length}(c_1([0,1]))$.
\end{itemize}
See Figure \ref{fig:rectangle_grid3} for illustration. Construct $c_\ell$ corresponding to each $\gamma_\ell$ with $\ell=1,\dots,L$. 

We move on to approximating functions using a square grid.
Consider the function $v_1$, which is piecewise constant in the triangulation of $D$, restricted to a tubular neighborhood $\Gamma_{\ell,\varrho}$. The line segment $\widetilde{\gamma}_\ell$ divides $\Gamma_{\ell,\varrho}$ into two components, and $v_1$ has a constant value in each component. See Figure \ref{fig:skewdomain}(b).

Now define another function $v_2$, which is piecewise constant on the square grid of size $2^{-j_1}{\times} 2^{-j_1}$. Outside the tubular neighborhoods $\Gamma_{\ell,\varrho}$  we set $v_1\equiv v_2$. Inside $\Gamma_{\ell,\varrho}$ we define the values of $v_2$ uniquely as follows. Consider the dyadic polygonal chains $c_\ell$  that are located close to the edges $\widetilde{\gamma}_\ell$. Now $v_2$ is allowed to have discontinuities only along the curves $c_\ell$ or at the cluster boundaries $\partial B_r$. See Figure  \ref{fig:skewdomain}(c).

By the construction of the dyadic polygonal chains the area between an edge $\widetilde{\gamma}_\ell$ and the polygonal chain is bounded by $3G2^{-j_1}$. Therefore we can estimate 
\begin{equation}\label{ruutuest1}
  \|v_1-v_2\|_{L^1(D)} \leq 6MG2^{-j_1},
\end{equation}
since we have $|v_1(x)-v_2(x)|\leq 2M$ for almost all $x\in D$. 

Now both $v_1$ and $v_2$ are piecewise constant, and the domains of constant value are polygons. The expressions  $V(v_1)$ and $V(v_2)$ can thus be calculated by summing the following quantities over the linear boundary segments: the jump in function values over the segment multiplied by the {\em Manhattan length of the segment}. Note that in $|V(v_1)-V(v_2)|$ all of the terms cancel each other except those coming from boundary segments located in the tubular neighborhoods $\Gamma_{\ell,\varrho}$. Due to the construction of $c_\ell$, the Manhattan and Euclidean lengths of $c_\ell$ are the same. Further, the Manhattan length of $\widetilde{\gamma}_\ell$ equals the length of $c_\ell$, apart from the difference between the two red intervals $I_\ell$ and $I_\ell^\prime$ shown in Figure \ref{fig:skewdomain}(a). By (\ref{curveclose}) we have 
$$
  \big| | I_\ell |- | I_\ell^\prime| \big| \leq   |I_\ell|+|I_\ell^\prime| \leq 6\cdot 2^{-j_1}.
$$
This together with the maximal jump size $2M$ allows us to estimate
\begin{equation}\label{ruutuest2}
|V(v_1)-V(v_2)|\leq  12ML2^{-j_1}.
\end{equation}
Using (\ref{ruutuest1}) and (\ref{ruutuest2}) we see that for large enough $j_1$ we have
$\rho_{1,1}(v_1,v_2)\leq  \e/3$, and the proof is complete.\hfill$\blacksquare$

\clearpage

\section{Multi-resolution convergence analysis}\label{sec:convergence}

\noindent
Define a functional $S_\infty: BV(D)\ra \R$ by
\begin{equation}\label{mainfunctional}
S_\infty(u)=\|\Areal u-m\|_{L^2(\Omega)}^2+\a_1 \|u\|_{L^1(D)}+\a V(u)
\end{equation}
with positive regularization parameters $\a_1>0$ and $\a>0$, and $V$ given by (\ref{def:V}). 
The functional $S_\infty$ defines an infinite-dimensional  minimization problem which we will approximate by discrete minimization problems. For this we define functionals
$S_j: BV(D)\to \R\cup\{\infty\}$ by
\begin{equation}\label{mainfunctionaln}
S_j(u) = 
\left\{\begin{array}{cl}
S_\infty(u), & \mbox{for }u\in \mbox{Range}(T_j),\\
\infty, & \mbox{for }u\not \in \mbox{Range}(T_j).
\end{array}\right.
\end{equation}
Denote
\begin{equation}\label{def:infS}
s_\infty := \inf_{u\in BV} S_\infty(u), \qquad  s_j := \inf_{u\in BV} S_j(u).
\end{equation}
We study the existence and properties of minimizers $\uMAP_j,\uMAP_\infty\in BV(D)$ satisfying $S_j(\uMAP_j)=s_j$ and $S_\infty(\uMAP_\infty)=s_\infty$, respectively. Such minimizers may not be unique. Therefore, we use the notations 
\begin{eqnarray*}
\argmin(S_j) &=& \{u\in BV(D):\ S_j(u)=s_j\},\\
\argmin(S_\infty) &=& \{u\in BV(D):\ S_\infty(u)=s_\infty\}.
\end{eqnarray*}

\noindent
We are now ready to state our main theorem.
\begin{theorem}\label{thm:main}
Assume either (A) or (B) about the linear  operator $\Areal$:
\begin{itemize} 
\item[{(A)}] $\Areal:L^2(D)\to L^2(\Omega)$ is compact and $\Areal:L^1(D)\to \mathcal D^\prime(\Omega)$ is continuous with some open and bounded set $\Omega\subset\R^2$,
\item[{(B)}] $\Areal:L^1(D)\to \R^M$ is bounded.
\end{itemize}
Let the functionals $S_j$ and $S_\infty$ be given by (\ref{mainfunctionaln}) and (\ref{mainfunctional}), respectively, and define $s_j$ and $s_\infty$ as in (\ref{def:infS}). Then 
\begin{itemize}
\item[(i)] $\displaystyle \lim_{j\to \infty}s_j= s_\infty$. 
\item[(ii)] There exists a minimizer $\uMAP_j\in\argmin(S_j)$ for all $j=1,2,3,\dots$
\item[(iii)] There exists a minimizer $\uMAP_\infty\in\argmin(S_\infty)$. 
\item[(iv)] Any sequence $\uMAP_j\in\argmin(S_j)$ of  minimizers has a subsequence $\uMAP_{j(\ell)}$ that converges weakly in $BV(D)$ to some limit $w\in BV(D)$, or more briefly $\uMAP_{j(\ell)}\weakto \,w$ as $\ell\to\infty$. Furthermore, $\displaystyle \lim_{\ell\to\infty}V(\uMAP_{j(\ell)})=V(w)$,
 that is, the subsequence $\uMAP_{j(\ell)}$  converges to $w$  in the intermediate topology of $BV(D)$.
\item[(v)] The limit $w=:\uMAP_\infty$ in (iv) is a minimizer: $\uMAP_\infty\in \argmin(S_\infty)$.
\end{itemize}
\end{theorem}

\noindent
The set $\Omega$ in (\ref{mainfunctional}) can be an open rectangle with $L^2(\Omega),$ modeling sinogram space, or $\Omega=\{1,2,\dots,M\}$ equipped with the counting measure, leading to $L^2(\Omega)=\R^M$. These two options correspond to assumptions (A) and (B) in Theorem \ref{thm:main}.

The connection between Theorem \ref{thm:main} and the computational explanation in the Introduction is this. In the computational part of this paper we take $\a_1=0$, which roughly corresponds to having a really small $\alpha_1>0$ in the theory. The minimizing vector $\fan\in\R^{n{\times}n}$ defined in (\ref{def:discreteTV}) can be seen as an $n{\times}n$ matrix with $n= 2^j$. The matrix $\fan$ is a 2D array of the values of a piecewise constant minimizer function $\uMAP_j\in\argmin(S_j)$ in the pixels (\ref{def:dyadicsquares}). Thus we see that for any choice of regularization parameter $\alpha$ in (\ref{def:discreteTV}), the total variation norms of $\fan$ with $n= 2^j$ converge to a limit as $j\ra\infty$. The measurement operator discussed in the Introduction is as in assumption (B).

Remark: our convergence  result is quite delicate. Even a slight change of norms will destroy the argumentation. See Appendix \ref{Appendix_example} for more details about this. 

Before proving the main theorem we prove a couple of lemmata.

\begin{lemma}\label{lemma:Aterm}
Let the assumptions of Theorem \ref{thm:main} hold. Assume that $u_j\to u$ as $j\ra\infty$ in the weak topology of $BV(D)$ given in  Definition \ref{def:weakBVconv}. Then 
\begin{equation}\label{Areal_conv}
\lim_{j\to \infty} \Areal u_j=\Areal u 
\end{equation}
with convergence in $L^2(\Omega)$ or in $\R^M$ in the cases of Assumption (A) or (B) in Theorem \ref{thm:main}, respectively. Moreover, it follows that
\begin{eqnarray}
&& \nonumber
\|\Areal u-m\|_{L^2(\Omega)}^2+\a_1 \|u\|_{L^1(D)}+\a V(u)\\
&\leq& \label{upperboundineq_expanded2}
\liminf_{j\ra \infty}\left(\|\Areal u_j-m\|_{L^2(\Omega)}^2+\a_1 \|u_j\|_{L^1(D)}+\a V(u_j)\right).
\end{eqnarray}
\end{lemma}

\noindent
{\bf Proof.} Assume (A)  of Theorem \ref{thm:main}. 
By \cite[Thm. 10.1.3]{Attouch2014} we have
$$
\|u\|_{L^2(D)}\leq C_1 \|u\|_{BV(D)}.
$$
Moreover, there is $c_1$ such that $ \|u_j\|_{BV(D)}\le c_1$ for all $n$. Hence,
$ \|u_j\|_{L^2(D)}\le C$ for all $j$. To show that $\Areal u_j$
converges in $L^2(\Omega)$
to $\Areal u$ as $j\to \infty$, assume the opposite. Then there is $\e>0$
and a subsequence
$j(\ell)\to \infty$ as $\ell\to \infty$ such that
\beq\label{e1}
\|\Areal u_{j(\ell)}-\Areal u\|_{L^2(\Omega)}\ge \e.
\eeq
Now $\|u_{j(\ell)}\|_{L^2(D)}\leq C$ and as $\Areal:L^2(D)\to L^2(\Omega)$ is compact,
there is a subsequence $u_{j(\ell_p)}$ and $w\in L^2(\Omega)$ such that
\beq\label{e2}
\lim_{p\to \infty} \|\Areal u_{j(\ell_p)}-w\|_{L^2(\Omega)}=0.
\eeq
On the other hand, as $u_j\weakto u$,
we have that $u_j\to u$ in the strong topology  of $L^1(D)$ by
definition. This implies that
\beq\label{e3}
\lim_{p\to \infty} \Areal u_{j(\ell_p)}=\Areal u ,\quad\hbox{in }\mathcal D^\prime(\Omega).
\eeq
The above two limits (\ref{e2}) in $L^2(\Omega)$ and  (\ref{e3})  in
$\mathcal D^\prime(\Omega)$ imply that
$\Areal u=w$. This, (\ref{e1}), and (\ref{e2}) are in contradiction. Hence
we have (\ref{Areal_conv})  in the case of assumption (A) in Theorem \ref{thm:main}. 

{\mnewtext Inequality  (\ref{upperboundineq_expanded2}) follows from the above and \cite[Prop.\ 10.1.1]{Attouch2014}.}

Assumption (B) can be handled by adding a projection operator similar to (\ref{def_Tn}).
\phantom{m}\hfill$\blacksquare$

\begin{lemma}\label{lemma:Gamma}
Let the assumptions of Theorem \ref{thm:main} hold. Consider the functionals $S_j:(BV(D),\mathcal{T}^{\prime\prime})\ra[0,+\infty]$ defined in (\ref{mainfunctionaln}) and $S_\infty:(BV(D),\mathcal{T}^{\prime\prime})\ra[0,+\infty]$ defined in (\ref{mainfunctional}). Here $\mathcal{T}^{\prime\prime}$ denotes the intermediate topology of Definition \ref{def:intermediate}. 

Then the sequence $S_j$ of functionals $\Gamma$-convergences to $S_\infty$ as $j\ra\infty$. 
\end{lemma}
{\bf Proof.} Note that one can use the distance function
$$
  d^{\prime\prime}(u_1,u_2) := \|u_1-u_2\|_{L^1(D)} + \left| \int_D|\nabla u_1|dx -\int_D|\nabla u_2| dx\right|
$$
to show that $\mathcal{T}^{\prime\prime}$ is metrizable. Because of the metrizability, showing $\Gamma$-convergence can be done by proving two things \cite[Proposition 8.1]{DalMaso1993}. The {\em lower-bound inequality}
\begin{equation}\label{lowerboundineq}
 S_\infty(u)\leq \liminf_{j\ra \infty}S_j(u_j)
\end{equation}
should hold whenever there is intermediate convergence of $u_j\in BV(D)$ to the limit $u\in BV(D)$. Further, we must show that for every $u\in BV(D)$ there is a sequence $u_j\in BV(D)$ converging intermediately to $u$ and satisfying the {\em upper-bound inequality }
\begin{equation}\label{upperboundineq}
 S_\infty(u)\geq \limsup_{j\ra \infty}S_j(u_j).
\end{equation}

\noindent
{\em Proving the upper-bound inequality (\ref{upperboundineq}).} Take any $u\in BV(D)$. Recall the orthogonal projection operators $T_j:L^2(D)\to L^2(D)$ defined in (\ref{def_Tn}), and the spaces $Y_j=\mbox{Range}(T_j)$ of piecewise constant functions. As $Y_j\subset Y_{j+1}$, by Lemma \ref{approxlemma} there are $v_j\in  Y_j$  such that $\displaystyle\lim_{j\to \infty}\rho_{1,1}(v_j,u)=0$. This implies the following:
\begin{eqnarray}
&&\label{unL1conv1}
\lim_{j\to\infty}\|v_j-u\|_{L^1(D)}=0,\\
&&\label{apu limit 1}
  \lim_{j\to \infty}V(v_j)=V(u),\\
&&\label{C0bounded} 
V(v_j) \leq C_0 \mbox{ for all }j=1,2,3,\dots, \mbox{ and for some }C_0<\infty.
\end{eqnarray}
We claim that there is intermediate convergence of a subsequence of $v_j$ to the limit $u$. The condition (\ref{def:intermediate1}) was already established in (\ref{unL1conv1}). Regarding the  condition  (\ref{def:intermediate2}), note that by (\ref{C0bounded}) the total variation of $v_j$ is bounded: $|\nabla v_j|\leq C_0$, and we can view the measures $|\nabla v_j|$ as a bounded sequence in the dual space $(C(D))^*$. Thus, by the Banach-Alaoglu theorem, there is a convergent subsequence $\nabla u_{j(\ell)}$ with $\ell=1,2,3,\dots,$ so that $\lim_{\ell\ra \infty}\nabla u_{j(\ell)}= \eta\in(C(D))^*$. Now since $v_j$ tends to $u$ in $L^1(D)$ as $j\ra\infty$, we have that $\nabla u_{j(\ell)}$ tends to $\nabla u$ in the sense of distributions as $\ell\ra\infty$. Therefore, we must have $\eta=\nabla u$, and we have established intermediate convergence of $u_{j(\ell)}$ to  $u$ as $\ell\to\infty$.

Since $Y_j\subset Y_{j+1}$, we can construct a sequence $u_j\in Y_j$ for $j=1,2,3,\dots,$ by setting $u_m=v_{j(\ell)}$ for all $j(\ell)\leq m<j(\ell+1)$. Now there is intermediate convergence of $u_j$ to $u$ as $j\ra \infty$. Thus we have $\lim_{j\ra\infty} \|u_j\|_{L^1(D)}= \|u\|_{L^1(D)}$ and $\lim_{j\ra\infty} V(u_j)= V(u)$. Furthermore, since intermediate convergence implies weak convergence, Lemma \ref{lemma:Aterm} shows that $\lim_{j\to \infty} \Areal u_j=\Areal u$ in $L^2(\Omega)$.
Now the fact that $u_j\in  Y_j$ and formula (\ref{mainfunctionaln}) together imply $S_j(u_j)=S_\infty(u_j)$. Therefore  the upper-bound inequality (\ref{upperboundineq}) follows from
\begin{eqnarray}
&& \nonumber
\|\Areal u-m\|_{L^2(\Omega)}^2+\a_1 \|u\|_{L^1(D)}+\a V(u)\\
&=& \nonumber
\lim_{j\ra \infty}\left(\|\Areal u_j-m\|_{L^2(\Omega)}^2+\a_1 \|u_j\|_{L^1(D)}+\a V(u_j)\right).
\end{eqnarray}

{\em Proving the lower-bound inequality (\ref{lowerboundineq}).} Assume that $u_j$ converges intermediately to $u$ in $BV$ as $j\ra \infty$. Then by definition we have $\lim_{j\ra \infty}\|u_j\|_{L^1(D)}=\|u\|_{L^1(D)}$, and (\ref{Areal_conv}) shows that $\lim_{j\ra \infty}\|\Areal u_j-m\|_{L^2(\Omega)}=\|\Areal u-m\|_{L^2(\Omega)}$. Moreover, intermediate convergence implies $V(u)\leq \liminf_{j\to \infty}V(u_j),$ see \cite[Proposition 10.1.2]{Attouch2014}.
Now estimate
\begin{eqnarray}
S_\infty(u)
&=& \nonumber
\|\Areal u-m\|_{L^2(\Omega)}^2+\a_1 \|u\|_{L^1(D)}+\a V(u)\\
&\leq& \nonumber%
\lim_{j\ra \infty}\|\Areal u_j-m\|_{L^2(\Omega)}^2+\a_1 \lim_{j\ra \infty} \|u_j\|_{L^1(D)}+\a \liminf_{j\ra \infty}V(u_j)\\
&=&\nonumber
\liminf_{j\ra \infty}\|\Areal u_j-m\|_{L^2(\Omega)}^2+\a_1 \liminf_{j\ra \infty} \|u_j\|_{L^1(D)}+\a \liminf_{j\ra \infty}V(u_j)\\
&\leq&\nonumber
\liminf_{j\ra \infty}\left(\|\Areal u_j-m\|_{L^2(\Omega)}^2+\a_1 \|u_j\|_{L^1(D)}+\a V(u_j)\right)\\
&=&\nonumber
\liminf_{j\ra \infty} S_\infty(u_j)\\
&\leq&\label{lowerboundineq_expanded}
\liminf_{j\ra \infty} S_j(u_j).
\end{eqnarray}
\hfill$\blacksquare$

\bigskip

\noindent
{\bf Proof of Theorem \ref{thm:main}.} 
Let us  restrict $S_j$ and $S_\infty$ to the set 
$$
  \mathcal B_R=\{u\in BV(D):\  (\alpha_1\|u\|_{L^1(D)}+\alpha V(u))\leq R\}
$$ 
and call the resulting functionals $S_j^{(R)}$ and $S_\infty^{(R)}$, respectively.  
Denote by $\mathcal{T}^\prime$  the topology in $\mathcal B_R\subset BV(D)$ induced by
the weak topology of $BV(D)$. Since $C(D)$ is separable, the space $(B_R,\mathcal{T}^\prime)$ is metrizable and, in particular, satisfies the 1st axiom of countability as a topological space.

Next we show that the restricted  functionals $S_j^{(R)}$  $\Gamma$-converge to the restricted limit functional $S_\infty^{(R)}$. When $u_j\in \mathcal B_R$  converges in the weak topology to the limit $u\in \mathcal B_R$,
the lower-bound inequality
(\ref{lowerboundineq}) follows from Lemma \ref{lemma:Aterm}.
 Furthermore, Lemma \ref{lemma:Gamma} yields that for all $u\in \mathcal B_R$ there is
  some sequence $u_j$ converging to $u$ in the intermediate topology of $BV$ for which the 
   the upper-bound inequality (\ref{upperboundineq}) holds. Then,
   we see that  the sequence
   \begin{equation}\label{discTrick}
      v_j(x):=
      \frac { V(u)+ \|u\|_{L^1(D)}} {V(u_j)+\|u_j\|_{L^1(D)}+\frac 1j} u_j(x)
   \end{equation}
satisfies $v_j\in \mathcal B_R$, $v_j$ converges to $u$ in the intermediate topology of $BV$,
 and that the upper-bound inequality (\ref{upperboundineq}) is
valid also for the sequence $v_j$ and $u$. 
As the intermediate convergence in $BV$ implies the weak convergence in BV, these facts
show that $S_j^{(R)}$  $\Gamma$-converge to $S_\infty^{(R)}$, in $(\mathcal B_R,\mathcal{T}^\prime)$,  as $j\ra \infty$ for any fixed $R>0$,
see e.g.\ \cite{Attouch1989,Attouch2014}.

%

Let us prove the existence of minimizers for $S_\infty^{(R)}$ and $S_j^{(R)}$. Note first that
\begin{itemize}
\item by Lemma \ref{lemma:Aterm} the functionals $S_\infty$ and $S_j$ are lower semicontinuous with respect to the weak topology of $BV(D)$, and 
\item by the Banach-Alaoglu theorem and \cite[Theorem 10.1.4]{Attouch2014} the space $(\mathcal B_R,\mathcal{T}^\prime)$ is compact.
\end{itemize}
Therefore, the functionals $S_\infty$ and $S_j$ have minimizers for any $R>0$.

Set $R=s_\infty+1$ with $s_\infty$ defined in (\ref{def:infS}). Then there exists a minimizer $\uMAP_\infty\in \argmin(S_\infty^{(R)})$ belonging to the open interior of   $\mathcal B_R$. 

Using (\ref{discTrick}) we see that there are $ u_j\in Y_j\cap \mathcal B_R$ for which  $u_j\ra \uMAP_\infty$ in the intermediate topology and (\ref{lowerboundineq_expanded}) holds. For all large enough $j$ we have 
\ba
S_\infty(u_j)\leq S_\infty(\uMAP_\infty)+1=s_\infty+1.
\ea
As $S_j(u_j)=S_\infty(u_j)$ and $S_m(u)\leq S_j(u)$  for $m\geq j$,
we see that for all large enough $j$ the  minimizers of the functionals $S_j^{(R)}$ coincide with 
the minimizers of $S_j$. Therefore we do not lose any generality in our minimization analysis by restricting to  $\mathcal B_R$. 

{\mnewtext Let us now study the convergence properties of minimizers. As we saw above, for large enough $j$ all minimizers $\uMAP_j\in \mbox{argmin}(S_j) =\mbox{argmin}(S_j^{(R)})$ 
and $\uMAP_\infty\in \mbox{argmin}(S_\infty)=\mbox{argmin}(S_\infty^{(R)})$ 
are in the ball $\mathcal B_R$.
The compactness of $(\mathcal B_R,\mathcal{T}^\prime)$ implies that  any sequence $(\uMAP_j)_{j=1}^\infty$ of minimizers has
a subsequence $(\uMAP_{j(\ell)})_{\ell=1}^\infty$ converging to some $w\in \mathcal B_R$ in the weak topology of $BV$.
By \cite[Cor.\ 7.20]{DalMaso1993}, the limit $w$ is then a minimizer of $S_\infty^{(R)}$ and
we have
\begin{equation}\label{M1}
  \lim_{\ell\to \infty} S_{j(\ell)}(\uMAP_{j(\ell)})=   \lim_{\ell\to \infty} S_{j(\ell)}^{(R)}(\uMAP_{j(\ell)})=S_\infty^{(R)}(w)= S_\infty(w)=s_\infty.
\end{equation}
Moreover, by \cite[Cor.\ 7.20]{DalMaso1993},
%
the equation (\ref{M1}) holds for an arbitrary  sequence $u_{j(\ell)}\in \mbox{argmin}(S_{j(\ell)})$, $\ell=1,2,\dots$  that converge in the weak topology of $BV$ to some limit $w$, and $w$ is a minimizer of $S_\infty$.


Let us next consider a subsequence $\uMAP_{j(\ell)}\in \mbox{argmin}(S_{j(\ell)})
$ converging to $\uMAP_\infty\in \mbox{argmin}(S_\infty)
$ in the weak topology
of $BV$.
Then we have 
\beq\label{inequality epi limit BV}
V(\uMAP_\infty) \leq \liminf_{\ell\to \infty}V(\uMAP_{j(\ell)}).
\eeq
Assume that we would have here a strict inequality.

By definition we have for any subsequence $\uMAP_{j(\ell)}$ converging to $\uMAP_\infty$ in the weak topology of $BV$ 
\ba
\|\uMAP_\infty\|_{L^1(D)} = \lim_{\ell\to \infty}\|\uMAP_{j(\ell)}\|_{L^1(D)}.
\ea
{\mnewtext By Lemma \ref{lemma:Aterm}, we have}
\beq\label{epi limit likelihood}
\|\Areal\uMAP_\infty-m \|_{L^2(\Omega)}^2 = \lim_{\ell\to \infty}\|\Areal\uMAP_{j(\ell)}-m \|_{L^2(\Omega)}^2 .
\eeq
These and a strict inequality in (\ref{inequality epi limit BV})   would imply that  $\liminf_{\ell\to \infty} S_{j(\ell)}(u_{j(\ell)})<s_\infty$, but that is not possible in view of (\ref{M1}). Therefore,
\beq\label{epi limit BV}
V(\uMAP_\infty) = \liminf_{\ell\to \infty}V(\uMAP_{j(\ell)}).
\eeq
Similar arguments with $\liminf$ replaced by $\limsup$ yield the desired equality. }

\phantom{m}\hfill$\blacksquare$

\section{Computational models and algorithms}\label{sec:compmodels}

\subsection{X-ray attenuation model}
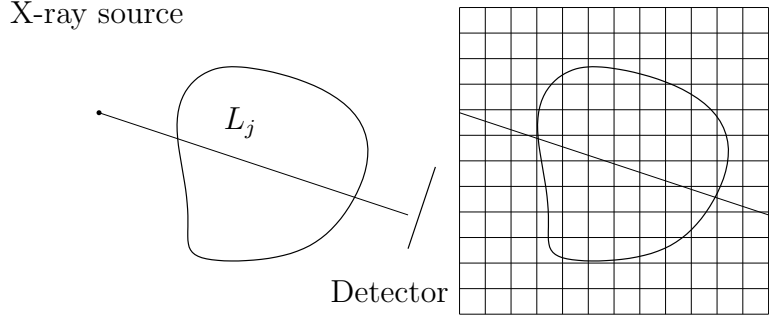
\begin{figure}[t]
\begin{picture}(150,170)
\put(-50,0){\includegraphics[width=9cm]{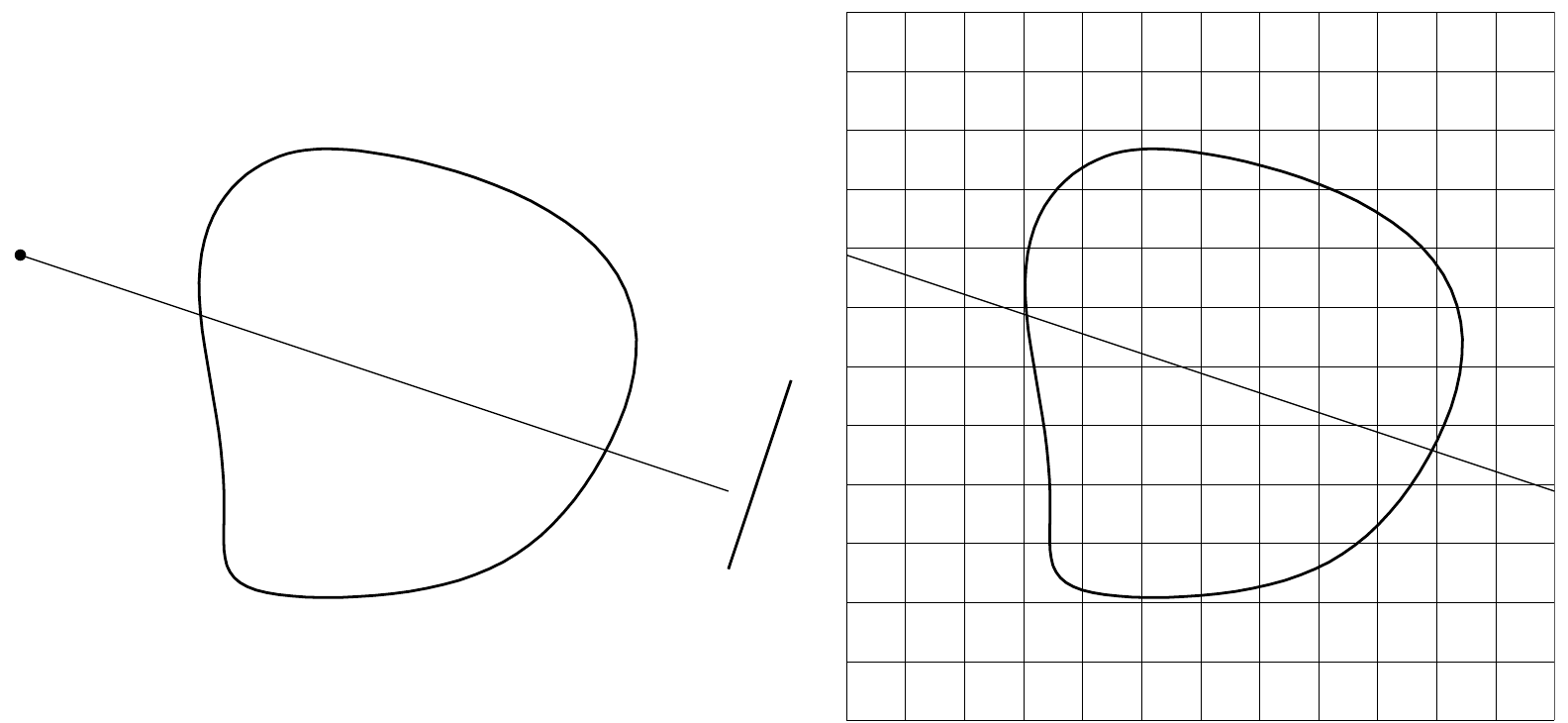}}
\put(-80,110){X-ray source}
\put(0,70){$L_j$}
\put(40,5){Detector}
\end{picture}
\caption{\label{fig:XrayIllustarion} Left: Illustration of the
  pencil-beam attenuation model of X-ray tomography. Right:
  Domain $\Omega$ discretized into a lattice of pixels/voxels $\Omega_i$}
\end{figure}

In X-ray imaging the radiation source is placed on one side of an object
and the detector at the opposite side. The rays pass through
the object and the attenuated signal is detected by a digital sensor
(array of almost point-like detectors), see Figure
\ref{fig:XrayIllustarion} for an illustration. We
model the (2D slice of the) object by a rectangle $\Omega \subset \mathbb{R}^{2}$ and
by a non-negative attenuation coefficient $\atten: \Omega \rightarrow
[0,\infty )$. The attenuation coefficient $\atten$ gives the relative
intensity loss of an X-ray travelling a small distance d$s$ at $s \in
\Omega$. This leads to the following linear model
\begin{equation}
  \label{eq:XrayMeasModel_1}
  g_{j} = -\log \left( \frac{I_j}{I_0}\right) = \int_{L_j} \atten(s)\mbox{d}s,
\end{equation}
where $L_j$ is the $j$th line of the X-ray, $g_{j}$ is the value of the projection measurement of the $j$th
source to detector line $L_j$, $I_j$ is the measured X-ray intensity
and $I_0$ is the initial intensity of the X-ray beam before entering
$\Omega$.

We discretize this model by dividing the domain $\Omega$ into a
lattice of pixels $\Omega_i$ and by computing the length of the
path $L_j$ inside each pixel $\Omega_i$, see figure
\ref{fig:XrayIllustarion}. Assuming that the attenuation function
$\atten(s)$ is constant inside each pixel $\Omega_i$, the projection
measurements $g_j$ can be approximated by
\begin{equation}
  \label{eq:XrayMeasModel_2}
  g_j = \int_{L_j}\atten(s)\mbox{d}s \approx \sum_{i} f_i \vert \Omega_i
  \cap L_j \vert,
\end{equation}
where $f_i$ is the constant value of the function $u$ in the $i$th pixel and $\vert \Omega_i  \cap L_j \vert$ denotes the length of the ray
$L_j$ through pixel $\Omega_i$. Arranging the set of $M$ measurements
into a vector $\g = (g_1,g_2,\ldots,g_M) \in \mathbb{R}^{M}$, we obtain
the model 
\begin{displaymath}
  \g = A\f,
\end{displaymath}
where $\f  \in \mathbb{R}^{N}=\mathbb{R}^{n\times n}$ is the vector
containing the attenuation values in the pixels and $A \in
\mathbb{R}^{M \times N}$ is the matrix that implements the transform
from the pixel values to the projection data, see equation (\ref{eq:XrayMeasModel_2}).

\subsection{Discrete anisotropic TV}\label{sec:DiscTV}

We consider the anisotropic total variation
\begin{equation*}
  \Vert \f \Vert _{\rm TV} = \Vert D_{H} \f \Vert_{1} + \Vert
  D_{V}\f \Vert _{1}
\end{equation*}
where the $N\times N$ matrices $D_{H}$ and $D_{V}$ implement horizontal and vertical differences in the $n\times n$ pixel image represented by vector $\f$. We use periodic boundary conditions. See Appendix \ref{Appendix_simple} for an example of the construction of $D_{H}$ and $D_{V}$. 

Now the TV regularized X-ray imaging problem (\ref{def:discreteTV}) can be rewritten as follows:
\begin{equation}
  \label{eq:TVregularizedXrayPb}
 \argmin_{\f\in \R^N_+} \left \{ \frac{1}{2} \Vert A \f - \greal \Vert^{2}_{2} +
    \alpha \left ( \Vert D_{H}\f \Vert_{1} + \Vert D_{V} \f \Vert_{1}
    \right ) \right \}. 
\end{equation}

\subsection{Quadratic programming}

The minimization problem of (\ref{eq:TVregularizedXrayPb}) can be
reformulated as a quadratic programming problem. Denoting $D_{H}
\f = \hh^{+} - \hh^{-}$ and $D_{V}\f = \vv^{+} - \vv^{-}$
(\ref{eq:TVregularizedXrayPb}) can be rewritten as follows
\begin{equation*}
  \argmin \left \{ \frac{1}{2} \f^{T}A^{T}A\f - \f^{T}A^{T} \greal + \alpha
    {\bf 1}^{T} \hh^{+} + \alpha {\bf 1}^{T} \hh^{-} + \alpha {\bf 1}^{T}
    \vv^{+} + \alpha {\bf 1}^{T} \vv^{-}  + \frac{1}{2} \greal^{T}\greal \right \},
\end{equation*}
where ${\bf 1}$ is a vector of all ones and
$\hh^{+},\hh^{-},\vv^{+},\vv^{-}\geq 0$. Further we can denote
$$ \zz = \left[ \begin{array}{c}\f\\ \hh^{+} \\ \hh^{-} \\ \vv^{+} \\ \vv^{-}
  \end{array}\right ],\ \  Q = \left[ \begin{array}{ccccc} A^{T}A & 0 &
    0&0&0\\0&0&0&0&0\\0&0&0&0&0\\0&0&0&0&0\\0&0&0&0&0
  \end{array}\right ],\ \   \cc = \left[ \begin{array}{c} -A^{T} \greal \\ \alpha
  {\bf 1} \\ \alpha  {\bf 1} \\ \alpha  {\bf 1} \\ \alpha  {\bf
    1}\end{array}\right ], \ \ \dd = \frac{1}{2}\greal^{T}\greal. 
$$ 
Now we can
write
\begin{eqnarray}
  \label{eq:QP}
\nonumber  \min_{\zz}&& \frac{1}{2} \zz^{T}Q \zz + \cc^{T} \zz + \dd\\
\mbox{s.t}&& \mathcal{B}\zz = \bb \\
\nonumber && \zz \geq 0
\end{eqnarray}
where $Q$ is symmetric and positive definite matrix and 
$$\mathcal{B} = \left [ \begin{array}{ccccc} D_{H} & -I & I& O &
    O\\D_{V} &O&O&-I&I \end{array} \right ],$$ $I$ denotes identity
matrix and $O$ denotes a matrix of all zeros.

\subsection{Primal-dual interior-point  (PD-IP) method}

 The formulation of the primal-dual interior-point algorithm starts by
rewriting the primal problem (\ref{eq:QP}) as a logarithmic barrier
problem \cite{Fiacco1968}
\begin{eqnarray}
  \label{eq:LogBar}
  \nonumber  \min_{\zz}&& \frac{1}{2} \zz^{T}Q\zz + \cc^{T}\zz + \dd - \mu \sum \log(\zz)\\
\mbox{s.t}&& \mathcal{B}\zz = \bb .
\end{eqnarray}
To solve this minimization problem, we introduce the Lagrangian function
\begin{equation}
  \label{eq:Lagrangian}
  \mathcal{L}(\zz,\yy;\mu) = \frac{1}{2}\zz^{T}Q\zz + \cc^{T}\zz + \dd - \mu \sum
  \log(\zz) - \yy^{T} (\mathcal{B}\zz -\bb),
\end{equation}
where $\yy$ is the Lagrangian multiplier. We note that the Lagrangian
multiplier $\yy$ is also the dual variable of the associated Lagrangian
dual problem. For further information see
\cite{Wright1997,Nocedal2006}, for example.

Now the minimization problem of (\ref{eq:LogBar}) can be solved by seeking a stationary
point for the Lagrangian function. Differentiating the Lagrangian with
respect to primal ($\zz \in \mathbb{R}^{n_{z}}$) and dual ($\yy \in \mathbb{R}^{n_{y}}$) variables yields
\begin{eqnarray*}
  0=\nabla_{\zz} \mathcal{L}(\zz,\yy;\mu) &=& Q\zz + \cc - \mu Z^{-1}{\bf 1} +
  \mathcal{B}^{T} \yy \\
0 = \nabla_{\yy}\mathcal{L}(\zz,\yy;\mu) &=& \bb - \mathcal{B}.
\end{eqnarray*}
These conditions are the first order necessary optimality conditions,  often referred to as the Karush-Kuhn-Tucker (KKT)
conditions. 

To derive the primal-dual (path-following) method we consider the
perturbed KKT conditions \cite{Nocedal2006,Wright1997}. Denoting $\mu Z^{-1}{\bf 1} =
\widetilde{\xx}$, the perturbed KKT
conditions can be written as a mapping
$\mathcal{F}:\mathbb{R}^{2n_{z}+n_{y}} \rightarrow
\mathbb{R}^{2n_{z}+n_{y}}$
\begin{equation}
  \label{eq:mapping}
  \mathcal{F}(\zz,\widetilde{\xx}, \yy;\mu) = \left [ \begin{array}{c}Q\zz -
      \mathcal{B}^{T}\yy - \widetilde{\xx} + \cc\\ \mathcal{B} \zz - \bb \\Z
      \widetilde{X} {\bf 1} - \mu {\bf 1},
    \end{array} \right ]
\end{equation}
where $Z = {\rm diag}(z_1,z_2,\ldots,z_{n_z})$, $\widetilde{X} =
{\rm diag}(\widetilde{x}_{1}, \widetilde{x}_{2},\ldots,\widetilde{x}_{n_z})$
and $\mu > 0$ is the central path parameter. The central path is
defined by the trajectory
$\mathcal{P}:\mathcal{P}\{\zz_{\mu},\yy_{\mu},\widetilde{\xx}_{\mu}|\mu >
0\}$. As $\mu \rightarrow 0$ the trajectory $\mathcal{P}$ converges to
optimal solution of both the primal and dual problems. Note that at
the optimal point $\mu = 0$.

Applying Newton's method to (\ref{eq:mapping}), we obtain a linear
system of the form 
\begin{equation}
  \label{eq:KKTmatrix}
  \left [ \begin{array}{ccc} -Q & \mathcal{B}^{T} & I \\ \mathcal{B} &
      0&0\\ I&0&\widetilde{X}^{-1}Z \end{array} \right ]\left
    [\begin{array}{c} \Delta \zz \\ \Delta \yy \\ \Delta \widetilde{\xx}
    \end{array} \right ] = \left [ \begin{array}{c} \pp_1 \\ \pp_2 \\
      \pp_3 \end{array}\right ],
\end{equation}
where $I$ is identity matrix, $\pp_1 = \cc + Q\zz - \mathcal{B}^{T}\yy -
\widetilde{\xx}$, $\pp_2 = \bb - \mathcal{B}\zz$ and $\pp_3 = \mu/ \widetilde{\xx} - \zz -\Delta \widetilde{\xx} \Delta \zz /  \widetilde{\xx}$. 
The variable $\Delta \widetilde{\xx}$ can be
removed without producing any off-diagonal entries in the remaining
system. Hence the KKT conditions can be written in a more compact way
as
\begin{equation}
  \label{eq:KKTreduced}
    \left [ \begin{array}{cc} -(Q+ Z^{-1} \widetilde{X}) & \mathcal{B}^{T} \\ \mathcal{B} &
      0 \end{array} \right ]\left
    [\begin{array}{c} \Delta \zz \\ \Delta \yy \end{array} \right ] =
  \left [ \begin{array}{c} \pp_1 - Z^{-1}\widetilde{X} \pp_3 \\ \pp_2 \end{array}\right ],
\end{equation}
with $\Delta \widetilde{\xx} = Z^{-1}\widetilde{X}(\pp_3 - \Delta \zz)$.

\subsection{The PD-IP algorithm}\label{sec:pd-ip-algorithm}

The PD-IP method presented here is an iterative method and is based on Mehrotra's predictor-corrector approach \cite{Mehrotra1992}. The resulting algorithm can
be summarized as follows.\\

{\bf Algorithm:}
\begin{itemize}
\item []{\bf for} k = 0,1,2,..
\begin{itemize}
\item [] Compute $(\zz^{0},\yy^{0},\widetilde{\xx}^{0})$, with
  $\zz^{0},\widetilde{\xx}^{0}>0$
\item [] Set $\mu = 0, \Delta \zz = 0$ and $\Delta \widetilde{\xx}=0$ and solve (\ref{eq:KKTreduced}) to compute the
  predictor steps $(\Delta \zz^{\rm pre},\Delta \yy^{\rm pre},\Delta
  \widetilde{\xx}^{\rm pre})$.
\item [] Compute a value for $\mu$.
\item [] Compute $\pp_3$ using  $(\Delta \zz^{\rm pre},\Delta \yy^{\rm pre},\Delta
  \widetilde{\xx}^{\rm pre})$ of step 2 and $\mu$ of step 3.
\item [] Solve (\ref{eq:KKTreduced}) for  $(\Delta \zz,\Delta \yy,\Delta
  \widetilde{\xx})$.
\item [] Compute step length multipliers $\lambda_{\rm primal}$
  and $\lambda_{\rm dual}$ 
\item [] Set
  \begin{eqnarray*}
    \zz^{k+1} &=& \zz^{k} + \lambda_{\rm primal} \Delta \zz \\
    \yy^{k+1} &=& \yy^{k} + \lambda_{\rm dual} \Delta \yy \\
    \widetilde{\xx}^{k+1} &=& \widetilde{\xx}^{k} + \lambda_{\rm dual} \Delta \widetilde{\xx} 
  \end{eqnarray*}
\end{itemize}
\item [] {\bf end for}
\end{itemize}

For further details of the algorithm see \cite{niinimaki2013}.

\section{Results}\label{sec:results}

\subsection{X-ray measurement data}

\begin{figure}[t]
\begin{picture}(340,170)
\put(0,0){\includegraphics[width=9cm]{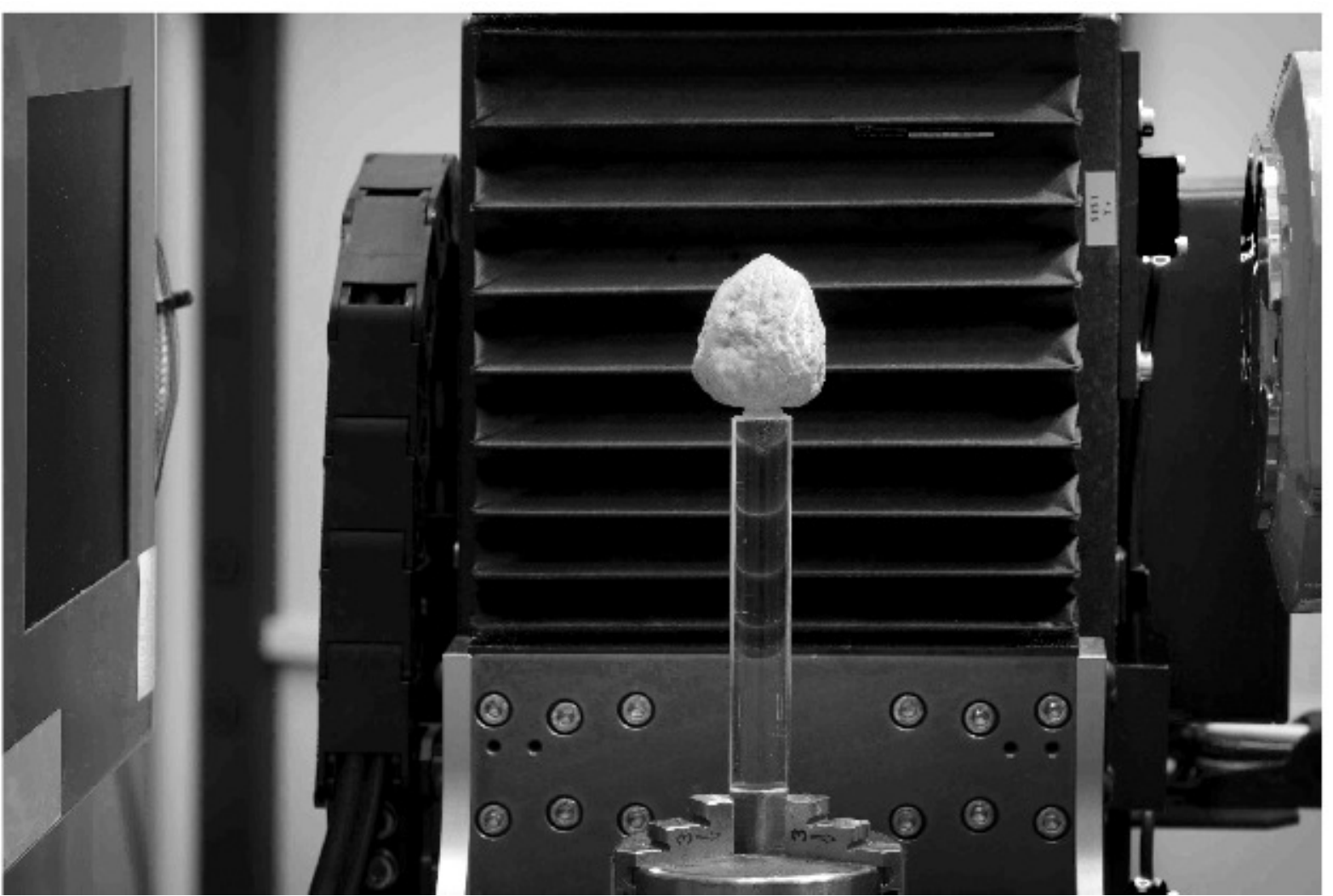}}
\put(275,0){\includegraphics[width=2.8cm]{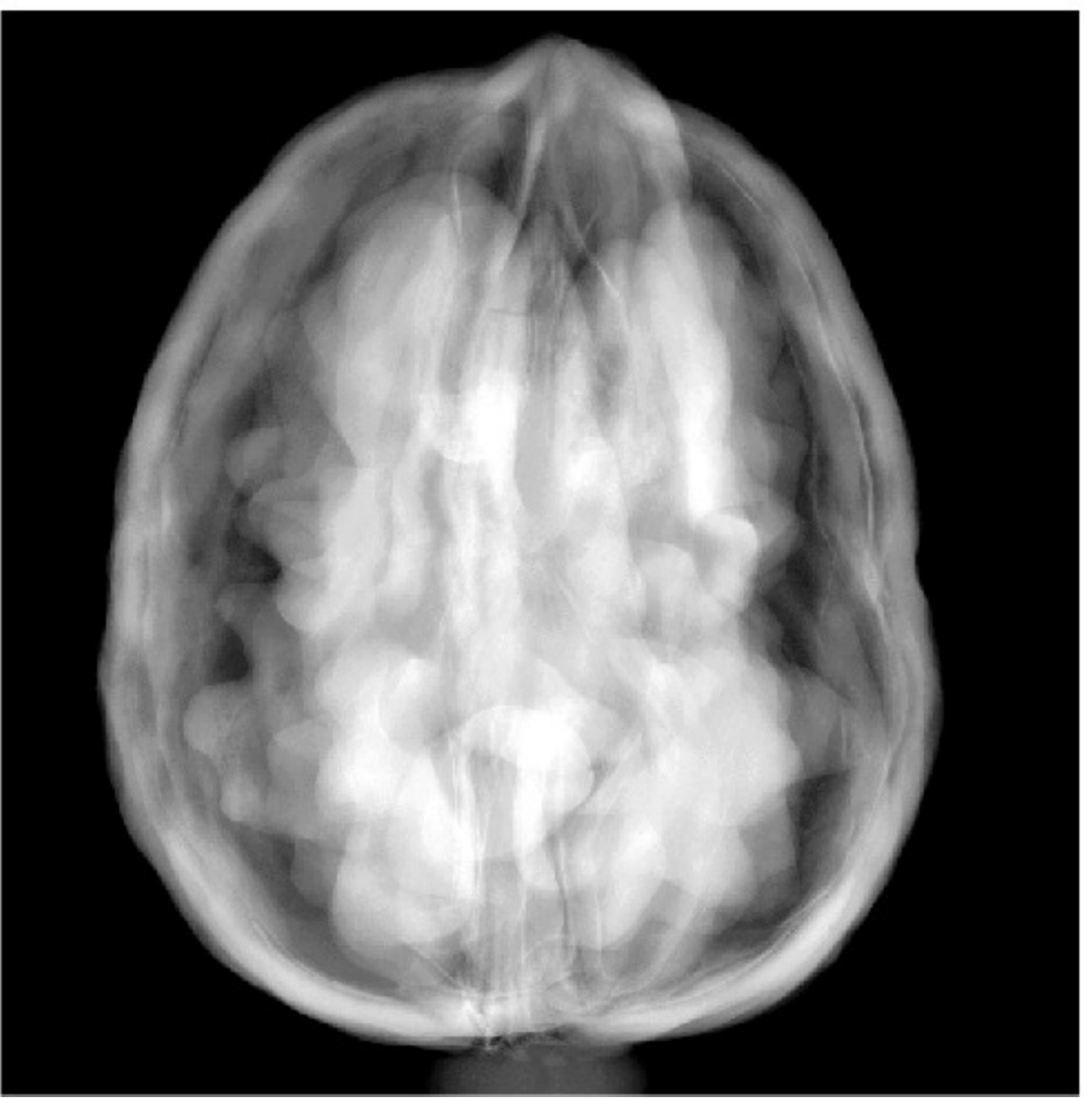}}
\put(275,90){\includegraphics[width=2.8cm]{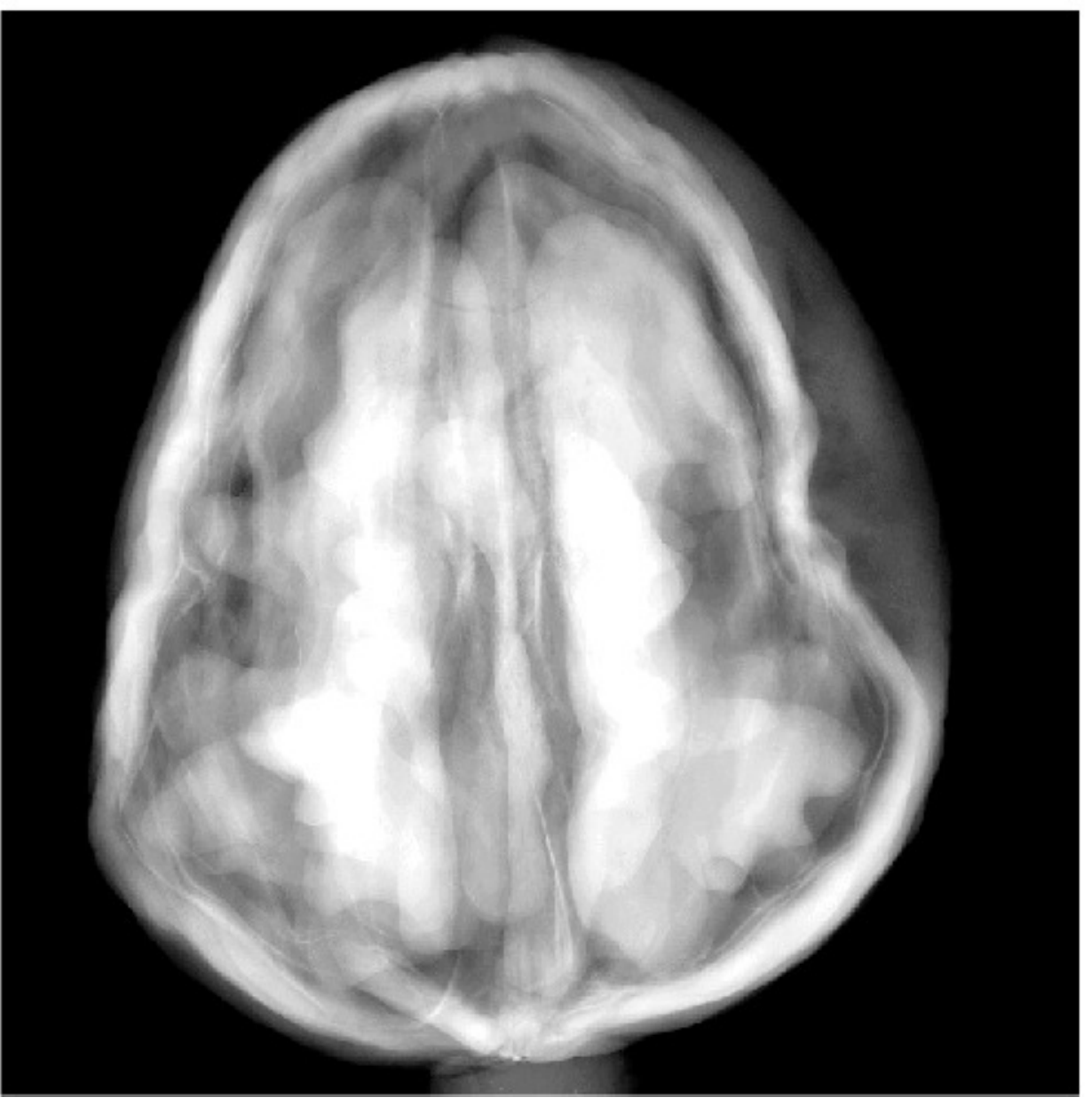}}
\end{picture}
\caption{\label{fig:KeijonLabra}Left: Experimental setup for collecting tomographic X-ray data of a walnut. The detector plane is on the left and the X-ray source on the right in the picture. The walnut is attached to a computer-controlled rotator platform. Right: Two examples of the resulting projection images.}
\end{figure}

The parameter selection method was tested using experimental X-ray tomography data
from a walnut\footnote{The data is made open access at: {\texttt http://www.fips.fi/dataset.php}}.
The projection data of the walnut was acquired with a
 custom-built $\mu$CT device nanotom 180 supplied by Phoenix Xray Systems +  
 Services GmbH (Wunstorf, Germany). The
 measurement set-up is presented in figure \ref{fig:KeijonLabra}. The chosen geometry  
 resulted in a magnification with resolution of 18.33 $\mu$m/pixel. The  
 X-ray detector is a 12-bit CMOS flat panel detector with $2304 \times 2284$  
 pixels of 50 $\mu$m size (Hamamatsu Photonics, Japan). A set of 90 
projection images were acquired over a 180 degree rotation   
with uniform angular step
of 2 degrees between projections. 
Each projection image was composed of an average of six  
 750 ms exposures.
 The
 X-ray tube acceleration voltage was 80 kV and tube current 200 $\mu$A, and  
 the full polychromatic beam was used for image acquisition.
For this work we chose only the projections corresponding to the
middle cross-section of the walnut, thus resulting to a 2D reconstruction
task.

We consider two test cases: 
\begin{itemize}
\item[i)] Reconstruction using the measured X-ray data  
\item[ii)] Reconstruction using the measured X-ray data with 5\% additional
additive random noise, i.e., using measurement $\widetilde{{\bf g}}+{\bf e}$ where
${\bf e}$ is realization from multivariate Gaussian with standard deviation equivalent to
5\% of the maximum value of the actual measurement ${\widetilde {\bf g}}$.
\end{itemize}

\subsection{Selection of $\alpha$ using the parameter choice rule}
\label{secwh2}

In this work we consider three levels of discretization: $n = 128,
192,256$, each resulting to a 2D reconstruction on a $n
\times n$ grid. TV regularized reconstructions $\fan$ were computed with
several values of regularization parameter $\alpha$ by solving the
minimization problem of (\ref{eq:TVregularizedXrayPb}) using the
algorithm of section  \ref{sec:pd-ip-algorithm}.  

\begin{figure}
 \begin{picture}(200,400)
 \put(-40,380){$n = 128$}
 \put(100,380){$n = 192$}
 \put(220,380){$n = 256$}
  \put(-70,260){\includegraphics[width=1.5in]{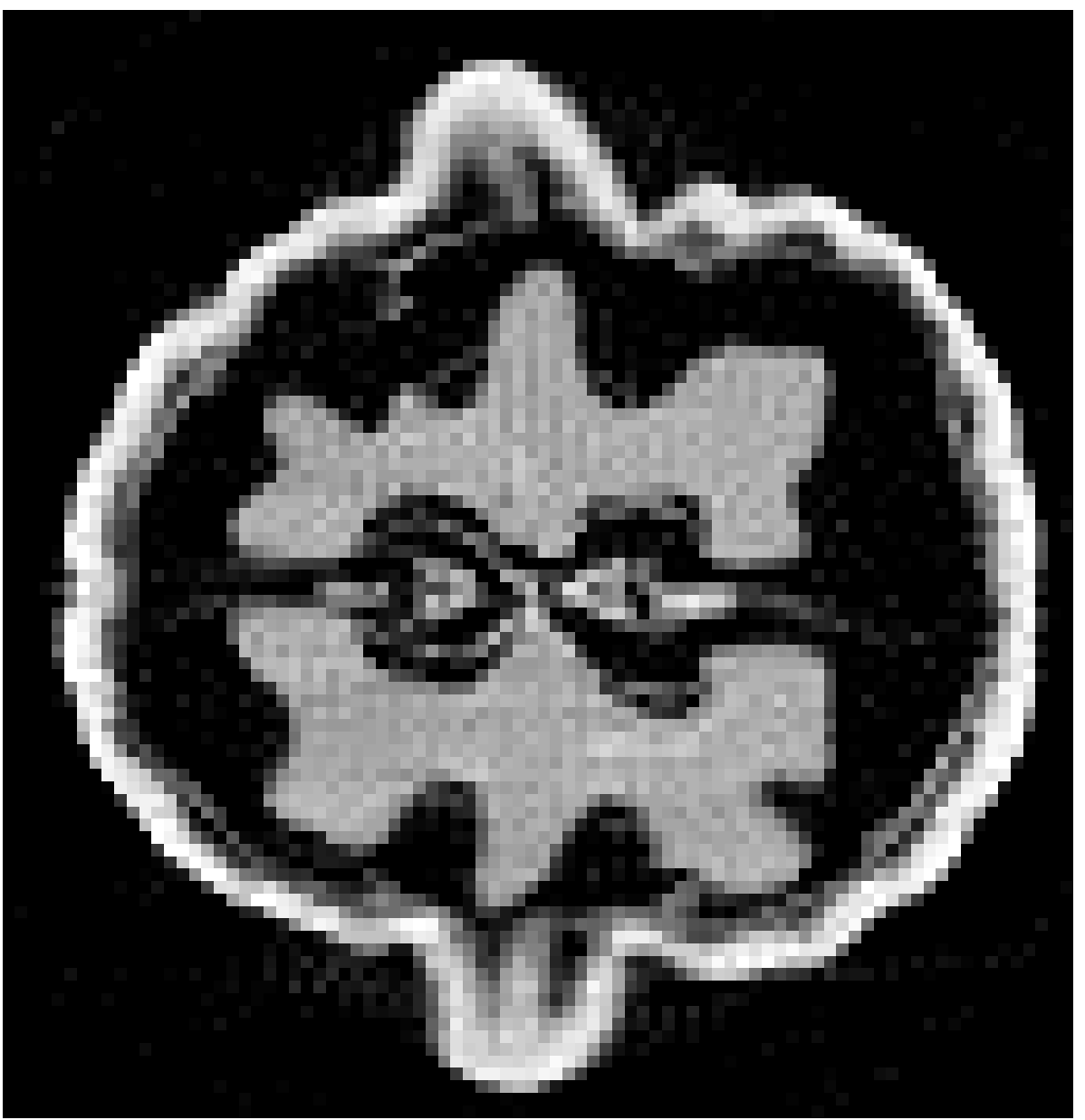}}
  \put(50,255.5){\includegraphics[width=1.5in]{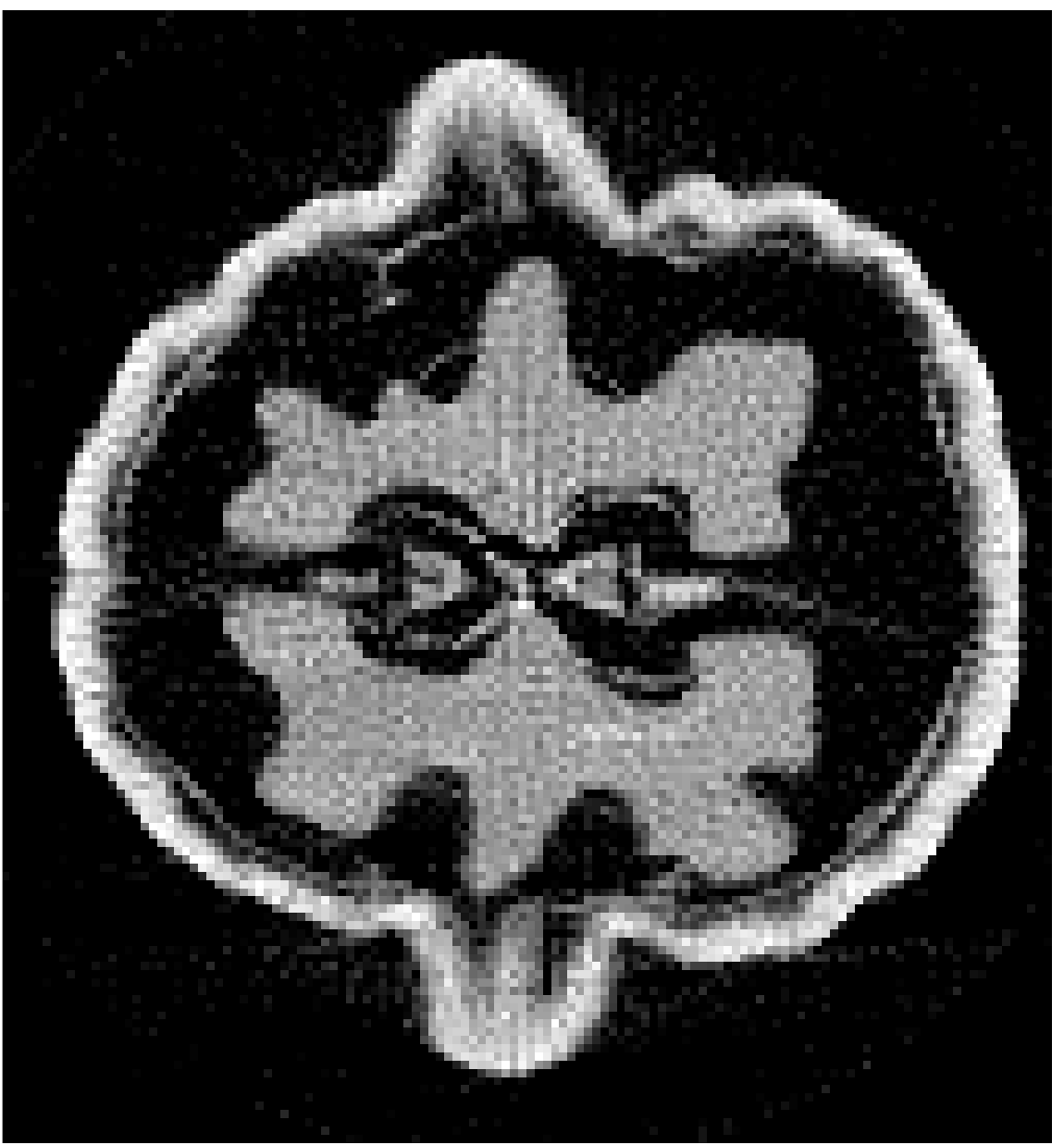}}
  \put(170,260){\includegraphics[width=1.5in]{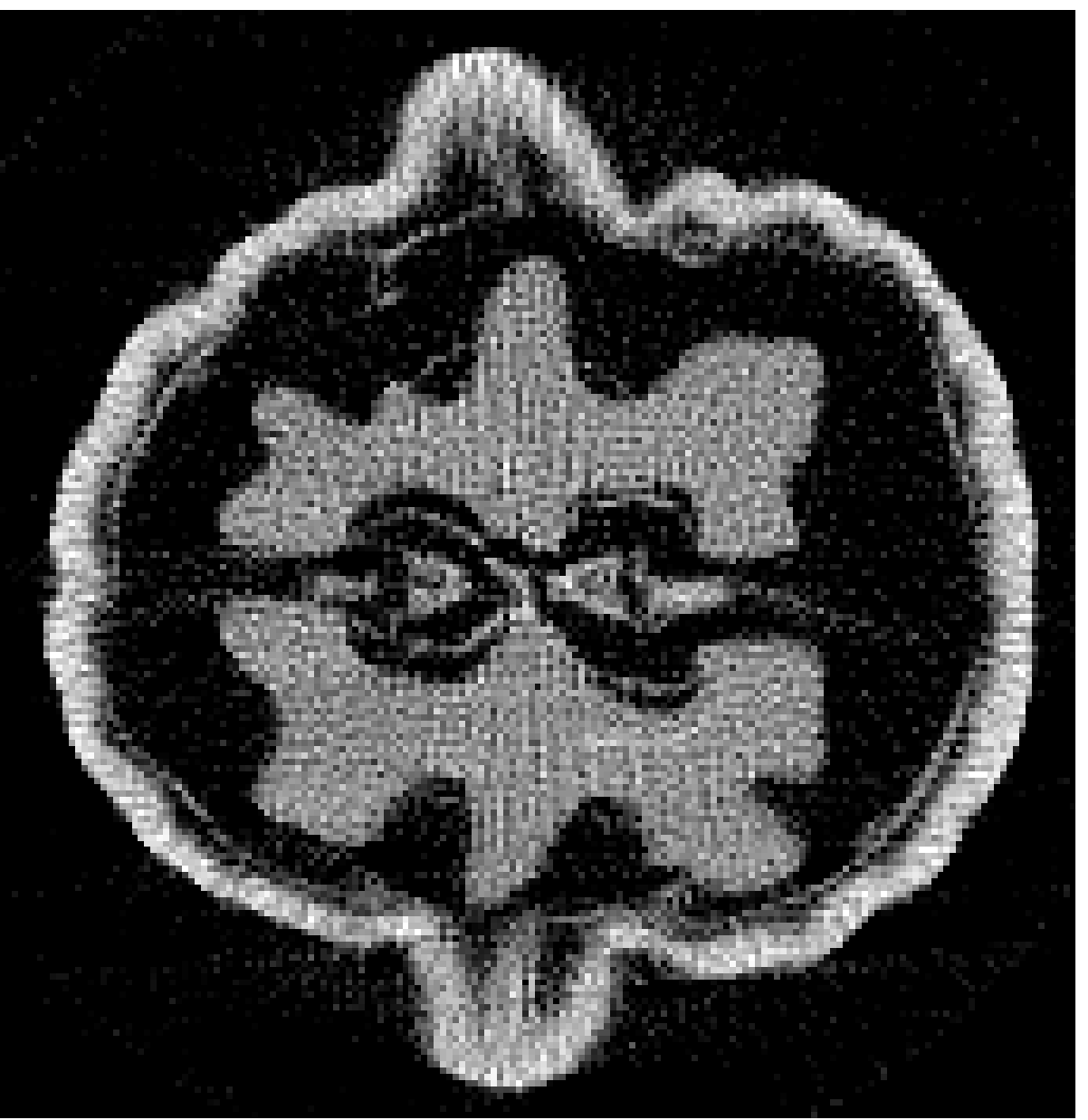}}
  \put(-70,140){\includegraphics[width=1.5in]{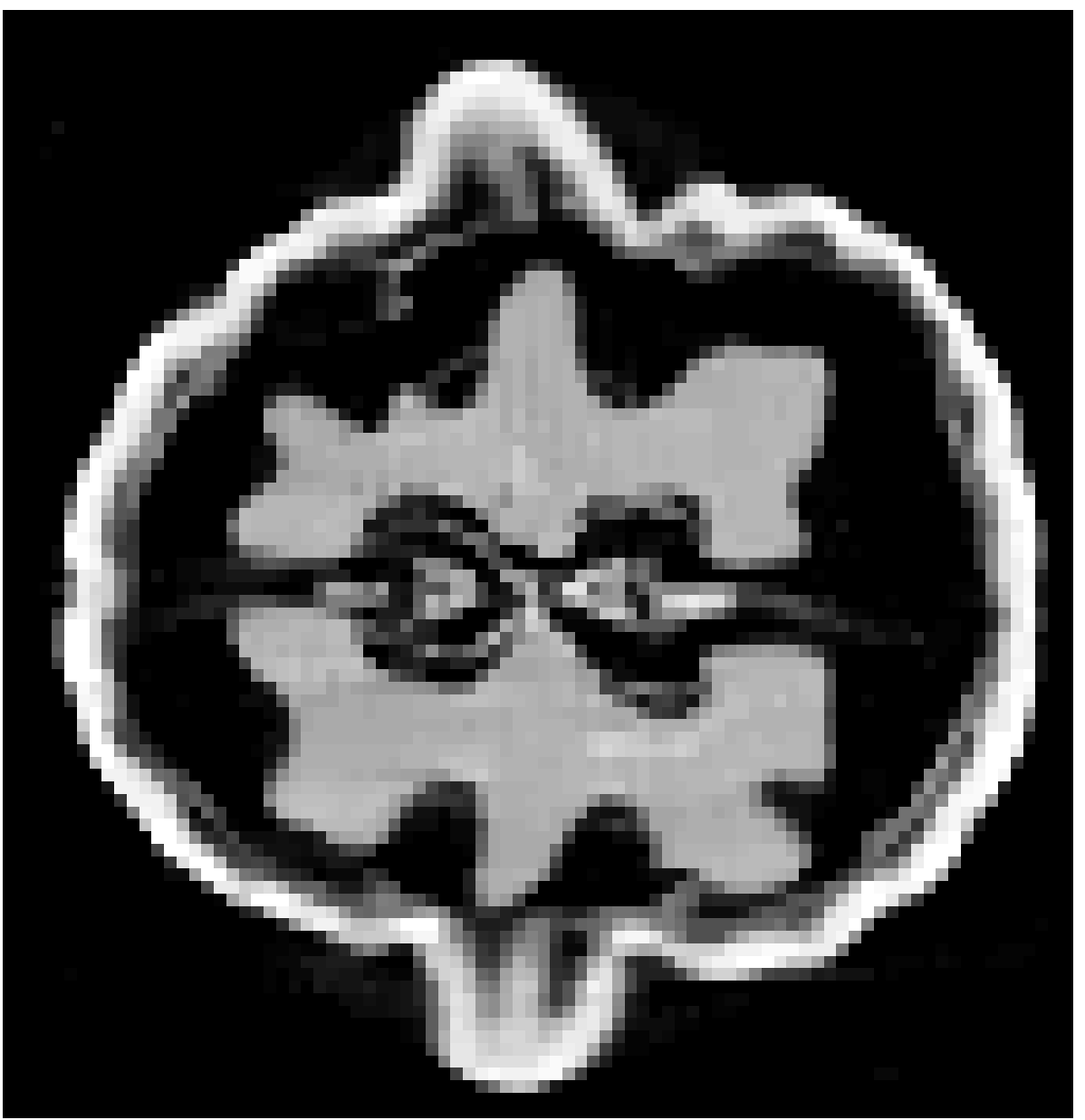}}
  \put(50,135.5){\includegraphics[width=1.5in]{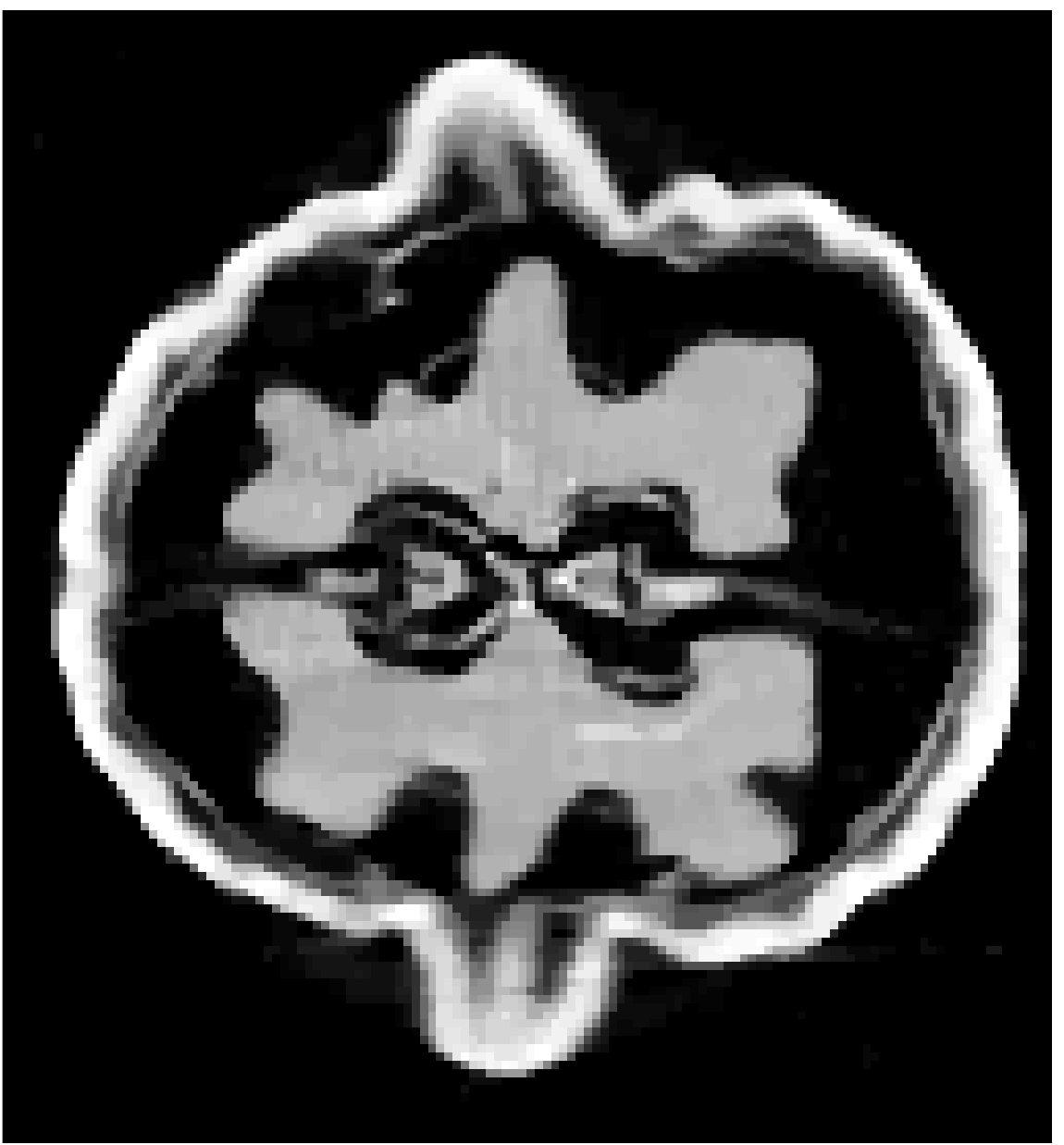}}
  \put(170,140){\includegraphics[width=1.5in]{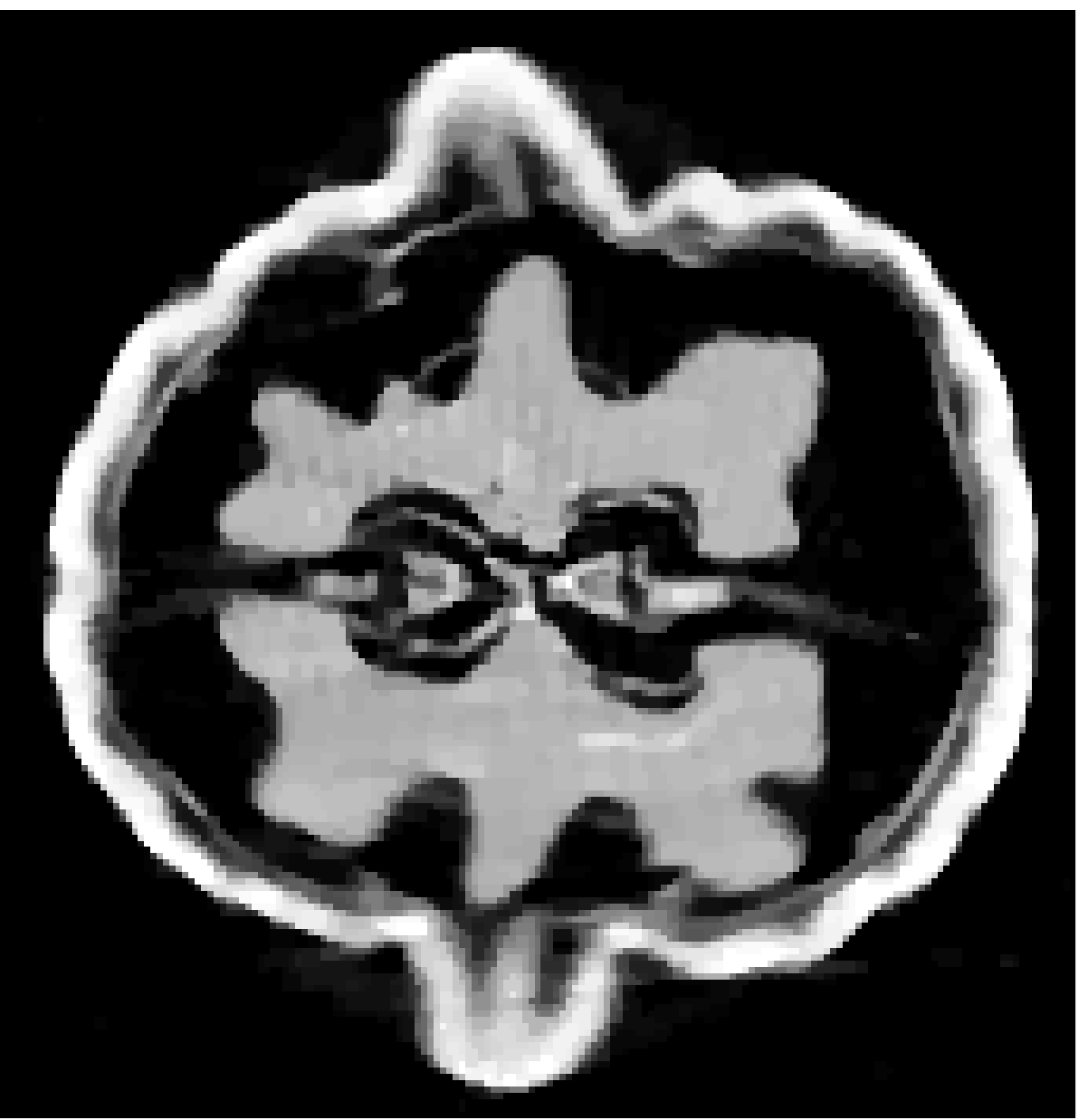}}
  \put(-70,20){\includegraphics[width=1.5in]{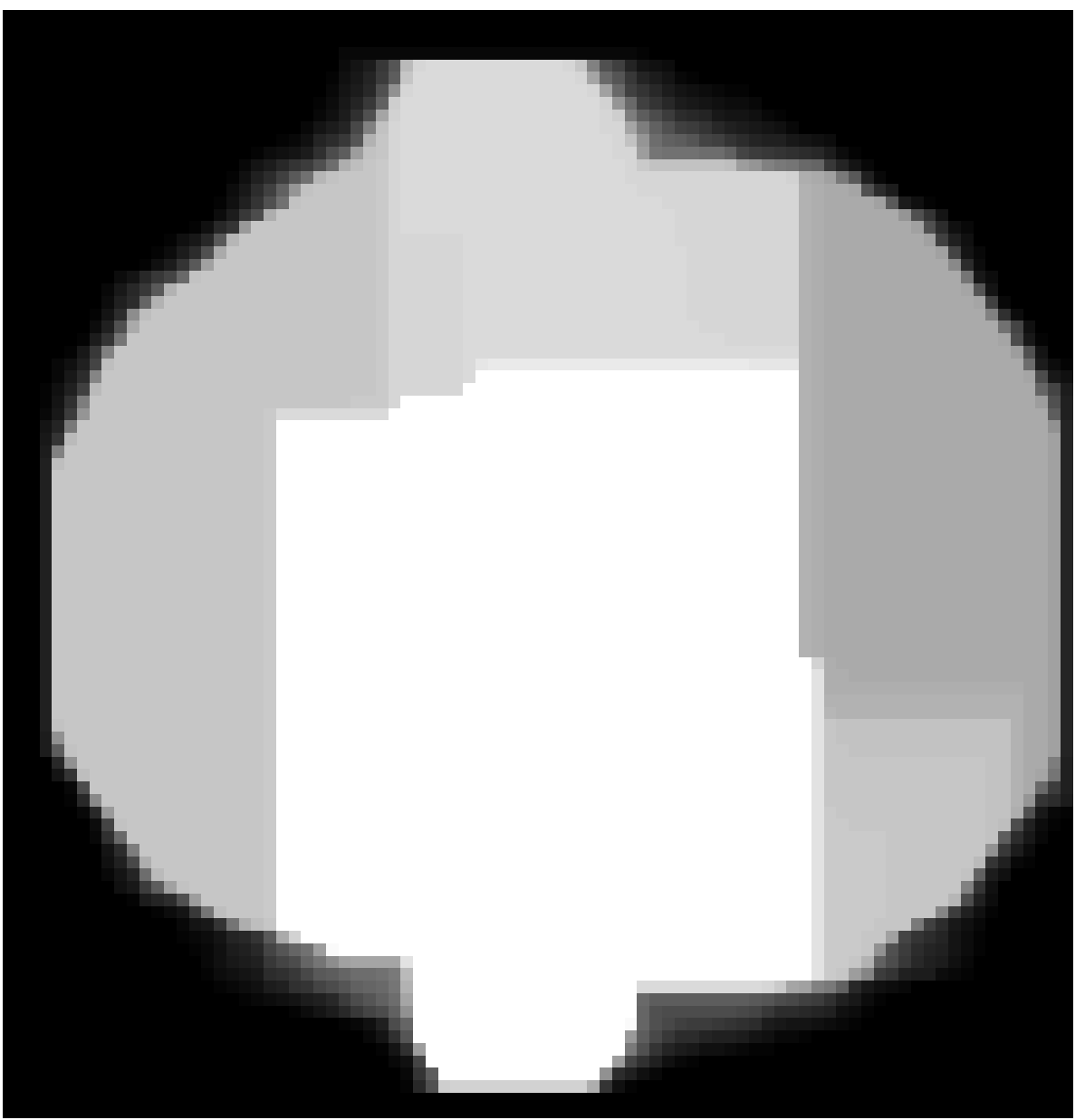}}
  \put(50,15.5){\includegraphics[width=1.5in]{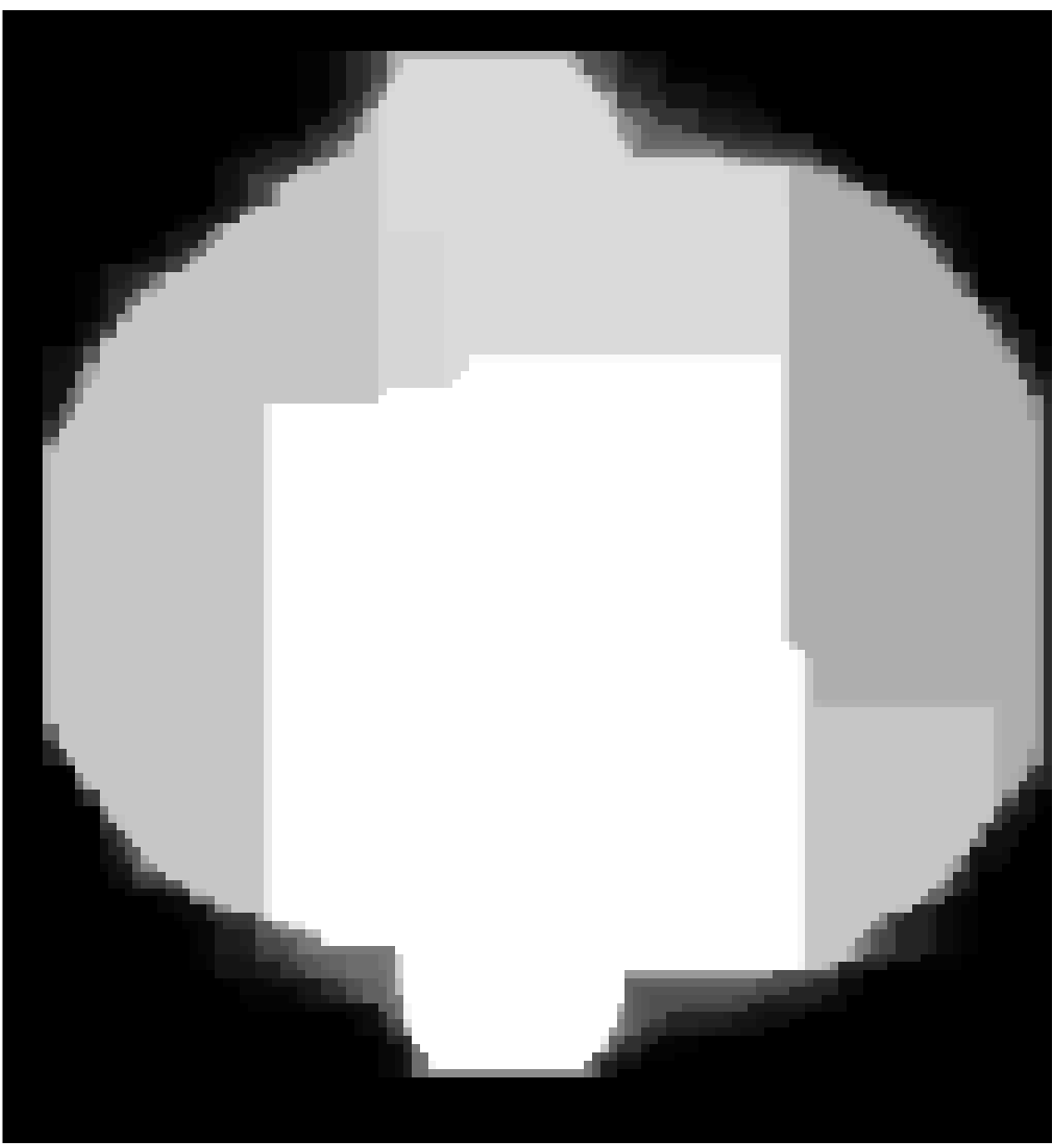}}
  \put(170,20){\includegraphics[width=1.5in]{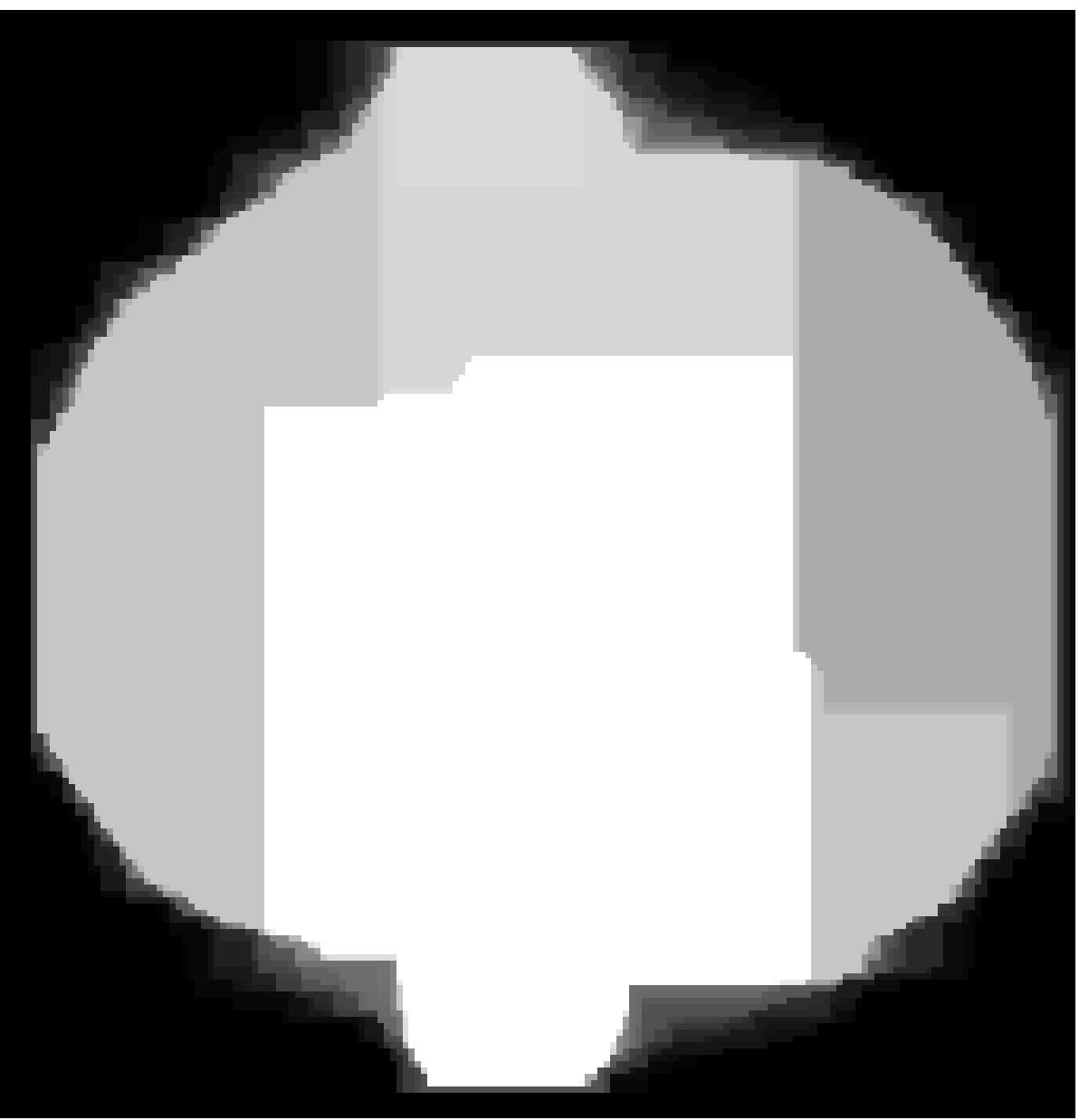}}
 \put(280,310){$\alpha = 10^{-3}$}
 \put(280,200){$\alpha = 10^{0}$}
 \put(280,80){$\alpha = 10^{3}$}
 \end{picture}
 \caption{\label{fig:TVrekot} TV regularized reconstructions $\fan$ from the projection data of the walnut computed with different values
   of regularization parameter $\alpha$. The values of $\alpha$ are
   indicated on the right. Discretization levels $n =
   128,192, 256$ are indicated on the top.}
 \end{figure}

\begin{figure}
 \begin{picture}(200,400)
 \put(-40,380){$n = 128$}
 \put(100,380){$n = 192$}
 \put(220,380){$n = 256$}
  \put(-70,260){\includegraphics[width=1.5in]{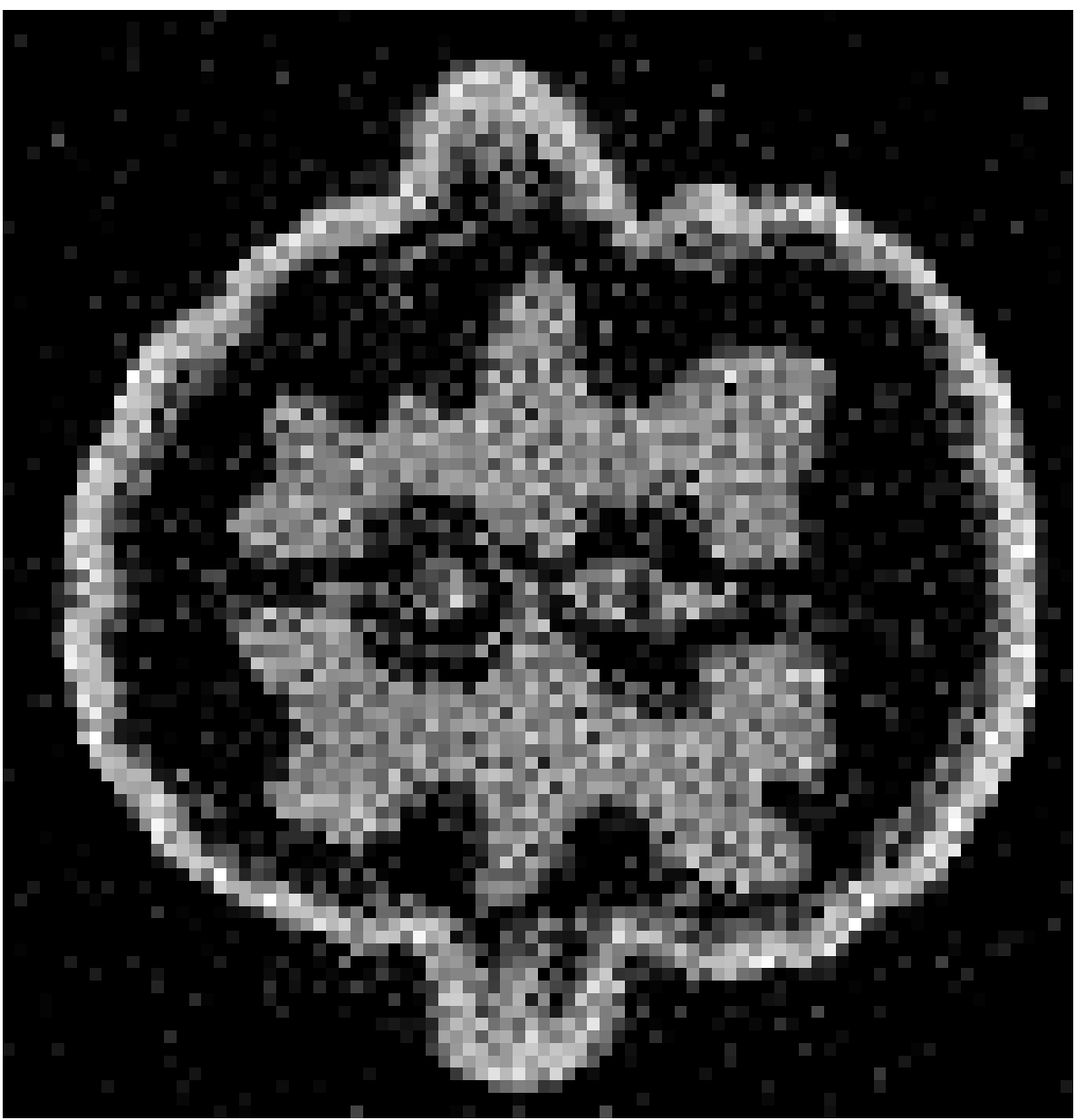}}
  \put(50,255.5){\includegraphics[width=1.5in]{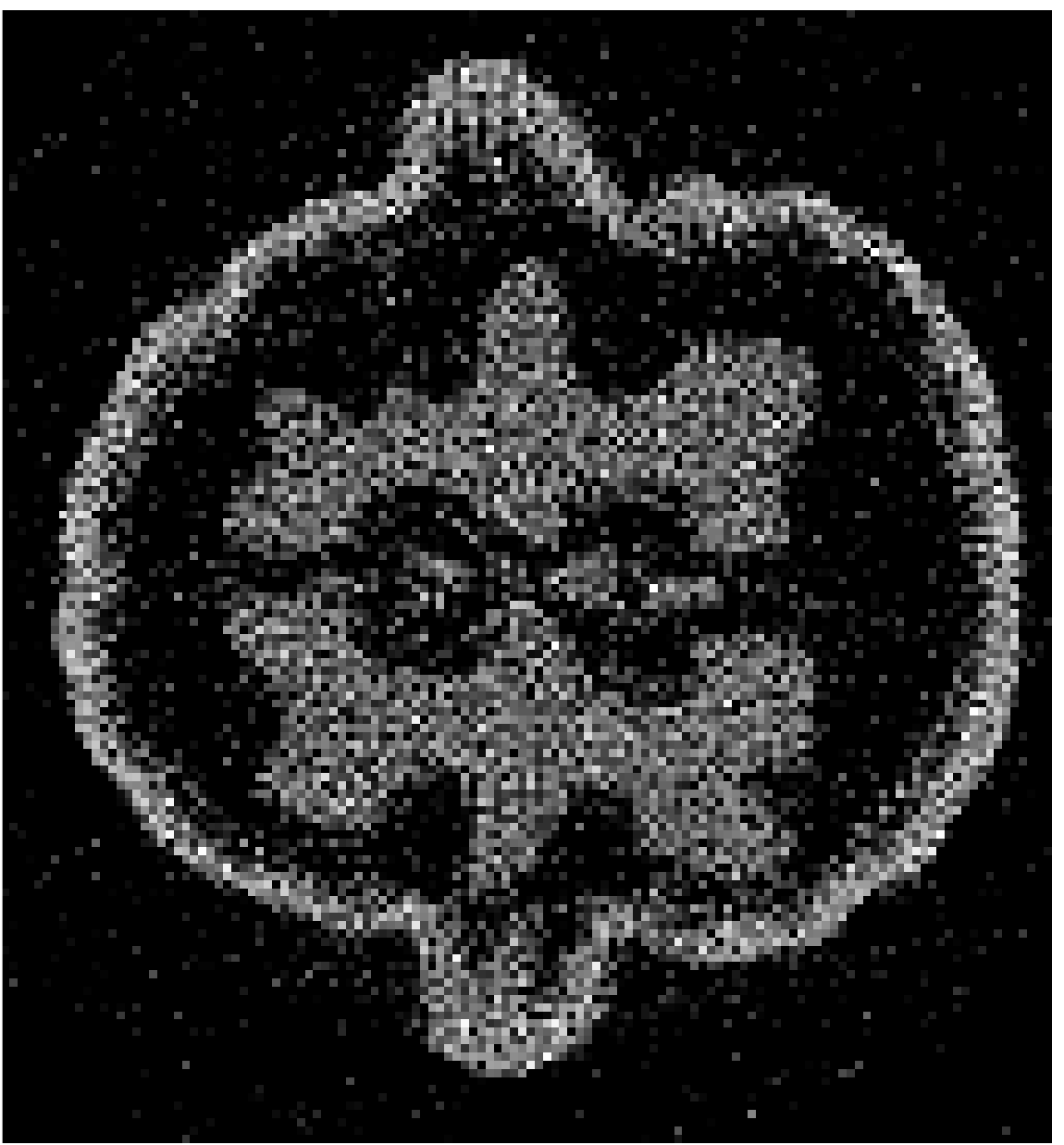}}
  \put(170,260){\includegraphics[width=1.5in]{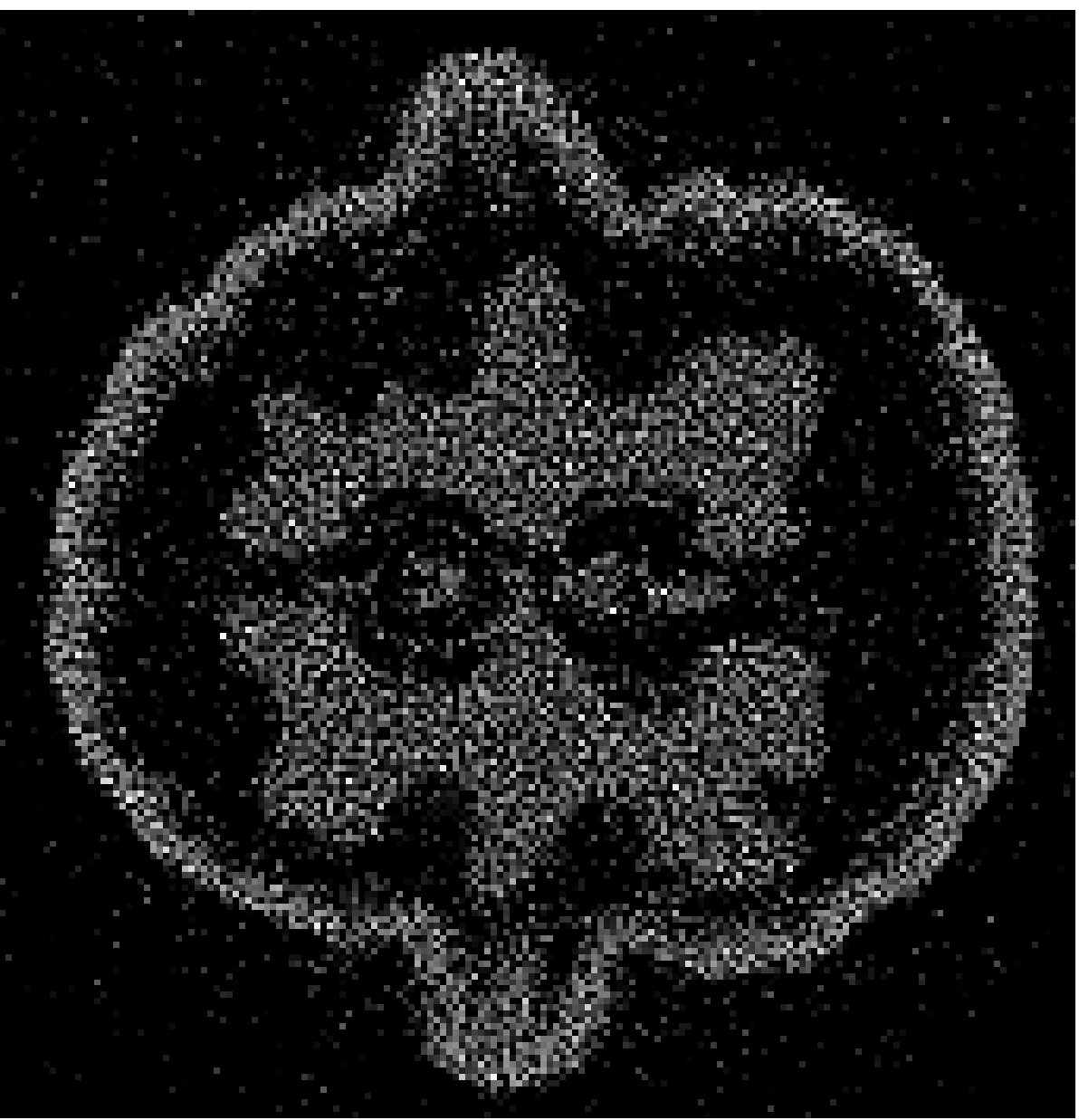}}
  \put(-70,140){\includegraphics[width=1.5in]{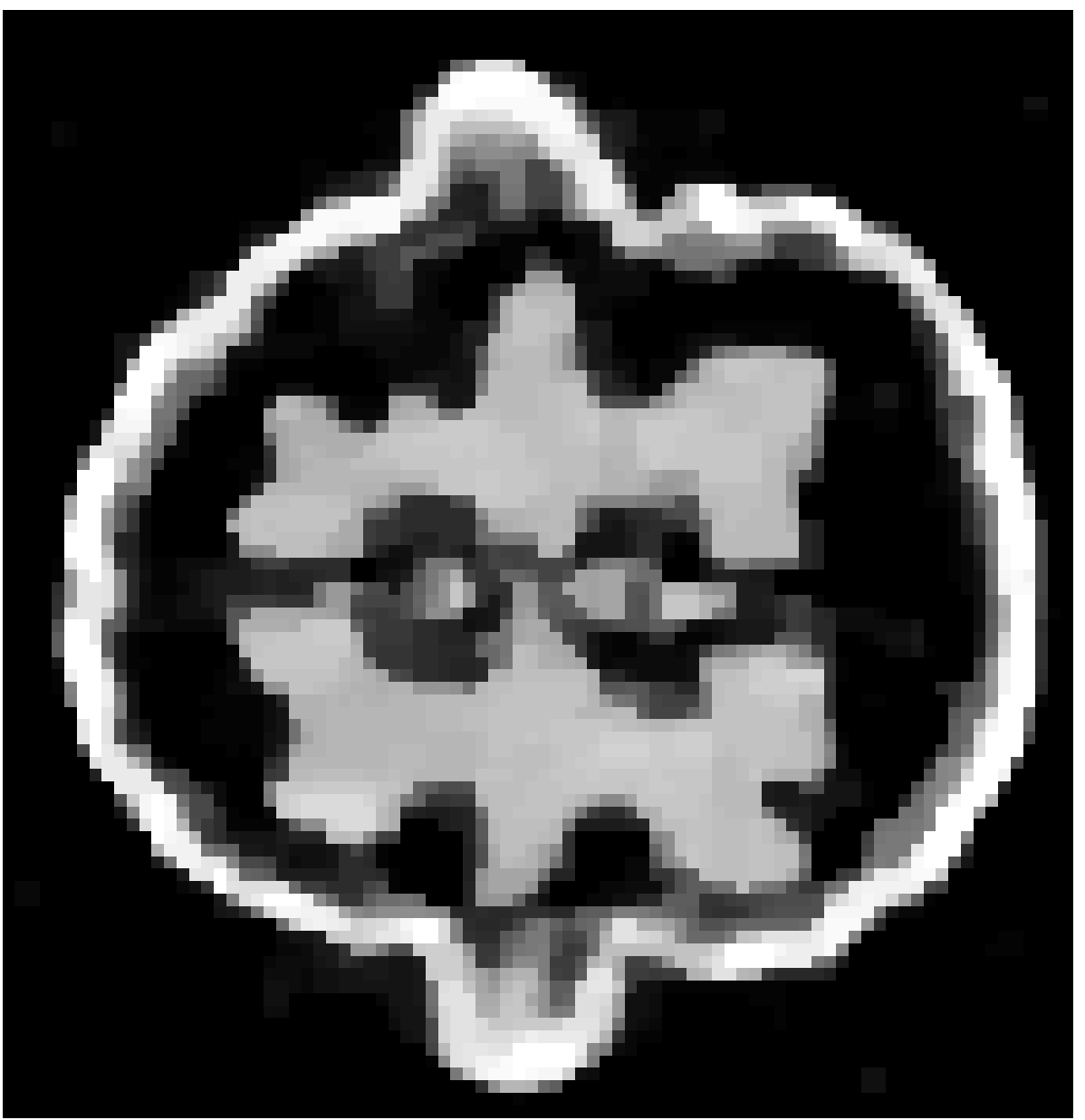}}
  \put(50,135.5){\includegraphics[width=1.5in]{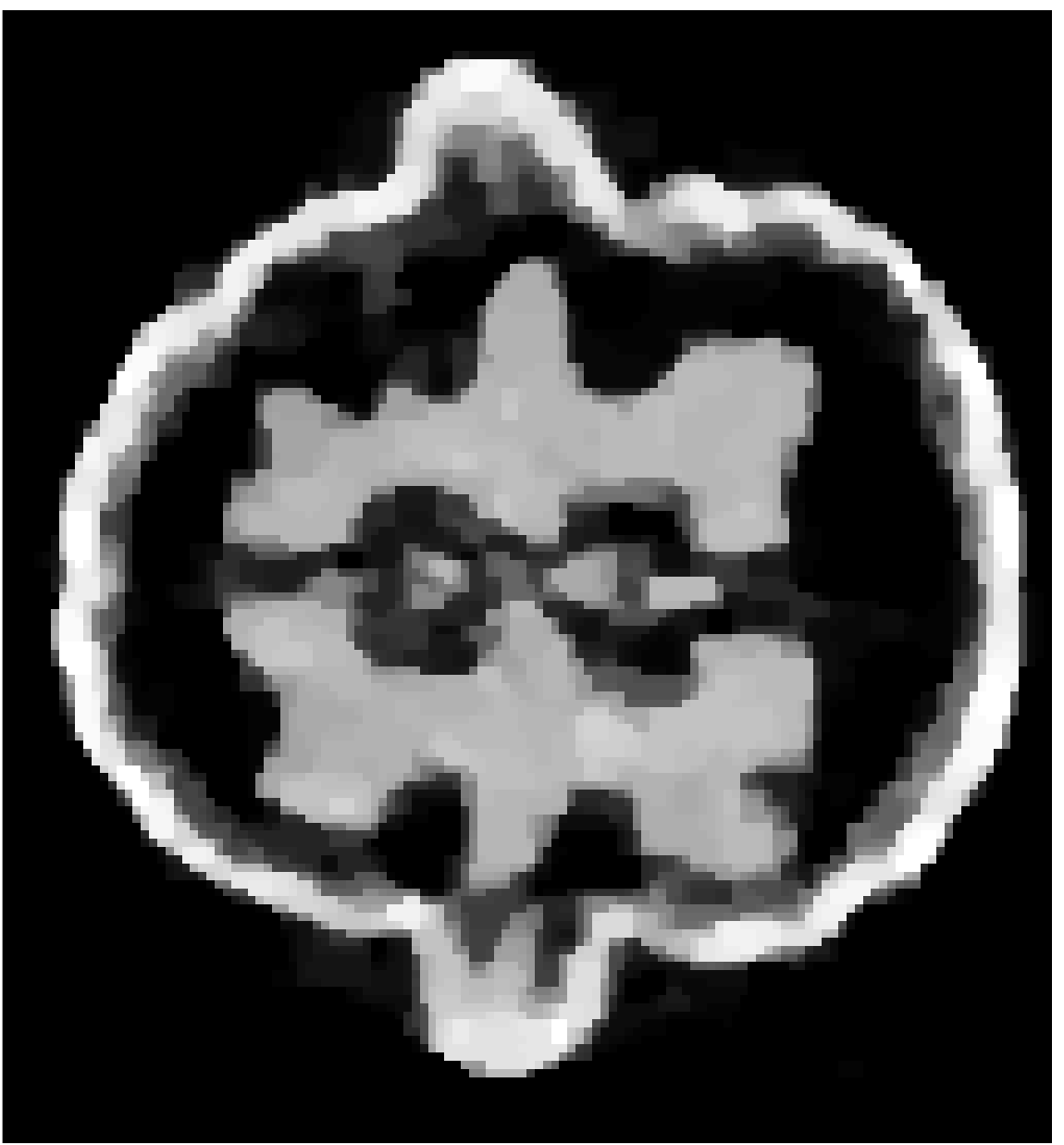}}
  \put(170,140){\includegraphics[width=1.5in]{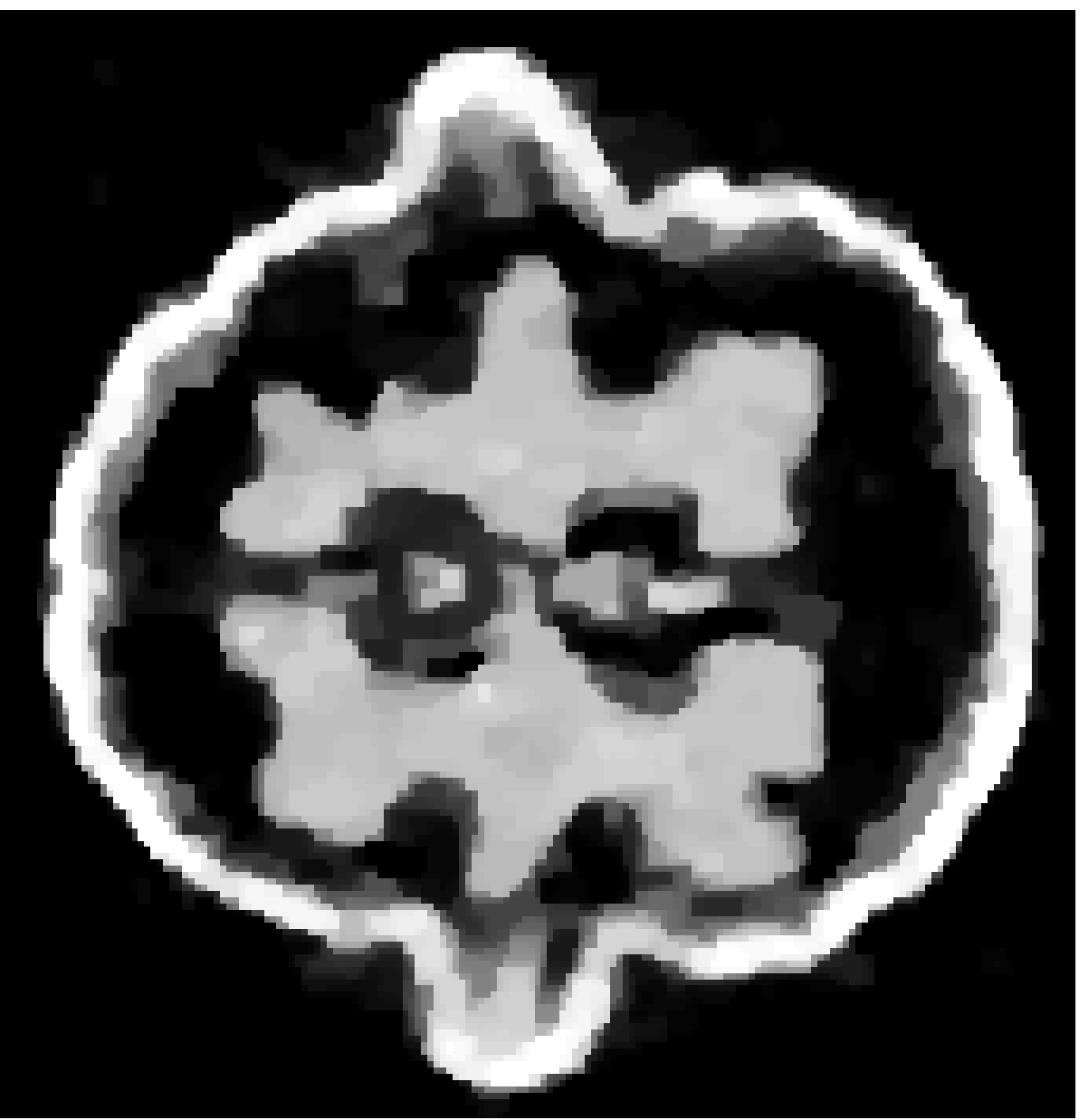}}
  \put(-70,20){\includegraphics[width=1.5in]{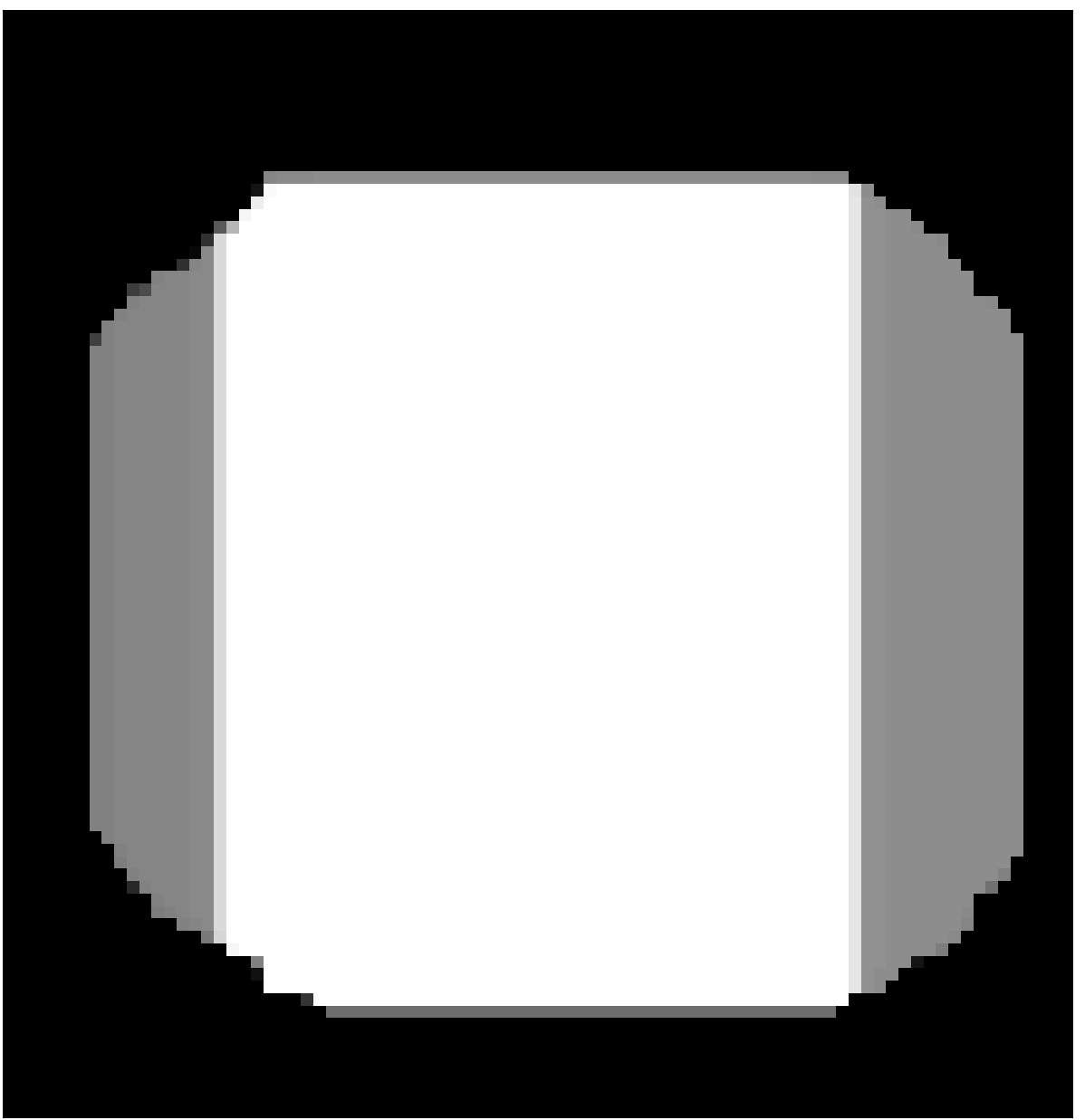}}
  \put(50,15.5){\includegraphics[width=1.5in]{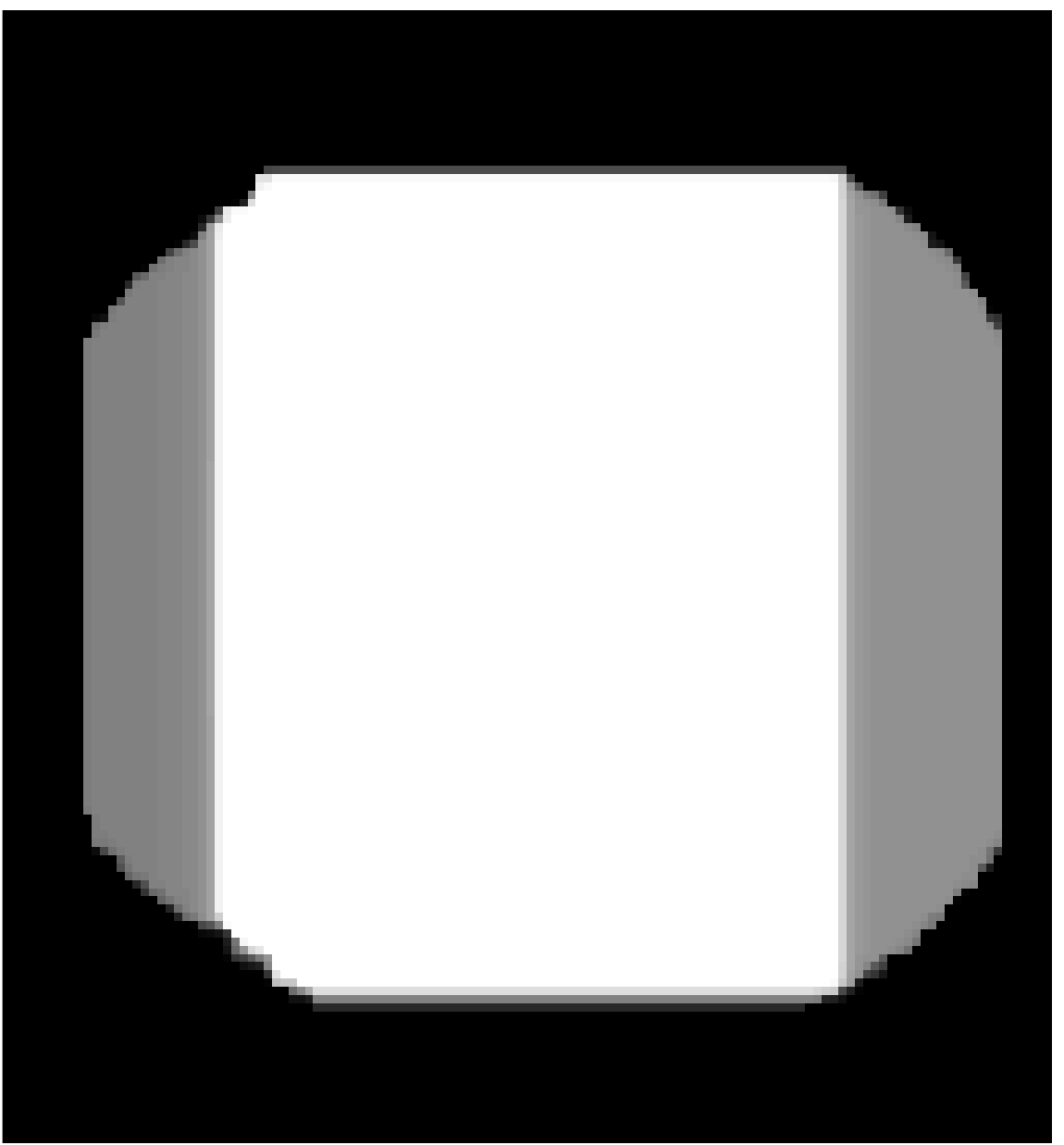}}
  \put(170,20){\includegraphics[width=1.5in]{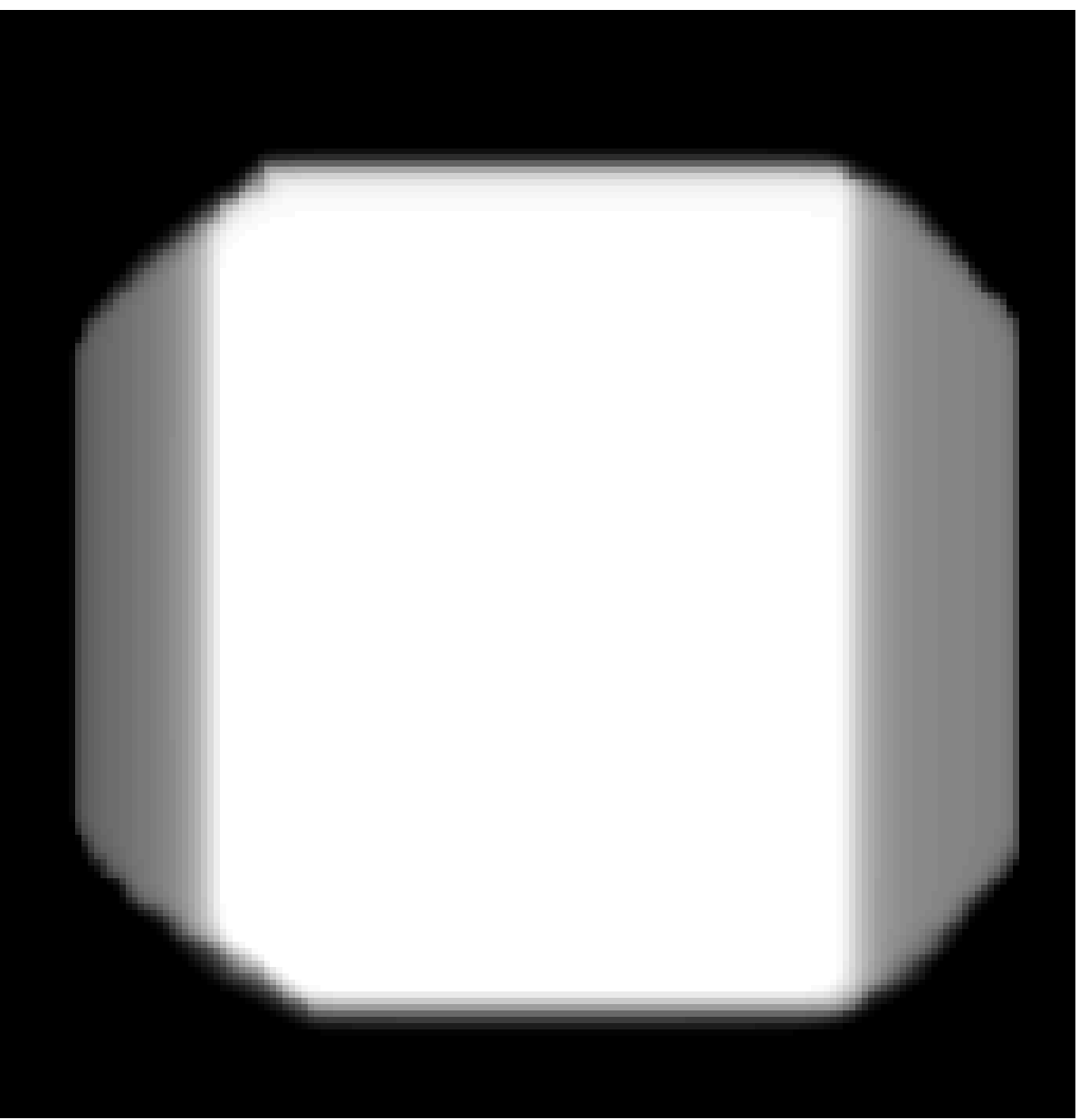}}
 \put(280,310){$\alpha = 10^{-3}$}
 \put(280,200){$\alpha = 10^{1}$}
 \put(280,80){$\alpha = 10^{4}$}
 \end{picture}
 \caption{\label{fig:TVrekotNoise5}TV regularized reconstructions $\fan$ from the measured projection data corrupted with 5\% additional noise computed with different values
   of regularization parameter $\alpha$. The values of $\alpha$ are
   denoted on the right. Discretization levels $n =
   128,192, 256$ are indicated on the top.}
 \end{figure}

Total variation norms $\| \fan \|_{{\rm TV}}$ of the reconstructed images were computed as a function of $\alpha$ at each discretization level $n$. Results from the measured data (case i)
and from the measured data corrupted with 5\% of random additional noise (case ii) are shown in
Table \ref{tab:doublenorms}. A few samples of the
corresponding reconstructions are shown in figures
\ref{fig:TVrekot} and \ref{fig:TVrekotNoise5} for test cases i) and
ii), respectively. 

The reconstructions in figures
\ref{fig:TVrekot} and \ref{fig:TVrekotNoise5} are selected as
follows. On the top row the value of the regularization parameter is too small
producing a "noisy" reconstruction and on the bottom row the
regularization parameter is too big producing an "over-regularized"
reconstruction. In the middle row the value of the regularization parameter
is selected as the smallest keeping the TV norms of reconstructions approximately independent of resolution. See Table
\ref{tab:doublenorms}.

\begin{table}[!h]
\caption{\label{tab:doublenorms} TV norms of 2D
  reconstructions $\fan$ computed onthree different
  discretization levels with several values of
  $\alpha$ ranging in the interval $[10^{-6}, 10^{6}]$. The discretization levels  were
  set to $n = 128, 192,  256$. The results are computed
  from projection data with 90 projection angles. Data is corrupted with
additional additive random noise (5$\%$).  Approximately resolution-independent TV norms are highlighted. The choice of $\alpha$ is the smallest possible leading to stable TV norms: $\alpha=1$ in the case of low noise, and $\alpha=10$ in the case of 5\% added noise.}
\begin{picture}(250,180) 
\setlength\fboxsep{0pt}
\put(16.5,-8){\fcolorbox{VeryLightOrange}{VeryLightOrange}{\framebox(90,96){}}}
\put(113.5,-8){\fcolorbox{VeryLightOrange}{VeryLightOrange}{\framebox(90,82){}}}
\put(13,165){{\bf Low noise}}
\put(112.5,165){{\bf 5\% noise}}
\put(-20,70){\begin{minipage}[t]{0.53\linewidth}
\begin{tabular}{l|lll|lll}
{$\alpha$}& {$128^2$} & {$192^2$} & {$256^2$} & {$128^2$} & {$192^2$} & {$256^2$} \\
\hline
{$10^{-4}$} &    1.51 & 2.29  & 3.64 & 2.42  &  5.05 &  8.71 \\
{$10^{-3}$} &    1.51 & 2.29  & 3.46 & 2.43  &  5.05 &  8.59 \\
{$10^{-2}$} &    1.50 & 2.23  & 2.97 & 2.42  &  5.01 &  8.59 \\
{$10^{-1}$} &    1.43 & 1.85  & 1.93 & 2.37  &  4.83 &  8.16 \\
{$10^{0}$} &     1.08 & 1.11  & 1.11 &  1.99  &  3.50 &  5.12 \\
{$10^{1}$} &     0.78 & 0.78  & 0.77 &  0.86  &  0.86 &  0.88 \\
{$10^{2}$} &     0.48 & 0.48  & 0.48 &  0.48  &  0.48 &  0.48 \\
{$10^{3}$} &     0.12 & 0.12  & 0.12 &  0.12  &  0.12 &  0.12 \\
{$10^{4}$} &     0.04 & 0.04   & 0.04 &  0.04  &  0.04 &  0.04 \\
{$10^{5}$} &     0      & 0       & 0     &     0      & 0       & 0      \\
{$10^{6}$} &     0      & 0       & 0     &     0      & 0       & 0      \\
\hline
\end{tabular}
\end{minipage}}
\end{picture}
\end{table}

\subsection{Selection of $\alpha$ using the S-curve method}

The S-curve method was originally introduced in
\cite{Kolehmainen2012,Hamalainen2013} within a Bayesian inversion framework. In
this work we apply the S-curve method for  TV regularized
tomography.
 
The S-curve was applied as follows. For test case i), we computed
TV regularized reconstructions $\fan$ with fixed $n$ for 14 different values of $\alpha$ 
ranging in the interval $[10^{-6}, 10^{6}]$  and computed the 
TV norms  
$$ S(\alpha) := \| \fan \|_{{\rm TV}}$$
for the reconstructions.
Then the data $\{\alpha, S(\alpha)\}$ was interpolated 
to get the S-curve. The S-curve was computed using the same 
three different discretization levels ($n=128,192,256$) that were
used in section \ref{secwh2}.

In the test case ii) we computed
reconstructions with 14 different values of $\alpha$ ranging in the
interval $[10^{-4}, 10^{8}]$, and applied the S-curve interpolation procedure similarly 
as in case i). 

\begin{figure}
\begin{picture}(340,120)
\put(0,0){\includegraphics[width=4cm]{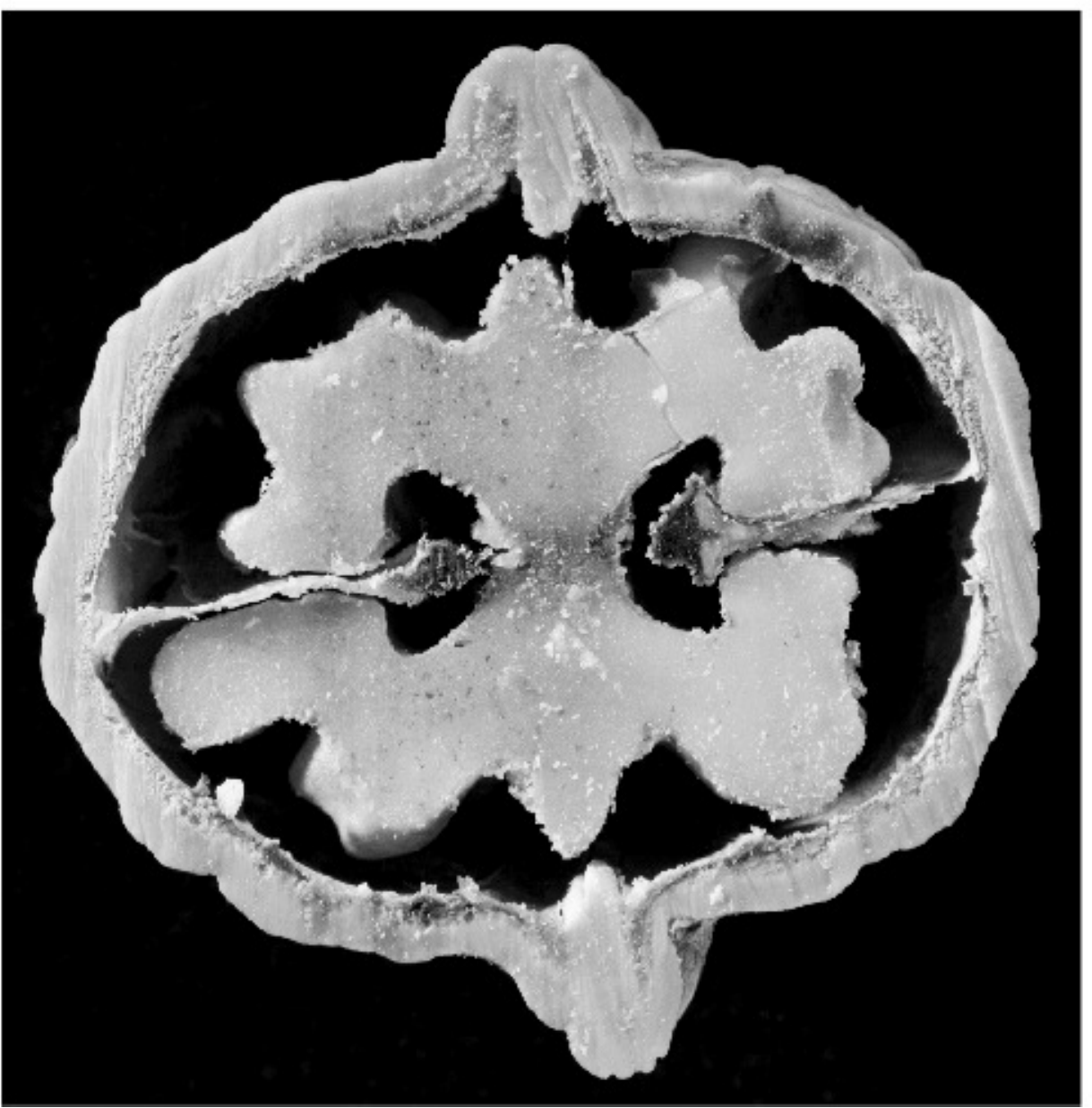}}
\put(120,0){\includegraphics[width=4cm]{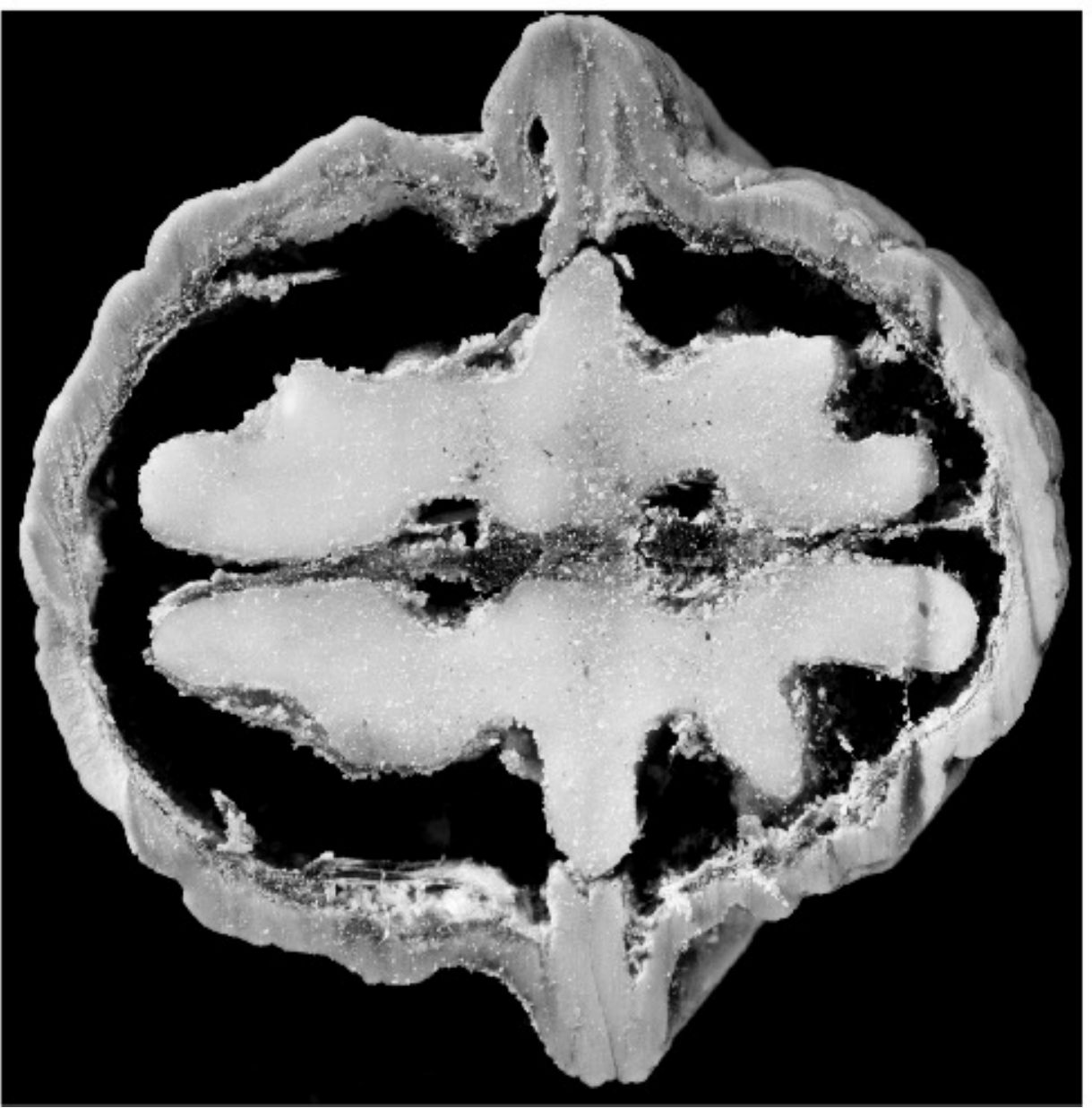}}
\put(240,0){\includegraphics[width=4cm]{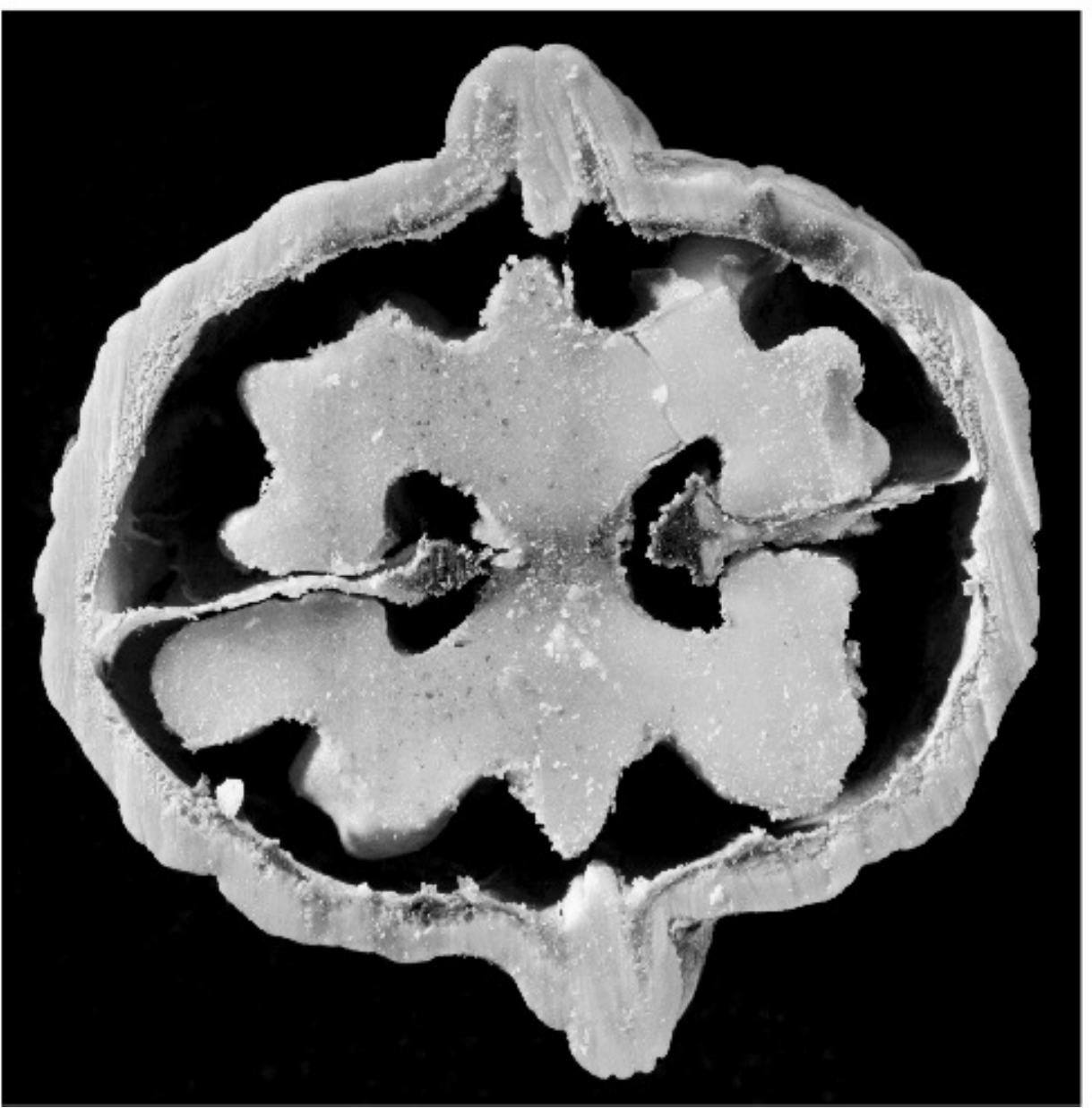}}
\end{picture}
 \caption{\label{fig:Walnut_photos}Photographs of three walnuts split in half. These images were used to provide {\em a priori} information about the expected total variation of the target. Of course, these are optical photographs and thus physically very different objects from the reconstructions (tomographic slices representing X-ray attenuation coefficient). 
However, we assume that the total variations of the photographs, when scaled by (\ref{photoscale}), are comparable to the total variation of the X-ray attenuation coefficient. }
 \end{figure}

As in \cite{Hamalainen2013} the {\it a priori} value for the sparsity level
$\hat{S}$ was estimated using three digital photographs of split
 walnuts. 
The photographs of the three walnuts are shown in figure
\ref{fig:Walnut_photos}. The walnut used to measure the X-ray data
is not included in the photos. 
 In reality, as is the case also in this work, the {\it a
   priori} information, that we use to estimate the value of
 $\hat{S}$, comes from a different modality (e.g. anatomical atlases) 
than the one we are
considering (here digital photographs versus X-ray attenuation
function). Therefore in order to compute the $\hat{S}$ for the total
variation regularized X-ray tomography, each of the digital photographs $\f_{\rm p}$ were
 scaled such that the norm of the computed X-ray projection data of the
 photograph is the same as the norm of the measured projection data
 $\greal$. This was obtained by
\begin{equation} \label{photoscale}
\widetilde \f_{\rm p} = \frac{\| \greal \|}{\| A \f_{\rm p} \|} \f_{\rm p}
\end{equation}
This scaling of the photographs is essential since the total variations of the 
photographs are not directly comparable to the
total variation of the X-ray attenuation function.

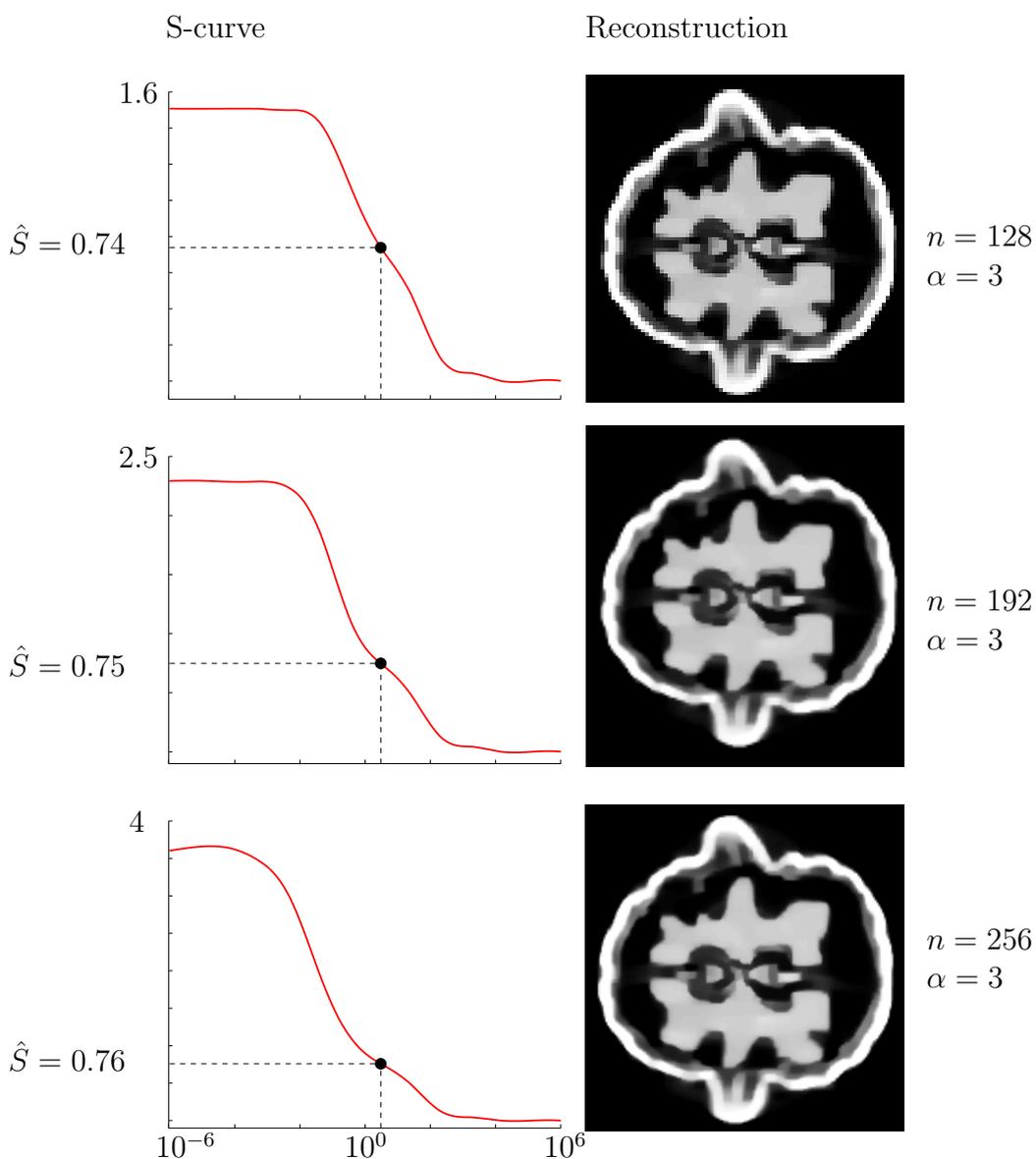
\begin{figure}
 \begin{picture}(200,450)
 \put(-40,430){S-curve}
 \put(120,430){Reconstruction}
  \put(-40,290){\includegraphics[width=2.1in]{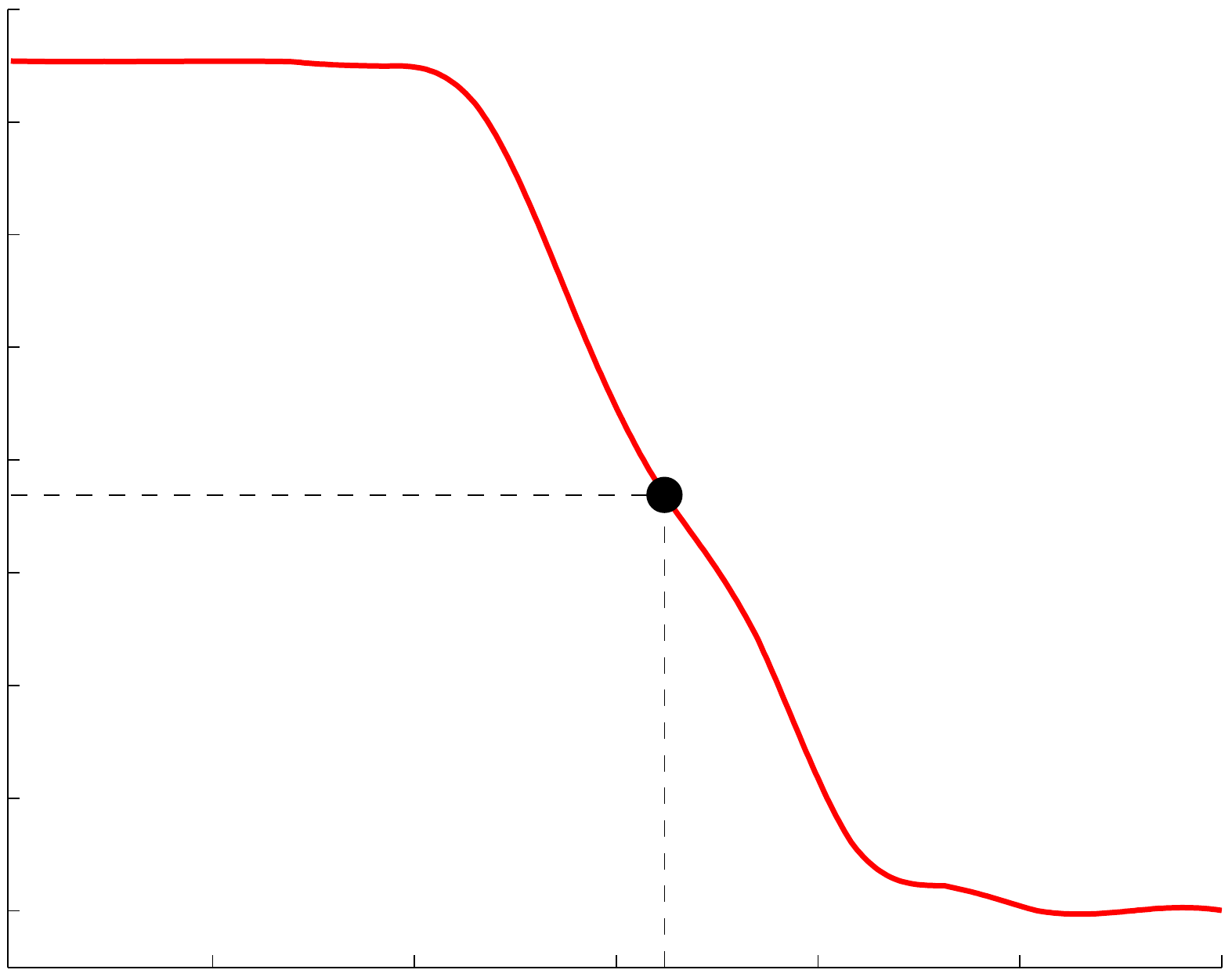}}
  \put(120,290){\includegraphics[width=1.7in]{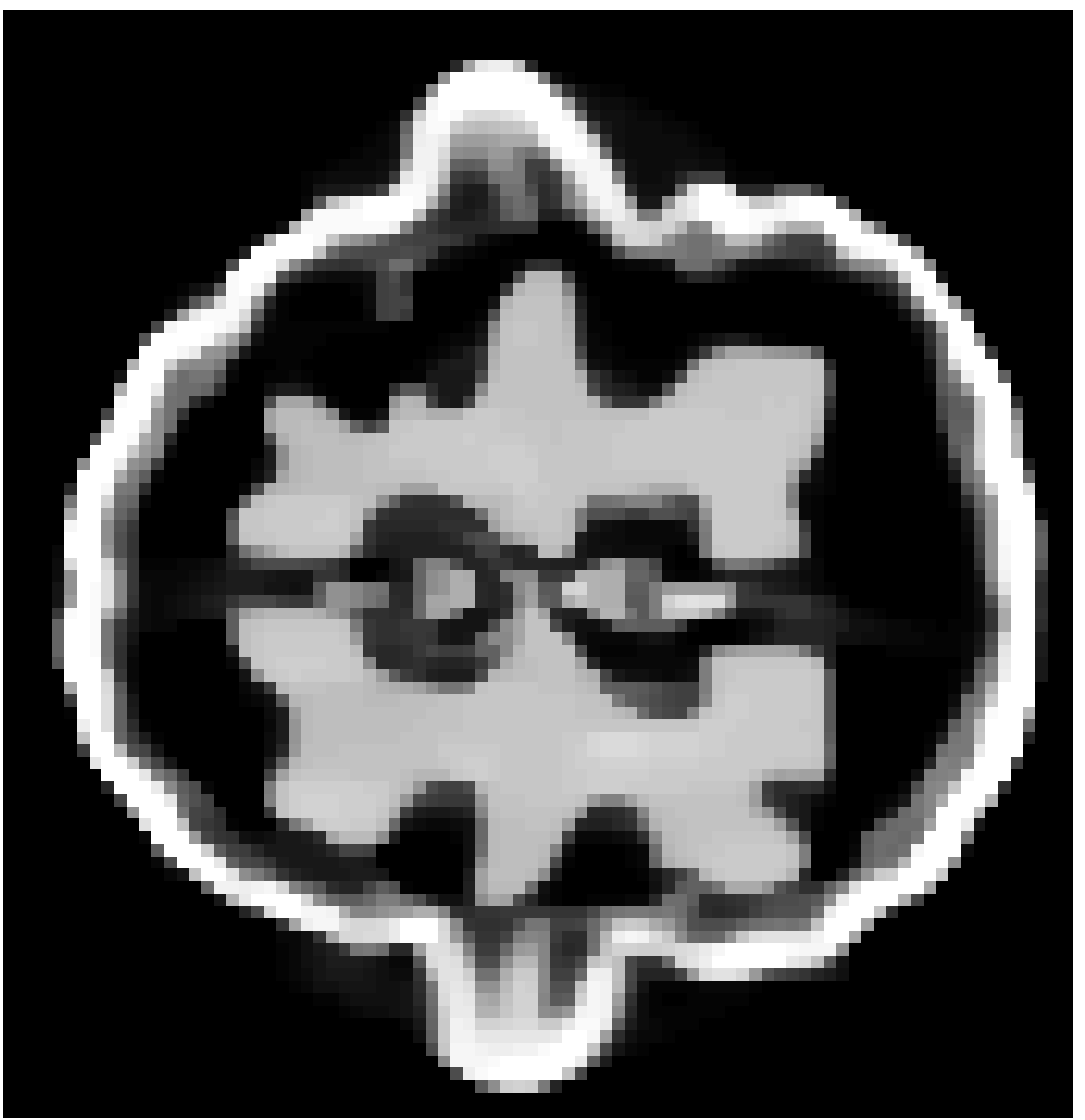}}
  \put(-40,150){\includegraphics[width=2.1in]{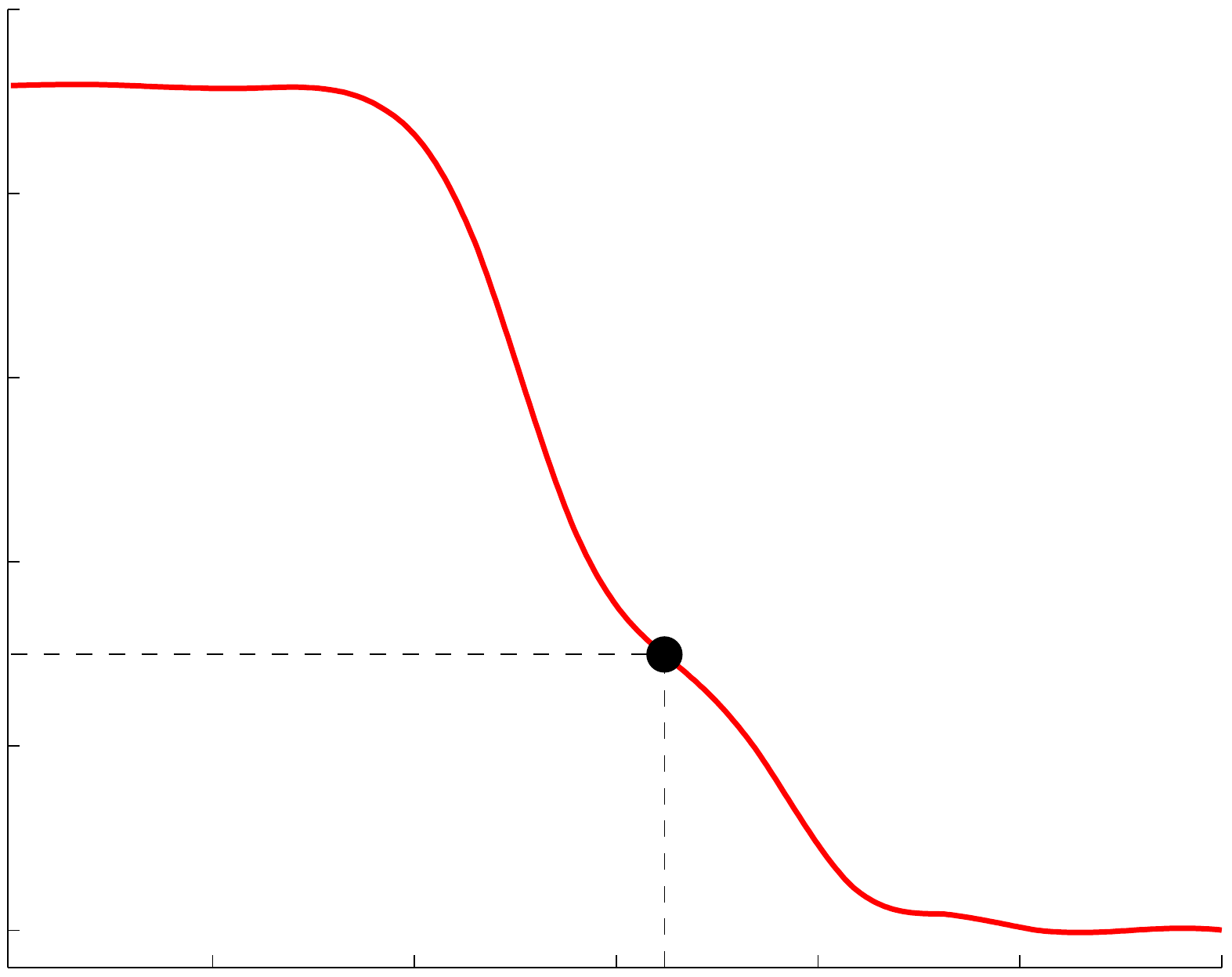}}
  \put(120,150){\includegraphics[width=1.7in]{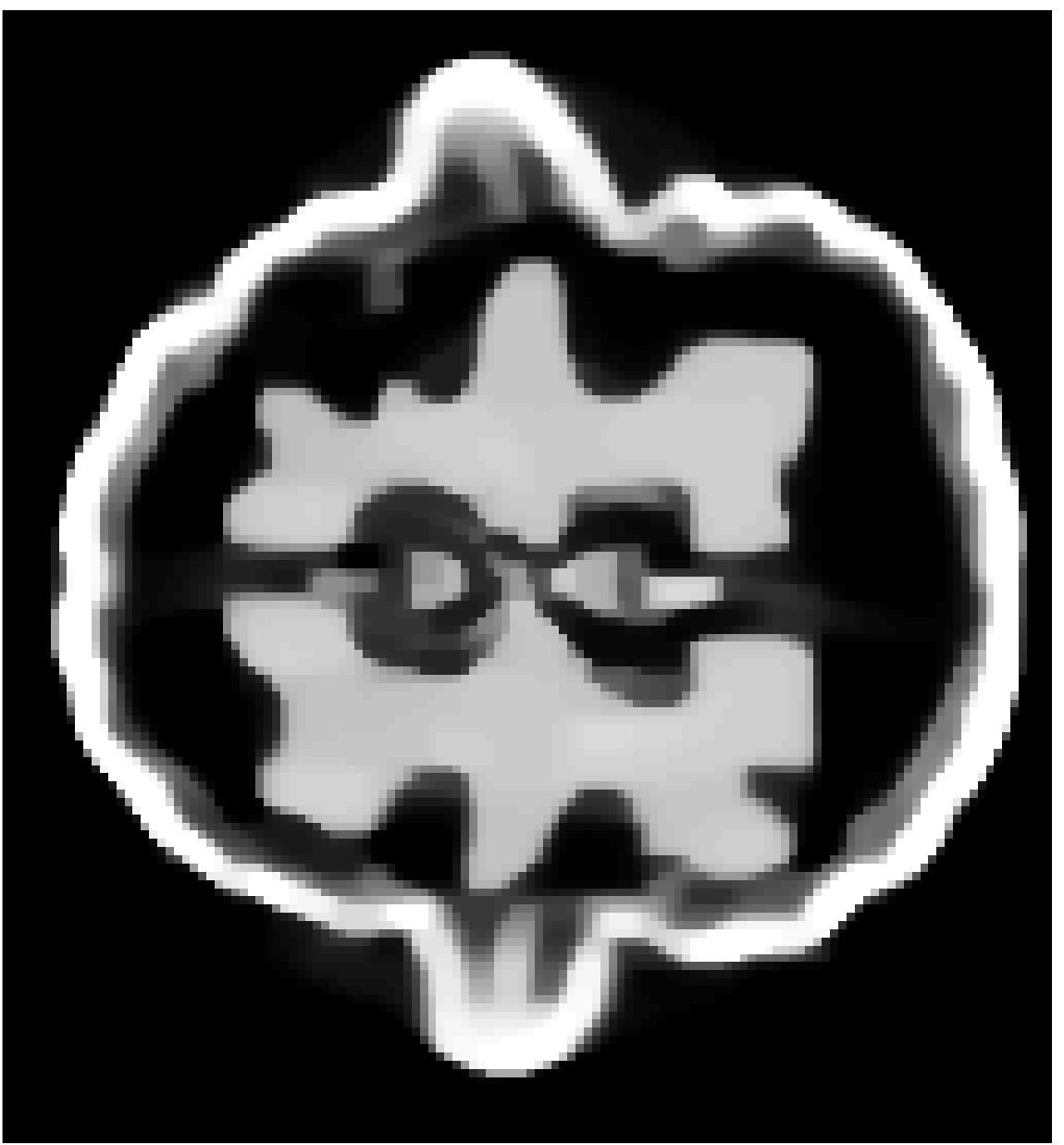}}
  \put(-40,10){\includegraphics[width=2.1in]{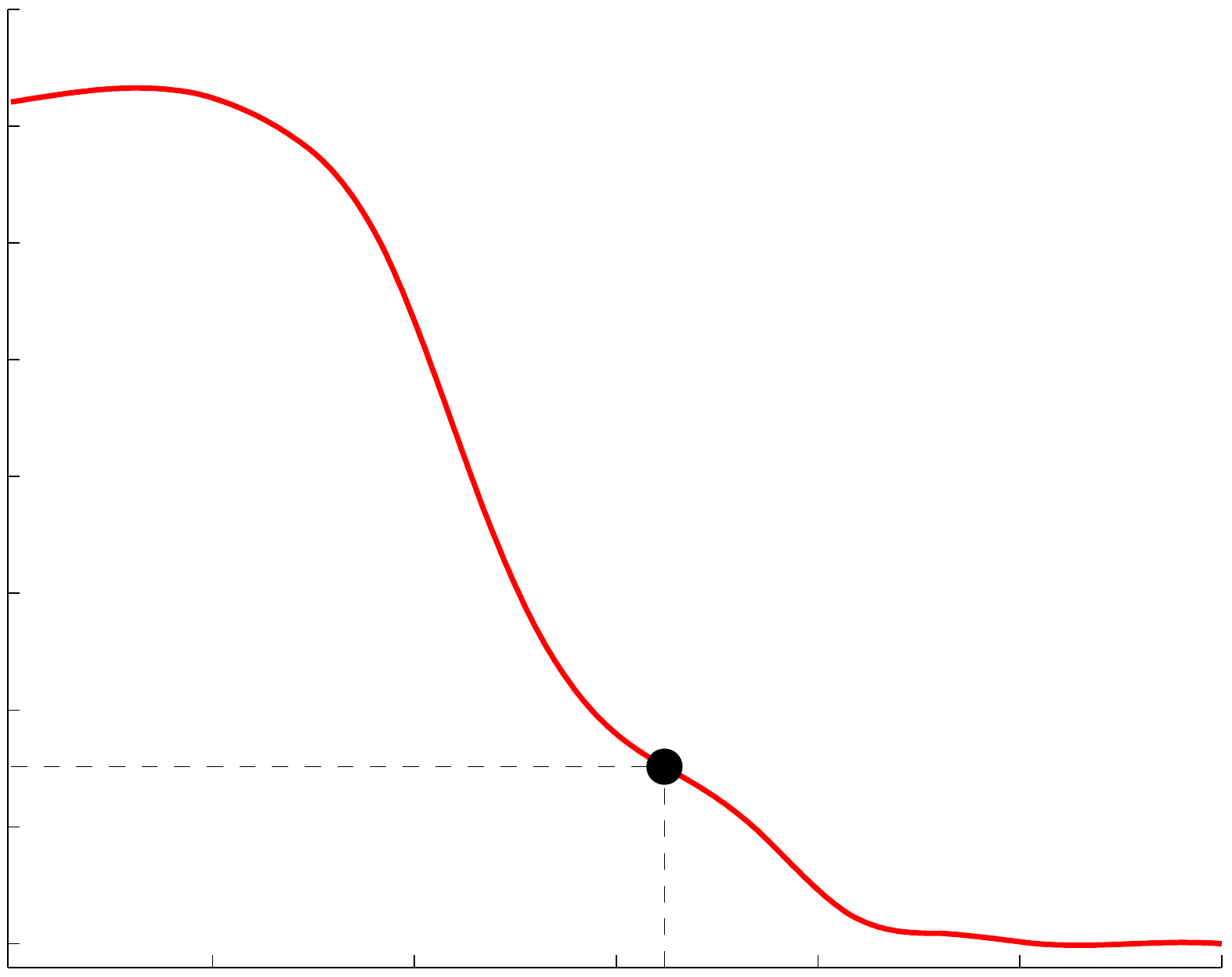}}
  \put(120,10){\includegraphics[width=1.7in]{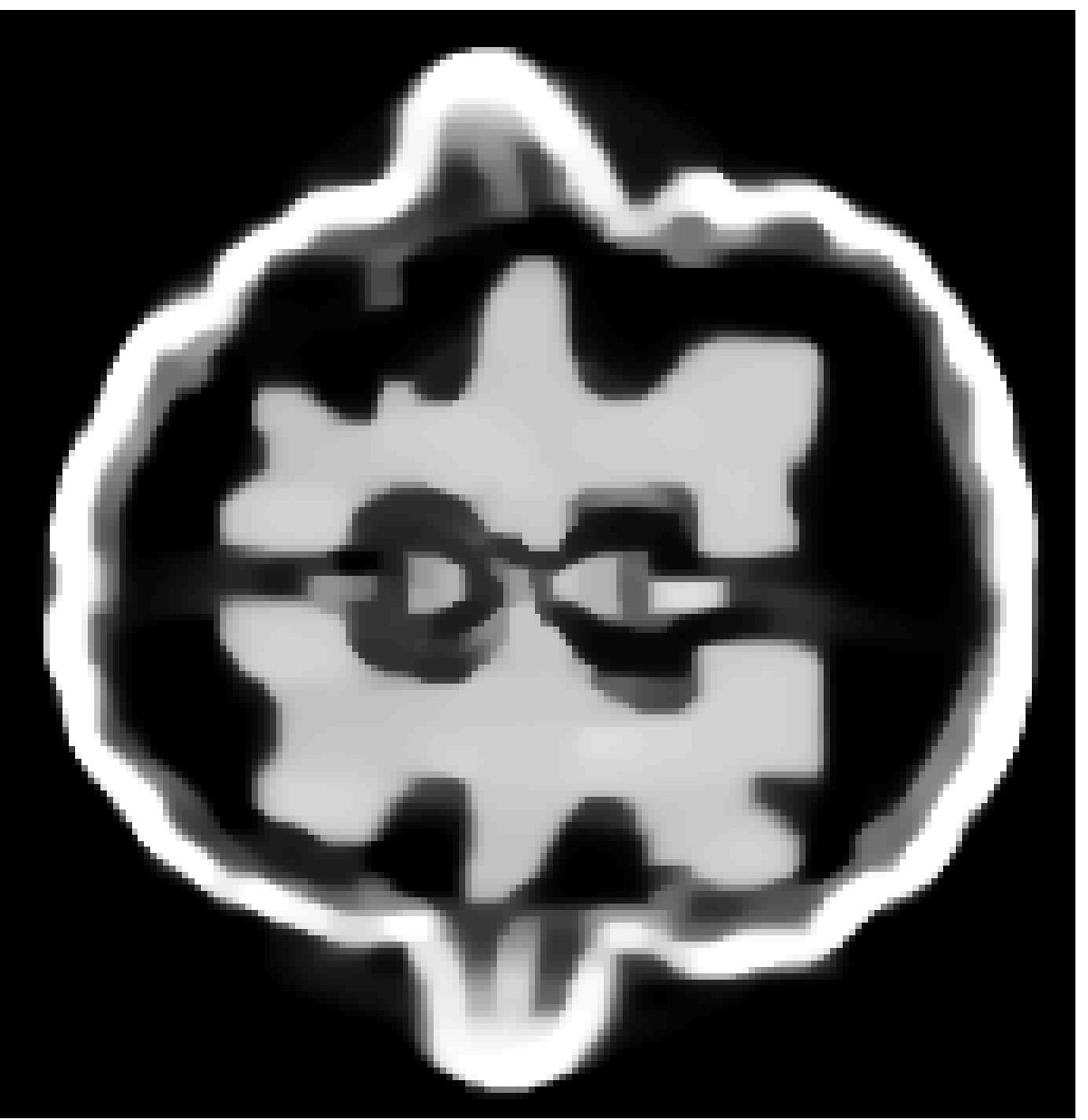}}
 \put(250,350){$n=128$}
\put(250,335){$\alpha = 3$}
 \put(250,210){$n=192$}
\put(250,195){$\alpha = 3$}
 \put(250,80){$n=256$}
\put(250,65){$\alpha = 3$}
\put(-100,347){$\hat{S} = 0.74$}
\put(-100,185){$\hat{S}= 0.75$}
\put(-100,33){$\hat{S}=0.76$}
\put(-44,0){$10^{-6}$}
\put(104,0){$10^{6}$}
\put(-58,405){$1.6$}
\put(-58,265){$2.5$}
\put(-54,125){$4$}
\put(28,0){$10^{0}$}
 \end{picture}
 \caption{\label{fig:Scurve} S-curve  method applied to the measured projection data with no    additional noise (case i). Left column: Plots of the S-curves used to    determine the values of $\alpha$. Right column: Reconstruction    computed with the selected $\alpha$.}
 \end{figure}

\begin{figure}
 \begin{picture}(200,450)
 \put(-40,430){S-curve}
 \put(120,430){Reconstruction}
  \put(-40,290){\includegraphics[width=2.1in]{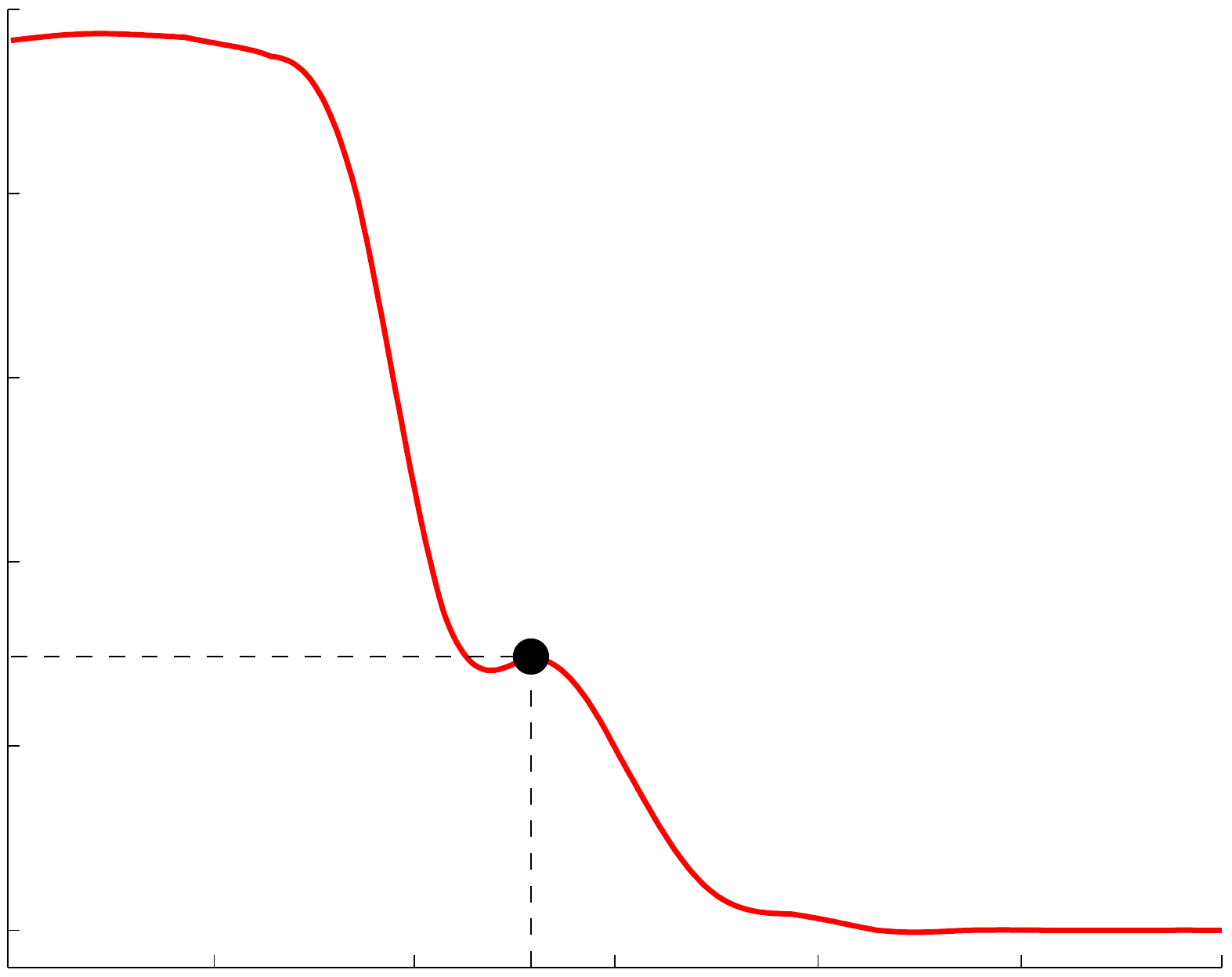}}
  \put(120,290){\includegraphics[width=1.7in]{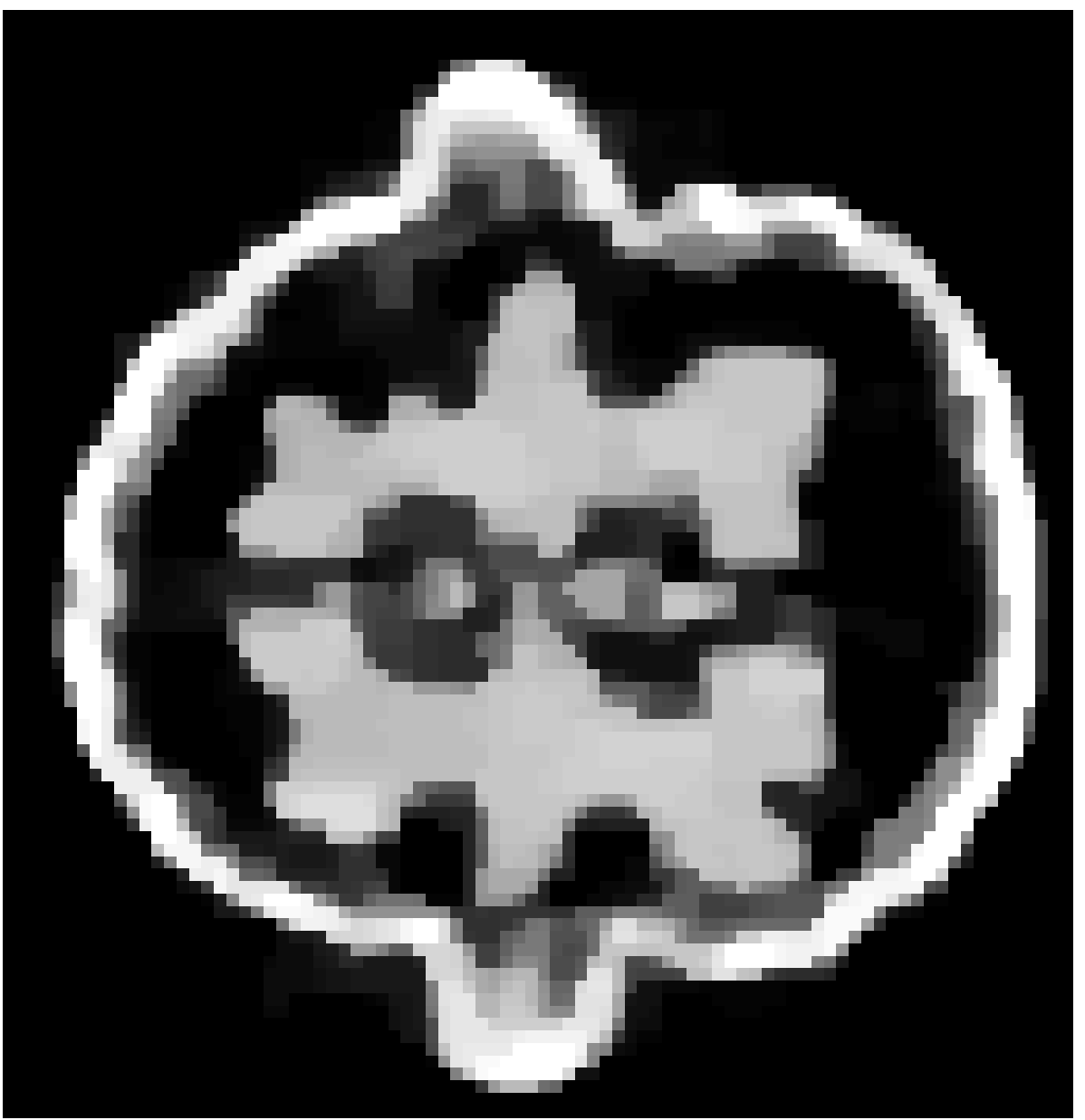}}
  \put(-40,150){\includegraphics[width=2.1in]{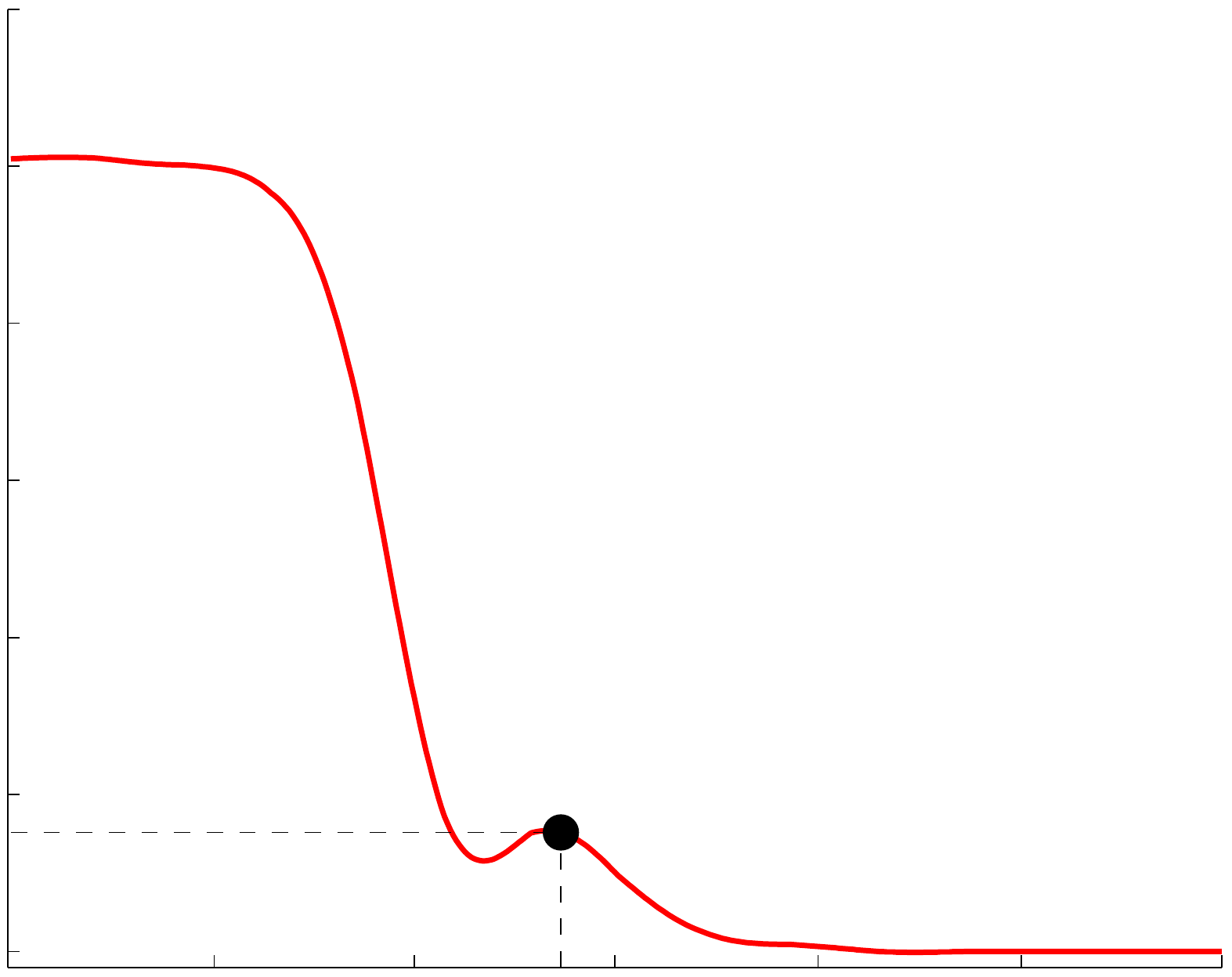}}
  \put(120,150){\includegraphics[width=1.7in]{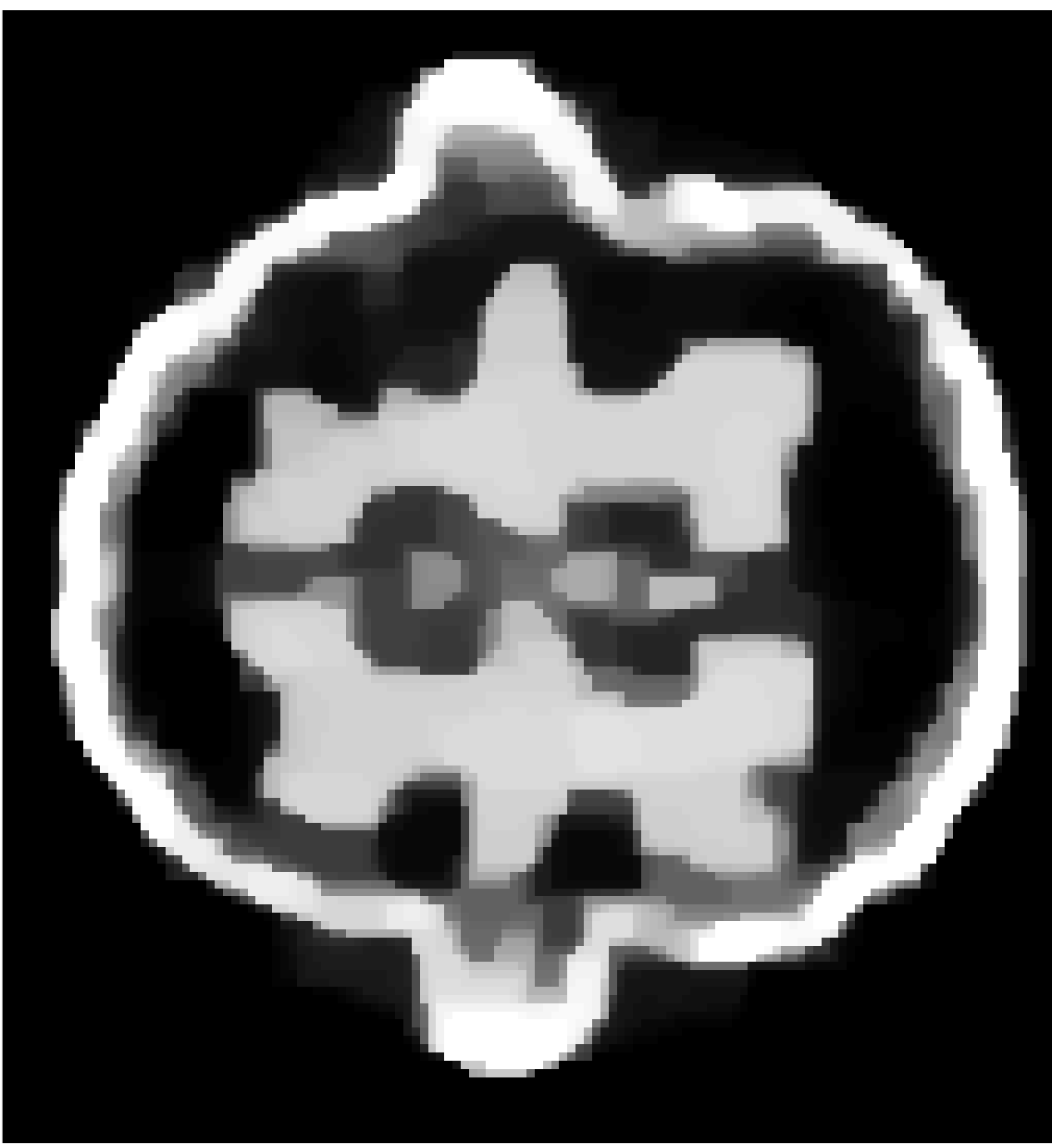}}
  \put(-40,10){\includegraphics[width=2.1in]{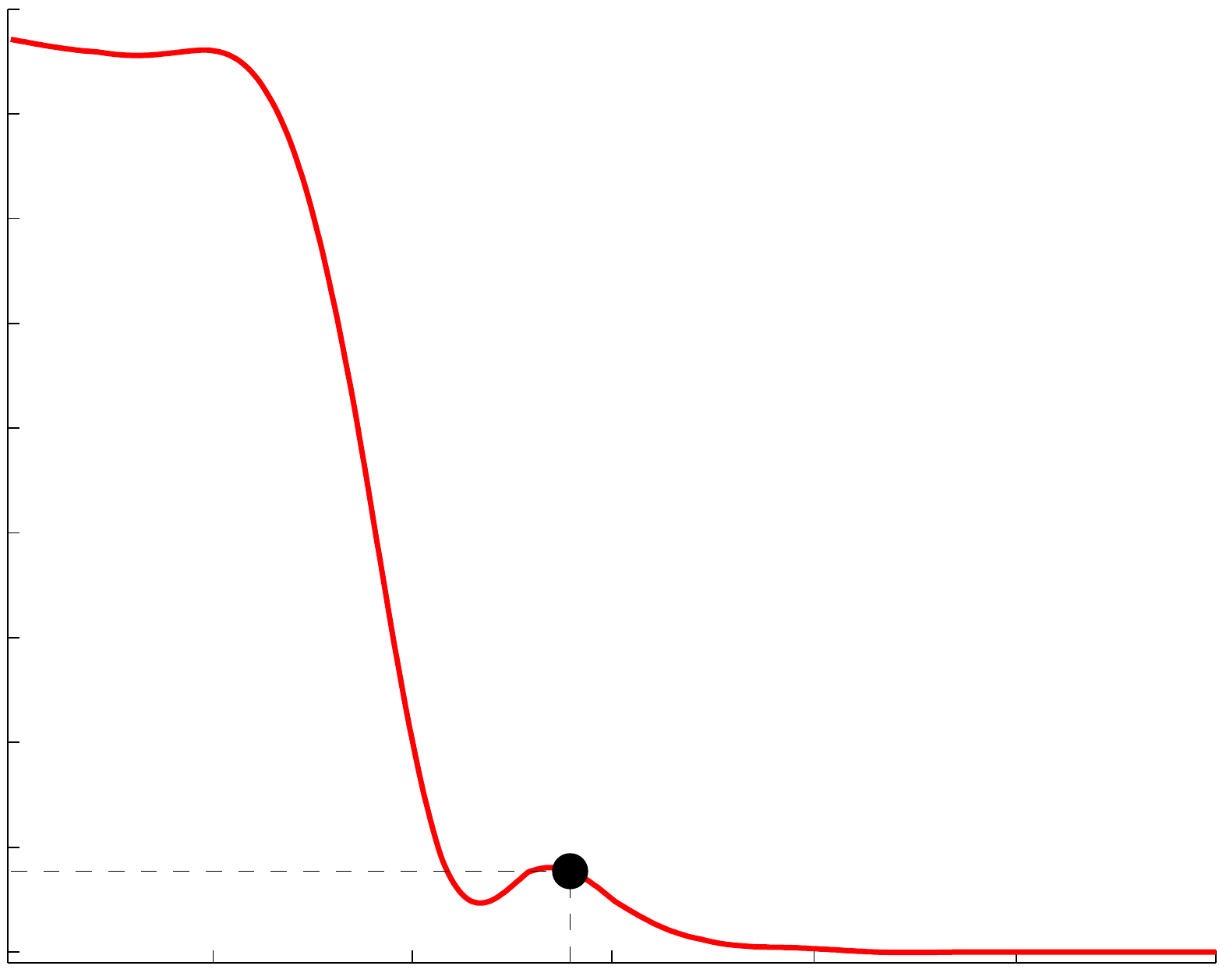}}
  \put(120,10){\includegraphics[width=1.7in]{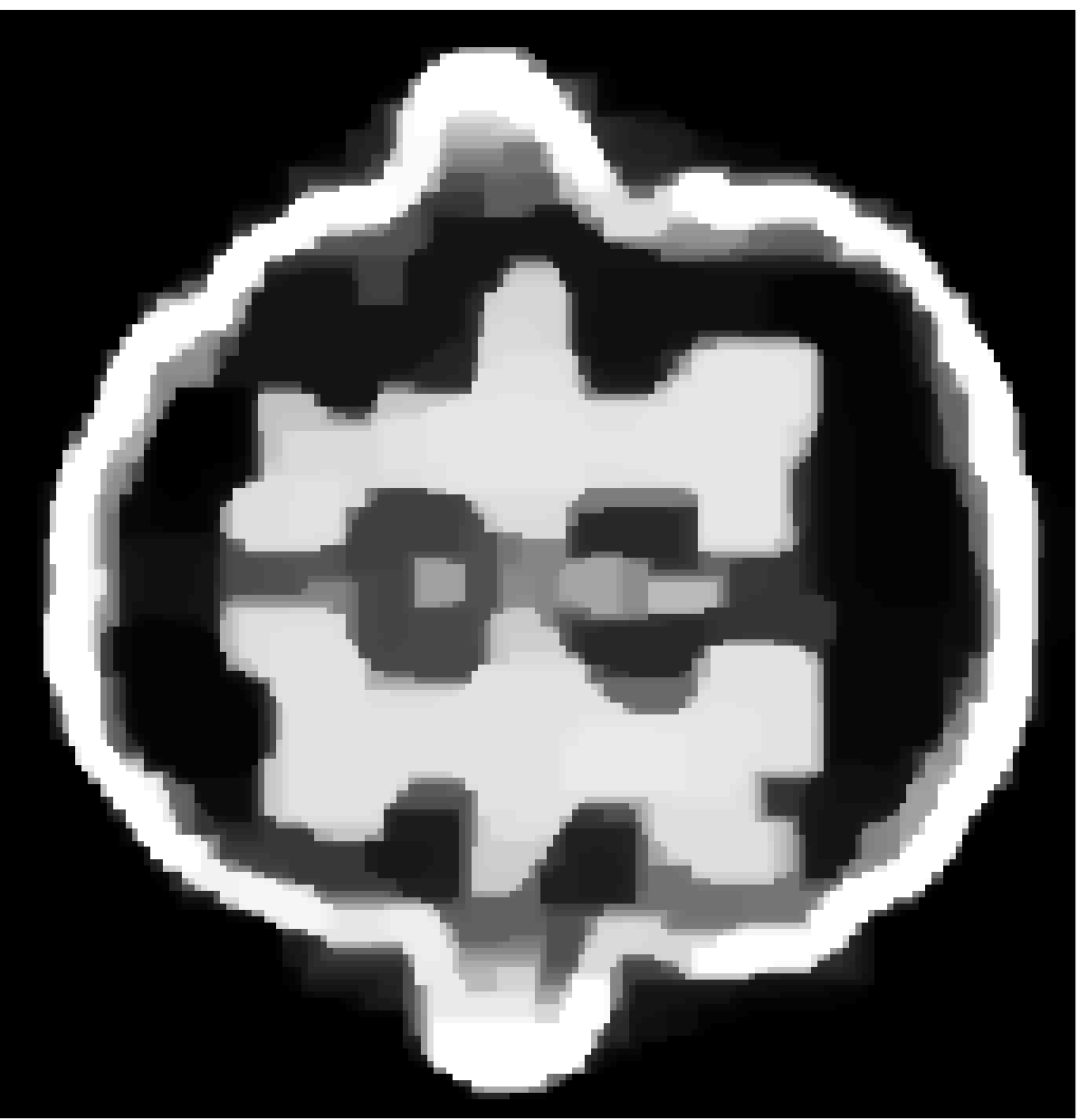}}
 \put(250,350){$n=128$}
\put(250,335){$\alpha = 14$}
 \put(250,210){$n=192$}
\put(250,195){$\alpha = 28$}
 \put(250,80){$n=256$}
\put(250,65){$\alpha = 37$}
\put(-100,327){$\hat{S}=0.74$}
\put(-100,165){$\hat{S}=0.75$}
\put(-100,22){$\hat{S}=0.76$}
\put(-54,0){$10^{-4}$}
\put(104,0){$10^{8}$}
\put(-58,405){$2.5$}
\put(-52,265){$6$}
\put(-52,125){$9$}
\put(28,0){$10^{2}$}
 \end{picture}
 \caption{\label{fig:ScurveNoise5} S-curve  method applied to the
   measured projection data    corrupted with 5\% additional noise
   (case ii). Left column: Plots of the S-curves which where used to
   determine the values of $\alpha$. Right column: Reconstruction
   computed using the selected $\alpha$.  }
 \end{figure}

The S-curve plots and the resulting reconstructions for each
discretization levels are presented in figure
\ref{fig:Scurve}. Corresponding results with 5 \% additional noise are
presented in figure \ref{fig:ScurveNoise5}.

\begin{figure}
 \begin{picture}(200,160)
 \put(-40,150){L-curve}
 \put(120,150){Reconstruction}
\put(120,140){$\alpha = 1358$} 
 \put(-50,10){\includegraphics[width=2.1in]{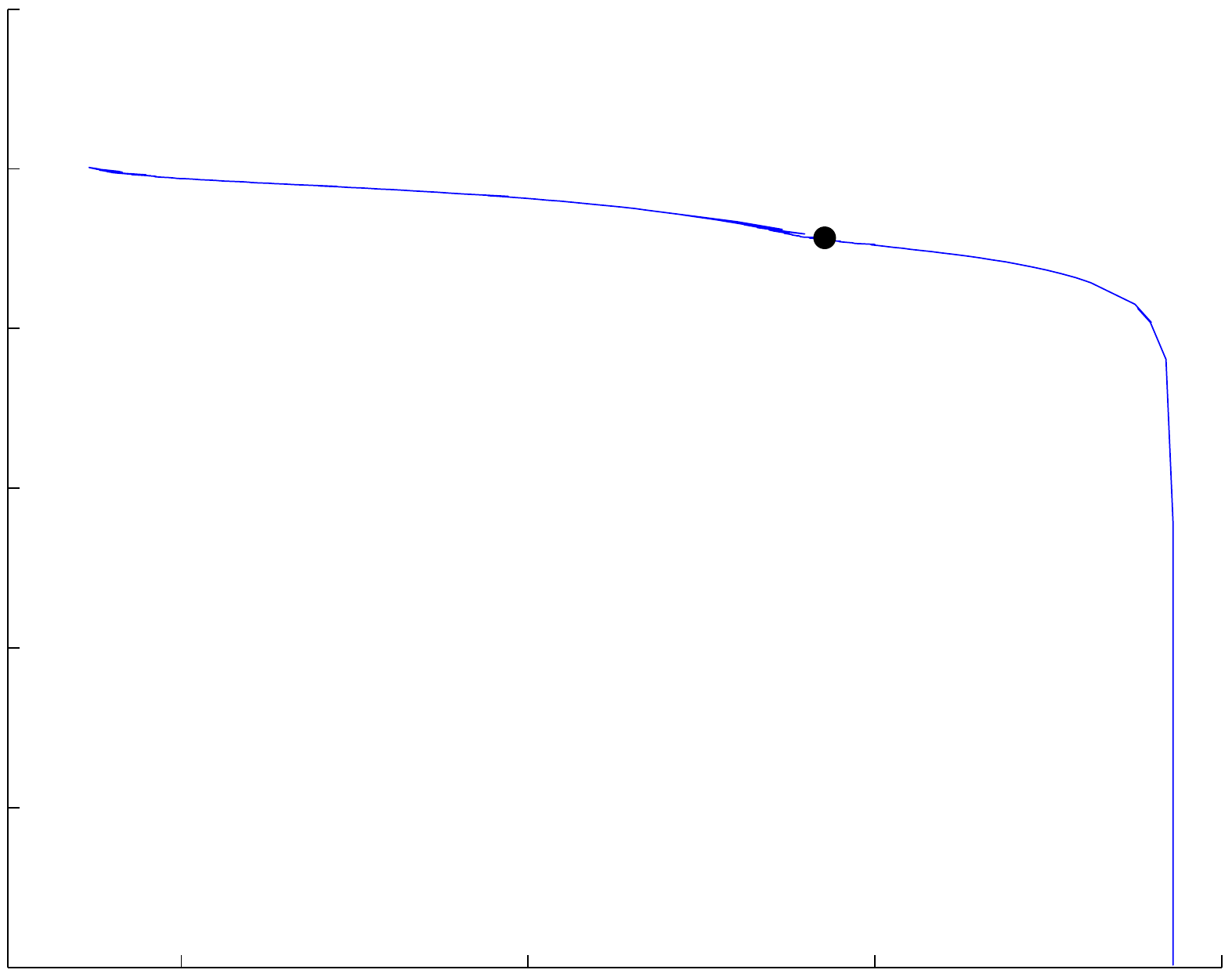}}
  \put(120,10){\includegraphics[width=1.7in]{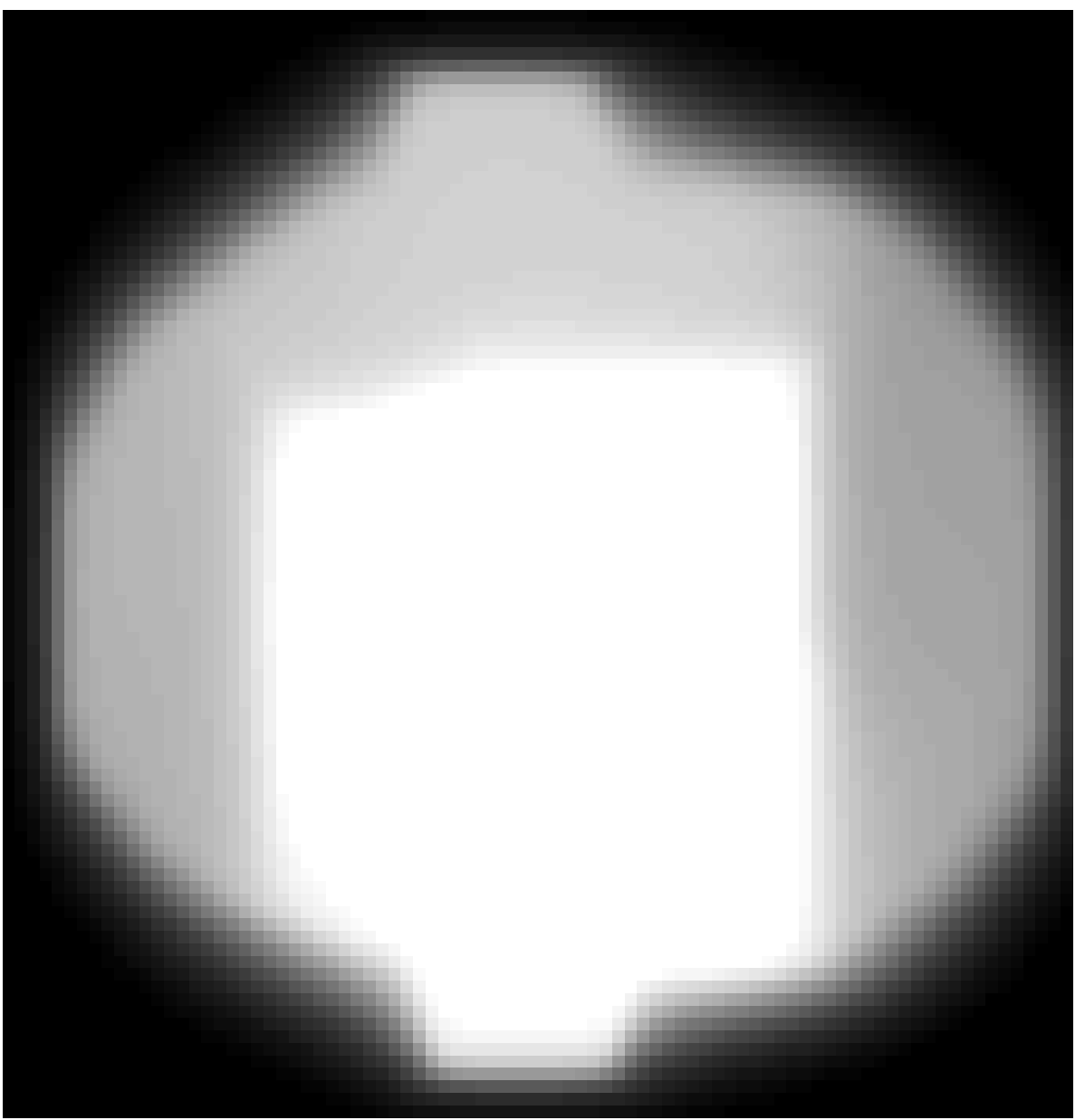}}
\put(-85,65){$TV(\f)$}
\put(-54,0){$10^{-5}$}
\put(94,0){$10^{6}$}
\put(0,0){$\Vert A \f - \greal \Vert$}
 \end{picture}
 \caption{\label{fig:Lcurve}L-curve  method applied to walnut
   data, no additional noise. Left: L-curve plot. Right: Reconstruction
   computed using the selected $\alpha = 1358$. The discretization level is $n=128$.}
 \end{figure}

\begin{figure}
  \begin{picture}(200,160)
  \put(-40,150){L-curve}
  \put(120,150){Reconstruction with}
 \put(120,140){$\alpha = 10.4$} 
  \put(-50,10){\includegraphics[width=2.1in]{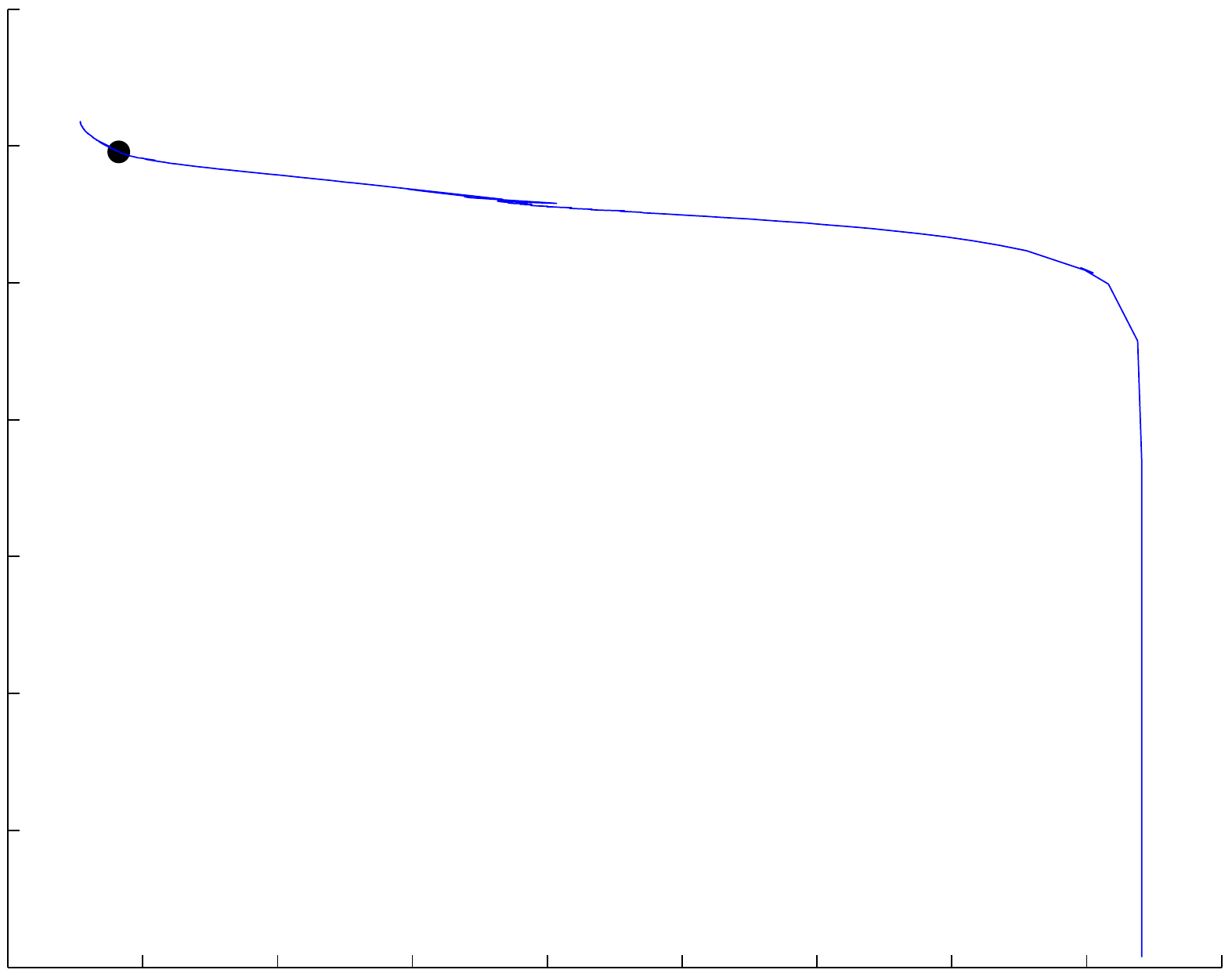}}
   \put(120,10){\includegraphics[width=1.7in]{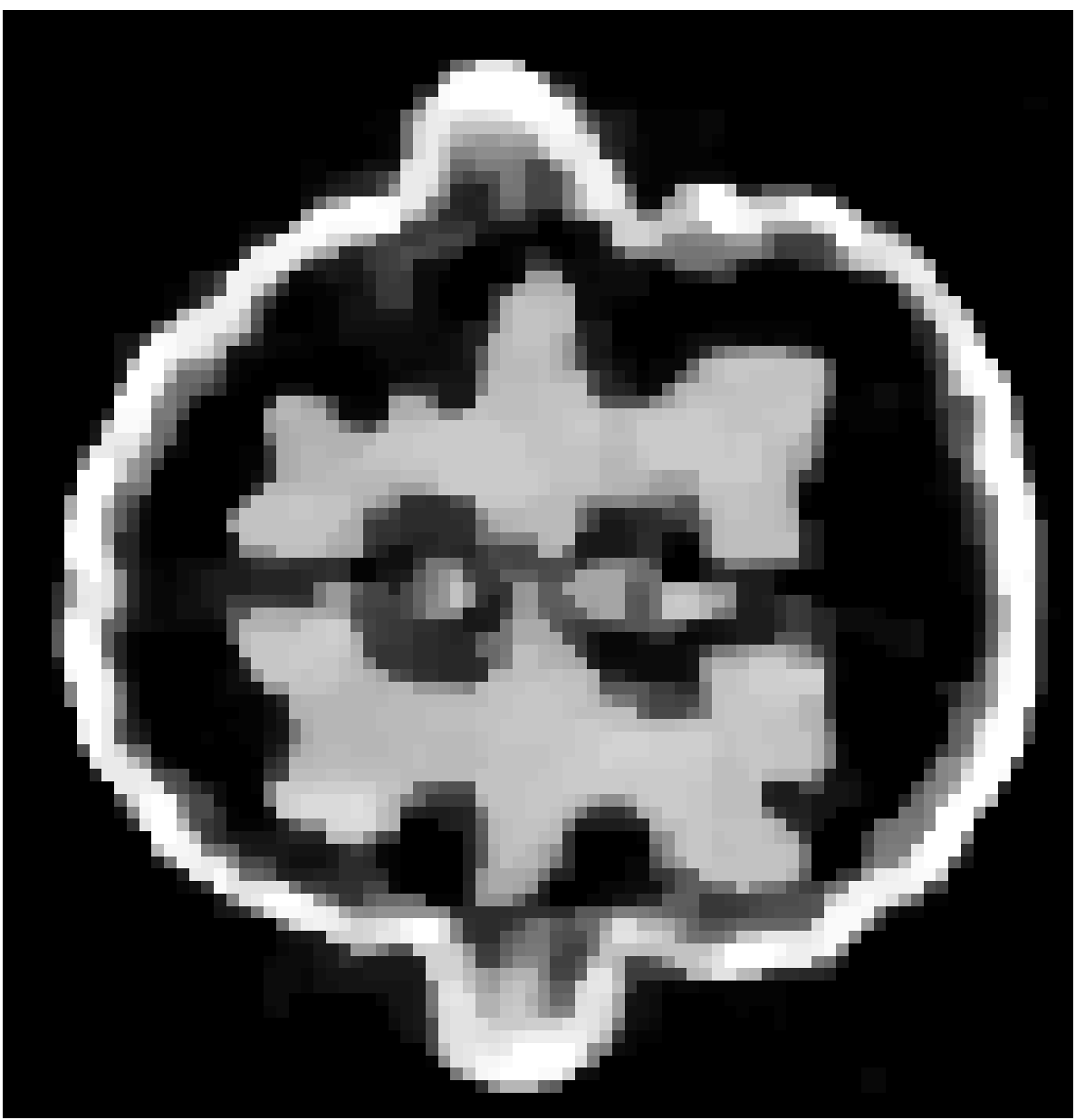}}
 \put(-85,65){$TV(\f)$}
 \put(-54,0){$10^{-4}$}
 \put(94,0){$10^{7}$}
 \put(0,0){$\Vert A \f - \greal \Vert$}
  \end{picture}
  \caption{\label{fig:LcurveNoise5_v1}L-curve  method applied to walnut
    data with 5\% additional noise. Left: L-curve plot. Right: Reconstruction
    computed using the selected $\alpha = 10.4$. The discretization level is $n=128$.}
  \end{figure}

As a reference, results with L-curve method were computed with both
test cases i) and ii). The L-curve plots and corresponding
reconstructions are here shown only for the discretization level
$n=128$, but the results for other discretization levels were similar. Figures \ref{fig:Lcurve} and
\ref{fig:LcurveNoise5_v1}, present the L-curve plots and the corresponding reconstructions
for the test cases i) and ii), respectively. 

\section{Discussion and conclusions}

\noindent
The results listed in Table \ref{tab:doublenorms} provide numerical evidence supporting the proposed parameter choice rule. It seems that the new multiresolution-based parameter choice rule finds useful and robust values for the Total Variation regularization parameter without the need for any {\em a priori} information. The fact that we use measured X-ray projection data (as opposed to simulated data) lends credibility to these initial findings. 

The precise mechanism behind the proposed method is still  unclear. We presented a mathematical proof that the reconstructions convergence for {\em any} choice of $\alpha$ as the resolution is refined. Therefore, it remains an open question whether the unstable norms in Tables \ref{tab:doublenorms} would converge with higher resolutions, or if those numbers are inaccurate, reflecting the numerical instability arising from the  ill-posedness of the underlying tomographic problem. In any case, in these real-data examples at least, the proposed method automatically adapts to the relevant  resolution and actual noise level. This is very promising indeed.

The S-curve method, introduced in
\cite{Kolehmainen2012,Hamalainen2013}, was found to work robustly and reliably for the present datasets. Curiously, the S-curve method produces parameter values of the same order of magnitude than the new resolution-based method. However, {\em a priori} information about the total variation of the unknown function is needed for applying the S-curve method, while the proposed multi-resolution based method does not need anything else than the data.

The proposed method has the drawback that many reconstructions need to be computed with various choices of resolution $n$ and regularization parameter $\alpha$. This can be computationally demanding. 

We stress that this work is just an initial feasibility study for the new parameter choice method. More comprehensive testing and comparisons to other methods are still needed to assess the properties of the new approach.

\appendix

\clearpage

\section{Construction of the difference matrices}\label{Appendix_simple}

\noindent
Let us illustrate the construction of difference matrices $D_{H}$ and $D_{V}$  discussed in Section \ref{sec:DiscTV} with a simple example. Consider the $3\times 3$ image matrix
$$
 \left[\begin{array}{ccc}
f_{11} & f_{12} & f_{13}\\
f_{21} & f_{22} & f_{23}\\
f_{31} & f_{32} & f_{33}\\
\end{array}\right],
$$
which is represented by the vector $\f=[f_{11},f_{21},f_{31},f_{12},f_{22},f_{32},f_{13},f_{23},f_{33}]^T\in\R^9$. Now the $9\times 9$ horizontal difference matrix $D_{H}$ takes  the form
$$
D_{H} =  \frac{1}{n}\left[\begin{array}{rrrrrrrrr}
-1 & 0 & 0 &  1 & 0 & 0 & 0 & 0 & 0\\
 0 & 0 & 0 & -1 & 0 & 0 & 1 & 0 & 0\\
 1 & 0 & 0 & 0 & 0 & 0 & -1 & 0 & 0\\
 0 & -1 & 0 & 0 & 1 & 0 & 0 & 0 & 0\\
 0 & 0 & 0 & 0 & -1 & 0 & 0 & 1 & 0\\
 0 & 1 & 0 & 0 & 0 & 0 & 0 & -1 & 0\\
 0 & 0 & -1 & 0 & 0 & 1 & 0 & 0 & 0 \\
 0 & 0 & 0 & 0 & 0 & -1 & 0 & 0 & 1 \\
 0 & 0 & 1 & 0 & 0 & 0 & 0 & 0 & -1 
\end{array}\right],
$$
The construction of the vertical difference matrix 
$D_{V}$ is similar.

\section{A note on $BV$ norms}\label{Appendix_example}

\noindent
What happened if one used the $BV$ norm 
$$
 \|u\|_{L^1(D)} + V_p(u) :=  \|u\|_{L^1(D)} +  \int_D\Big(\left|\frac{\partial u(x)}{\partial x_1}\right|^p+\left|\frac{\partial u(x)}{\partial x_2}\right|^p\Big)^{1/p}dm(x)
$$ 
with $p>1$ in the analysis of Section \ref{sec:convergence}? The norm (\ref{def:V}) used there corresponds to the choice $p=1$. First of all, taking $p>1$ would only give an inequality  in (\ref{apu limit 1}):
$$
  \lim_{j\to \infty}V_p(v_j)\geq V_p(u),
$$
and the proof would not work. 

\begin{figure}
\begin{picture}(340,80)
  \put(-20,12){\includegraphics[height=2cm]{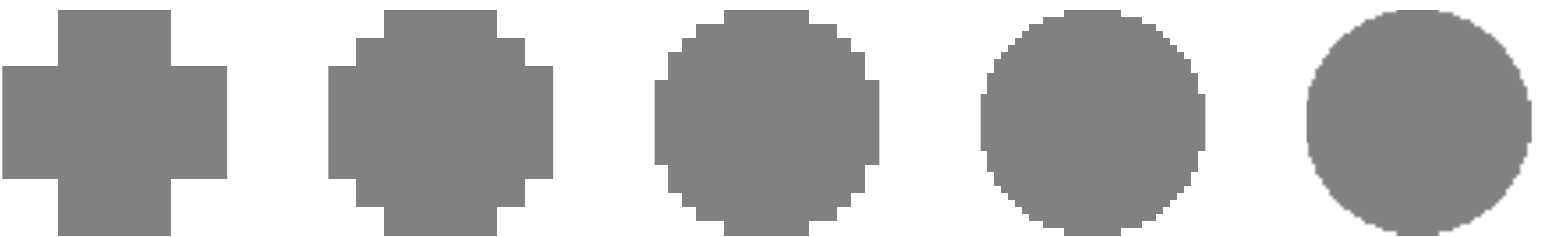}}
\end{picture}
 \caption{\label{fig:discdiscs}Pixel-based piecewise linear approximations to the characteristic function of the unit disc. Here white color denotes zero value and gray color denotes value one. The TV norm of all of these functions is exactly  8 (equal to the length of the boundary).}
 \end{figure}

Let us point out a specific example using $p=2$. Consider the characteristic function $\chi=\chi_{\{|x|<1\}}$ of the unit disc, whose BV norm is $2\pi$: length of the boundary  (which is $2\pi$) times the jump (which is $1$).  Now approximate $\chi$ using functions equal to one inside pixel-based set approximations to the unit disc, see Figure \ref{fig:discdiscs}. The total length of the boundary of each such set is $8$. Thus the BV norms of the approximations converge to $8$ instead of $2\pi$.

Now using the BV norm with $p=1$ leads to the (Manhattan) length of the unit circle to be $8$, ensuring the appropriate convergence.

\bibliographystyle{siam}
\bibliography{TVparamchoice_tod13}

\end{document}